\documentclass[10pt,reqno]{amsart}

\usepackage{amssymb, amscd, amsaddr, stmaryrd}
\usepackage{enumerate}
\usepackage{enumitem}
\setlist[itemize]{leftmargin=*}
\setlist[enumerate]{leftmargin=*}
\usepackage{cases}
\usepackage{amsfonts}
\usepackage[dvips,dvipdf]{graphicx}
\usepackage[usenames,dvipsnames]{color}

\usepackage{color}
\usepackage{color, colortbl}
\usepackage{multirow}
\usepackage[table]{xcolor}
\usepackage{comment}
\usepackage{subfigure}
\usepackage{epstopdf}
\usepackage{tikz}
\usepackage{chngcntr}

\usepackage{algorithm} 
\usepackage{algorithmic} 
\usepackage{pifont}

\usepackage{hyperref}
\usepackage{cleveref}

\hypersetup{
  bookmarksnumbered = true,
  bookmarksopen=false,
  pdfborder=0 0 0,         
  pdffitwindow=true,      
  pdfnewwindow=true, 
  colorlinks=true,           
  linkcolor=blue,            
  citecolor=magenta,    
  filecolor=magenta,     
  urlcolor=cyan              
}

\usepackage{geometry}
\usepackage{./macros}

\makeatletter
\newcommand{\definetitlefootnote}[1]{%
  \newcommand\addtitlefootnote{%
    \makebox[0pt][l]{$^{*}$}%
    \footnote{\protect\@titlefootnotetext}
  }%
  \newcommand\@titlefootnotetext{\spaceskip=\z@skip $^{*}$#1}%
}
\makeatother

\counterwithin*{equation}{section}



\Crefname{equation}{Eq.}{Eqs.}

\newtheorem{theorem}{Theorem}[section]
\newtheorem{alg}{Algorithm}[section]
\newtheorem{definition}{Definition}[section]
\newtheorem{proposition}[theorem]{Proposition}
\newtheorem{corollary}[theorem]{Corollary}

\newtheorem{lemma}[theorem]{Lemma}

\theoremstyle{definition}
\theoremstyle{definition}
\theoremstyle{remark}
\newtheorem{remark}[theorem]{Remark}
\newtheorem{example}[theorem]{Example}

\newtheorem{assumption}[theorem]{Assumption}

\begin{document}

\title[DLR-DG]{A semi-implicit dynamical low-rank discontinuous Galerkin method for space homogeneous kinetic equations. Part I: emission and absorption \addtitlefootnote}


\author[P. Yin, E. Endeve, C.D. Hauck, S.R. Schnake]{Peimeng Yin$^\dagger$, Eirik Endeve$^{\dagger}$$^{\ddagger}$, Cory D. Hauck$^{\dagger}$$^{\S}$, and Stefan R. Schnake$^{\dagger}$}

\address{$^\dagger$ Mathematics in Computation Section, Computer Science and Mathematics Division, \\ Oak Ridge National Laboratory, Oak Ridge, TN 37831, USA}
\address{$^\ddagger$ Department of Physics and Astronomy, \\  University of Tennessee, Knoxville, 1408 Circle Drive, Knoxville TN 37996, USA}
\address{$^\S$ Department of Mathematics, \\  University of Tennessee, Knoxville, 1403 Circle Drive, Knoxville TN 37996, USA}

\email{yinp@ornl.gov; endevee@ornl.gov; hauckc@ornl.gov; schnakesr@ornl.gov}



\definetitlefootnote{Notice:  This manuscript has been authored by UT-Battelle, LLC, under contract DE-AC05-00OR22725 with the US Department of Energy (DOE). The US government retains and the publisher, by accepting the article for publication, acknowledges that the US government retains a nonexclusive, paid-up, irrevocable, worldwide license to publish or reproduce the published form of this manuscript, or allow others to do so, for US government purposes. DOE will provide public access to these results of federally sponsored research in accordance with the DOE Public Access Plan (http://energy.gov/downloads/doe-public-access-plan).}



\date{\today}

\begin{abstract}
Dynamical low-rank approximation (DLRA) is an emerging tool for reducing computational costs and provides memory savings when solving high-dimensional problems.  
In this work, we propose and analyze a semi-implicit dynamical low-rank discontinuous Galerkin (DLR-DG) method for the space homogeneous kinetic equation with a relaxation operator, modeling the emission and absorption of particles by a background medium.  
Both DLRA and the DG scheme can be formulated as Galerkin equations.  
To ensure their consistency, a weighted DLRA is introduced so that the resulting DLR-DG solution is a solution to the fully discrete DG scheme in a subspace of the classical DG solution space.  
Similar to the classical DG method, we show that the proposed DLR-DG method is well-posed.  
We also identify conditions such that the DLR-DG solution converges to the equilibrium.
Numerical results are presented to demonstrate the theoretical findings.
\end{abstract}

\maketitle


\setcounter{tocdepth}{3}


\section{Introduction}

In this paper, we consider high-order approximation methods for solving kinetic equations using low-dimensional surrogates that capture their essential features.  
These methods have been demonstrated to be computationally cheaper for many high-dimensional dynamical systems (see, e.g., \cite{grasedyck2013}), and the dynamical low-rank approximation (DLRA) is one well-known method used for this purpose.  
Specifically, we design and analyze a dynamical low-rank discontinuous Galerkin (DLR-DG) method for solving the space homogeneous kinetic equation for modeling the emission and absorption of particles by a background medium.  

Kinetic models of particle systems consider the evolution of the particle distribution function $f(\boldsymbol{p},\boldsymbol{x},t)$, a phase-space density depending on the particle momentum $\boldsymbol{p}\in\mathbb{R}^{3}$, position $\boldsymbol{x}\in\mathbb{R}^{3}$, and time $t$.  
Kinetic equations, governing the evolution of $f$, are expressed as a balance between phase-space advection (e.g., due to inertia and external forces) and collisions (e.g., due to interparticle interactions or interactions with a background).  
In the absence of collisions, the distribution function can develop complex phase-space structures, while collisions tend to drive $f$ towards an equilibrium, characterized by (spatially) local conditions, in which the dynamics can be accurately described by fluid models (where variables depend only on $\boldsymbol{x}$ and $t$).  
As such, kinetic models are high-dimensional models that can exhibit low-dimensional structure under certain conditions (e.g., particle systems undergoing frequent collisions).  

DLRA methods can be traced back to the Dirac--Frenkel--McLachlan variational principle developed in the 1930s \cite{dirac1930,frenkel1934}.  
Essentially, the right-hand side of a matrix differential equation is projected onto the tangent space of the manifold of fixed rank matrices, which yields a set of differential equations that govern the factors of an SVD-like decomposition.  
As such, they can be suitable for modeling high-dimensional systems that exhibit dynamics in a lower-dimensional manifold (e.g., kinetic equations).  
Recently, they have been applied to simulate high-dimensional quantum systems, biological cellular systems \cite{jahnke2008, beck2000, lubich2005}, kinetic/transport equations \cite{einkemmer2018, peng2020, einkemmer2021efficient, einkemmer2021asymptotic, ding2021dynamical, peng2021high, peng2022sweep, einkemmer2022asymptotic}, hyperbolic problems with uncertainty \cite{kusch2022}, and neural network training \cite{schotthofer2022}.  

The discontinuous Galerkin (DG) method is a finite element method that uses a discontinuous piecewise polynomial space to approximate the numerical solution.  
The method offers several advantages, such as high-order accuracy on a compact stencil, compatibility with $hp$-adaptivity, and the ability to handle domains with complex geometry \cite{hesthaven2007, riviere2008, shu2009}.  
Its mathematical formulation makes it amenable to rigorous analysis.  
Moreover, DG methods are attractive for solving kinetic equations because of their ability to maintain structural properties (e.g., asymptotic limits \cite{larsenMorel_1989,adams_2001,guermondKanschat_2010} and conservation \cite{ayuso_etal_2011,cheng_etal_2013}) of the continuum model formulation, in part, because of flexibility in the approximation spaces.  
However, the use of the DG methods to solve kinetic equations in full dimensionality, without any form of adaptivity to reduce the total number of degrees of freedom, can be computationally expensive.  

The DLR-DG method studied in this paper applies DLRA to the matrix differential equation resulting from the semi-discretization of the kinetic equation using the DG method.  
The combination of DLRA and DG methods aims to leverage the benefits of both approaches by lowering the computational complexity relative to standard DG methods while retaining high-order accuracy.  
However, the direct application of DLRA can result in loss of fundamental solution properties, such as well-posedness and the ability to capture equilibria (steady-states), that are more easily captured with standard DG methods.  
In this work, we consider a model kinetic equation of relaxation-type in reduced dimensionality (by imposing axial symmetry in momentum space) and focus on establishing conditions for which the DLR-DG formulation possesses the same properties as the standard DG scheme.  
We use spherical-polar momentum space coordinates.  
As a result, a volume Jacobian appears in the inner product of the DG scheme, giving a matrix weight in the matrix differential equation. 

In the weak formulation, the DG scheme represents a Galerkin equation, for which the trial space and the test space are identical.  
Similarly, the DLRA is a minimization problem associated with a matrix differential equation, seeking a solution within a rank-$r$ manifold, $\mathcal{M}_r$, and can also be reformulated as a Galerkin equation.  
We aim to formulate the minimization problem for the DLRA such that its Galerkin equation is consistent with that of the semi-discrete DG scheme in its matrix formulation.  
Then, when the coefficient matrix of the DG solution possesses a rank-$r$ decomposition and evolves tangentially to $\mathcal{M}_{r}$, its formulation becomes identical to the Galerkin equation of the DLRA.  
This consistency establishes that the DLRA, under relatively mild conditions, inherits key properties of the DG scheme, including well-posedness.

In DLRA, the original $\mathbb{R}^{m\times n}$ matrix is approximated by a rank-$r$ SVD-like factorization into three matrices, including two bases in $\mathbb{R}^{m\times r}$ and $\mathbb{R}^{n\times r}$, and a square matrix in $\mathbb{R}^{r\times r}$, which are integrated separately.  
These matrices can be integrated in time using different methods, including the projector-splitting integrator \cite{lubich2014} and the unconventional integrator \cite{ceruti2022}. 
Both integrators can handle small singular values, but the unconventional integrator has the advantage of avoiding the unstable backwards in time integration substep of the projector-splitting method, making it suitable for dissipative problems \cite{ceruti2022}.  
Therefore, we use the unconventional integrator to integrate the DLRA derived in this study. 

For computational efficiency, implicit-explicit (IMEX) time discretization \cite{ascher1997implicit, pareschi2005implicit} is commonly used to integrate kinetic equations, combining implicit integration for collisions with explicit integration for phase-space advection.  
For collision operators, implicit time discretization is desired because the short time scales induced by collisions can render explicit methods inefficient.  
Many IMEX schemes, including diagonally implicit Runge--Kutta methods, can be decomposed into explicit updates followed by an implicit solve that is equivalent to a backward Euler update with a modified initial state and time step.  
This decomposition, combined with the fact that the collision term often incurs the highest computational cost, motivates our focus on the space homogeneous kinetic equation.   
Applying backward Euler time discretization to the unconventional integrator results in what we call the \emph{semi-implicit unconventional integrator} (SIUI).  

In the current paper, we demonstrate that the dynamical low-rank solution obtained with the SIUI at each time step is equivalent to the solution of a fully discrete DG scheme constructed from a subspace of the original DG solution space.  
This subspace is a function of the current state and therefore changes at each time step.
(We call the subspace the DLR-DG space and the corresponding DG solution the DLR-DG solution.)  
This result allows us to analyze the properties of the dynamical low-rank solution of the SIUI through the fully discrete DG scheme in the DLR-DG space.  
The well-posedness of the SIUI solution then follows from that of the DLR-DG solution.  

We also identify conditions for which the DLR-DG solution converges to the rank-$1$ equilibrium of our model equation.  
First, we solve the steady state equation to obtain the coefficient matrix of the equilibrium solution and factor it in $\mathcal{M}_r$.  
Next, we project the equilibrium solution to the current DLR-DG space and identify conditions to bound the projection error by the product of an arbitrarily small number and the $L^2$ error between the previous step DLR-DG solution and the equilibrium solution.  
Finally, based on the projection error of the equilibrium solution to the current DLR-DG space, we identify a sufficient condition on the time step so that the distance between the DLR-DG solution and the equilibrium solution decays geometrically.  
To the best of our knowledge, no prior works have shown the equilibrium convergence of the DRLA.

The rest of the paper is organized as follows.  
In Section \ref{sec-2}, we introduce the space homogeneous kinetic equation, the full-rank DG discretization, and summarize the properties of the full-rank DG solution.  
In Section \ref{sec-3}, we formulate the matrix differential equation associated with the DG scheme and introduce its weighted DLRA.  
In Section \ref{sect:fully-discrete-DLRA}, we introduce the SIUI and the equivalent DLR-DG scheme.  
By analysis of the DLR-DG scheme, we prove the well-posedness of the SIUI and the convergence of the DLR-DG solution to the equilibrium.  
Numerical examples illustrating the theoretical results are given in Section \ref{sec-5}.

\section{Background}\label{sec-2}

\subsection{Model Equations}

The space homogeneous kinetic equation modeling the emission and absorption of particles by a material background at rest can be written as (see, e.g., \cite{mihalasMihalas_1999})
\be\label{tp10}
\bal
    \partial_t f(x,\varepsilon,\vartheta,\varphi,t) = \mathcal C(f)(x,\varepsilon,\vartheta,\varphi,t),
\eal
\ee
where $f\geq 0$ is the phase-space distribution function depending on position $x\in D_{x}\subset\bbR^3$, and spherical-polar momentum coordinates $(\varepsilon,\vartheta,\varphi)$, and time $t\geq 0$.  
Here, $\varepsilon\geq 0$ is the particle energy, $\vartheta\in [0,\pi]$ is the latitudinal angle, and $\varphi\in[0,2\pi)$ the azimuthal angle.  
We also introduce the latitudinal angle cosine $\mu=\cos(\theta)\in [-1,1]$.  

Since we consider the space homogeneous case in this paper, we will suppress the explicit dependence on the position coordinate $x$ from hereon.  
Furthermore, we impose axial symmetry in the azimuthal direction in momentum space (i.e., $f$ is independent of $\varphi$).  
We then write Eq.~\eqref{tp10} as
\be\label{tp1}
\bal
    \partial_t f(\varepsilon,\mu,t) 
    = \mathcal C(f)(\varepsilon,\mu,t),
\eal
\ee
where the collision operator on the right-hand side is given by
\be\label{collisionterm}
    \mathcal C(f)(\varepsilon,\mu,t) = \eta(\varepsilon) - \chi(\varepsilon) f(\mu, \varepsilon,t),
\ee
where $\eta>0$ is the emissivity and $\chi>0$ is the opacity. 
Both the emissivity and opacity are assumed to be independent of the momentum space angle cosine $\mu$, as is often done when the particle-material coupling is modeled in the material rest frame \cite{mihalasMihalas_1999}.  
The specific dependence of the opacity $\chi$ on the particle energy $\varepsilon$ depends on the details of the particle-material interaction process. Generally, we make the following assumption on $\chi$.
 \begin{assumption}\label{remchi}  
There exists constants $\chi_{\min}, \chi_{\max}$ such that
\be
0<\chi_{\min}\leq \chi \leq \chi_{\max}.
\ee
\end{assumption}
The kinetic equation \eqref{tp1} is subject to the initial condition 
\be\label{tpinit}
    f(\mu,\varepsilon,t=0) = f_0(\mu,\varepsilon).
\ee
Then, \eqref{tp1} is well-posed if \Cref{remchi} is satisfied. 

There is no coupling across energies in the collision term on the right-hand side of Eq.~\eqref{tp1}, and $\varepsilon$ is simply a parameter of the model.  
However, we include the energy dimension in the DG discretization to develop a more general framework that can accommodate coupling in energy and angle --- either through the inclusion of inelastic scattering or external fields.  
In addition, discretizing in both energy and angle allows us to capture momentum space structures.

It can be verified that as as $t\rightarrow \infty$ the solution $f(\mu,\varepsilon,t)$ to \Cref{tp1} converges to the \emph{isotropic} equilibrium solution $f^\teq(\varepsilon)=\eta(\varepsilon)/\chi(\varepsilon)$, which is the solution of the steady equation
\be
    \mathcal{C}(f^\teq) = 0.
\ee
Since $f^{\teq}$ is low-dimensional (independent of $\mu$), it is expected that for problems with strong particle-material coupling, $f$ will tend to become independent of $\mu$ for large $t$.  

\subsection{DG discretization and matrix equations}

Given $\varepsilon_{\max} >0$, we denote the computational domain by $\Omega=\{(\mu,\varepsilon):\mu\in[-1,1],\varepsilon\in[0,\varepsilon_{\max}]\}$ 
with volume measure $\text{d}\Omega = \varepsilon^2\text{d}\varepsilon \text{d}\mu$\footnote{The Lebesgue measure for axially symmetric functions defined on a ball centered at 0 in $\mathbb{R}^3$ is $2\pi \varepsilon^2\text{d}\varepsilon \text{d}\mu$, but we drop the $2\pi$ as each integral will have it as a common factor.}.  
Let $L^2(\Omega)$ be the Hilbert space of square integrable functions defined on $\Omega$ with respect to the measure $\text{d}\Omega$, and inner product denoted by
\begin{equation}
\label{eq:inner-product-def}
(v,w;\varepsilon^2)_{\Omega} := \int_{\Omega} vw\, \text{d}\Omega =
 \int_0^{\varepsilon_{\max}} \int_{-1}^1 vw \,\varepsilon^2 d\mu d{\varepsilon}.
\end{equation}
The associated norm on $L^2(\Omega)$ is given by $\|\varepsilon w\|_{L^2(\Omega)}^2:=(w,w;\varepsilon^2)_\Omega$.  

We write $\Omega=\Omega_\mu \times \Omega_\varepsilon$, where $\Omega_\varepsilon =[0,\varepsilon_{\max}]$ and $\Omega_\mu = [-1,1]$ with measures $\varepsilon^2\text{d}\varepsilon$ and $\text{d}\mu$, respectively.  Let $(\cdot,\cdot;\varepsilon^2)_{\Omega_\varepsilon}$ and $(\cdot,\cdot)_{\Omega_\mu}$ be the $L^2$ inner products induced from the given measures. 

Given $N_\varepsilon\in\mathbb{N}$ and $N_\mu\in\mathbb{N}$, we partition $\Omega_\varepsilon$ and $\Omega_\mu$ into $N_\varepsilon$ and $N_\mu$ cells, respectively.  Denote these partitions by 
\begin{align}
0=\varepsilon_{1/2}<\varepsilon_{3/2}<&\ldots<\varepsilon_{N_\varepsilon-1/2}<\varepsilon_{N_\varepsilon+1/2}=\varepsilon_{\max}, \\
-1=\mu_{1/2}<\mu_{3/2}<&\ldots<\mu_{N_\mu-1/2}<\mu_{N_\mu+1/2}=1.
\end{align} 
We partition the domain $\Omega$ into logical rectangles given by
\begin{equation}
K_{\mathfrak{i} \mathfrak{j}} = \{(\mu, \varepsilon): \mu \in K^\mu_\mathfrak{i}, \ \varepsilon\in K^\varepsilon_\mathfrak{j}  \},
\end{equation}
where $K^\mu_\mathfrak{i}=[\mu_{\mathfrak{i}-1/2}, \mu_{\mathfrak{i}+1/2}]$ for $1\leq \mathfrak{i} \leq N_\mu$, and $K^\varepsilon_\mathfrak{j}=[\varepsilon_{\mathfrak{j}-1/2}, \varepsilon_{\mathfrak{j} +1/2}]$ for $1\leq \mathfrak{j} \leq N_\varepsilon$.  

We now define the discontinuous Galerkin finite element space in each direction as
\be\label{dgspace}
V_{z, h} := \{\phi \in L^2(\Omega_v): \phi|_{K^z_\mathfrak{i}} \in \cP_{k}(K^z_\mathfrak{i}), 1\leq \mathfrak{i} \leq N_z\},
\ee
where $z=\mu, \varepsilon$, and $\cP_{k}$ denotes polynomials of maximal degree ${k}$.
The discontinuous Galerkin finite element space is defined as
$$
V_h = V_{\mu,h}\otimes V_{\varepsilon,h} = \{v : v|_{K_{\mathfrak{i} \mathfrak{j}}} \in \cQ_k(K_{\mathfrak{i}\mathfrak{j}}), 1\leq \mathfrak{i} \leq N_\mu, 1\leq \mathfrak{j} \leq N_\varepsilon\},
$$
where $\cQ_k$ denotes the space of tensor-product polynomials of degree at most $k$ for each variable defined on $K_{\mathfrak{i} \mathfrak{j}}$.

Generally, for a scalar function $v$ and vector valued functions $\mathcal{V} = [v_1, \ldots, v_{m}]^\top\in\R^{m}$ and $\mathcal{W} = [w_1, \ldots, w_{n}]^\top\in\R^{n}$, defined on $D\subseteq\Omega$, we define
\be
\bal
(v,\mathcal{W};\phi)_{D} = & (\mathcal{W},v;\phi)_{D} = [ (v,w_j;\phi)_{D}]_{n\times 1} \in\R^{n},\\
(\mathcal{V},\mathcal{W}^\top ;\phi)_{D} = &
[(v_i,w_j;\phi)_{D}]_{m\times n}\in\R^{m\times n},
\eal
\ee
where $\phi = \phi(\varepsilon) > 0$ is a specified weighting function.

\subsubsection{Semi-discrete full-rank DG scheme}

The standard semi-discrete DG scheme, which we call the semi-discrete full-rank DG scheme for \Cref{tp1}, together with the initial data \eqref{tpinit}, is to find $f_{h}(\mu, \varepsilon, t)\in V_h$ such that
\begin{subequations}\label{semiDG}
\begin{align}
(\partial_t f_h, w_h;\varepsilon^2)_{\Omega} &= \mathcal{A}(f_h, w_h), \quad \forall w_h \in V_h \label{semiDGeq} \\
\label{semidginit}
(f_h|_{t=0},w_h;\varepsilon^2)_{\Omega} &= (f_0,w_h;\varepsilon^2)_{\Omega}, \quad \forall w_h \in V_h.
\end{align}
\end{subequations}
where $\mathcal{A}:L^2(\Omega)\times L^2(\Omega)\to\mathbb{R}$ is defined by
\be\label{bifs1}
\bal
\mathcal{A}(f_h, w_h)  = (\mathcal{C}(f_h),w_h;\varepsilon^2)_{\Omega} 
= \left(\eta, w_h;\varepsilon^2\right)_{\Omega}-\left(\chi f_h, w_h;\varepsilon^2 \right)_{\Omega}.
\eal
\ee

\begin{remark}
We use the term \textit{full-rank} throughout the paper to refer to a standard discontinuous Galerkin discretization with no low-rank techniques applied. 
\end{remark}
\begin{definition}
    The discrete equilibrium $f_h^{\teq}\in V_h$ is the solution to the variational problem 
\be\label{feql2proj}
    \mathcal{A}(f_h^{\teq}, w_h) = 0 \quad \forall w_h \in V_h.
\ee
\end{definition}
As long as \Cref{remchi} holds, Eq.~\eqref{feql2proj} admits a unique solution $f_h^{\teq}$, which can be shown to be a quasi-optimal approximation to $f^{\teq}$ in $L^2(\Omega)$. 

\subsubsection{Fully-discrete full-rank DG scheme}
We wish to employ implicit time discretization methods because the short time scales induced by collision operators can render explicit methods inefficient.  
For $\mathfrak{n}\geq 0$, let $f_h^\mathfrak{n} = f_h(\mu,\varepsilon,t^\mathfrak{n})\in V_h$ be an approximation of $f(\mu,\varepsilon,t^\mathfrak{n})$, where $t^\mathfrak{n}=\mathfrak{n}\Delta t$ and $\Delta t>0$ is a specified time step.
We apply a backward Euler time discretization to the semi-discrete full-rank DG scheme \eqref{semiDG}.
For simplicity, we denote 
\be
    D_t v^{\mathfrak{n}+1}
    = \frac{v^{\mathfrak{n}+1}-v^\mathfrak{n}}{\Delta t},
\ee
where $v$ can be any function (or matrix in the later sections).
Then, the first-order fully-discrete full-rank DG scheme for \Cref{tp1} is to find $f_h^{\mathfrak{n}+1} \in V_h$ such that
\begin{equation}
\label{fullDG1}
\begin{aligned}
    \left(D_t f_h^{\mathfrak{n}+1}, w_h;\varepsilon^2 \right)_{\Omega} &= \mathcal{A}(f_h^{\mathfrak{n}+1},w_h) \quad \forall w_h \in V_h.
\end{aligned}
\end{equation}

We now give the following result detailing the well-posedness, and convergence to the discrete equilibrium of the fully-discrete full-rank DG scheme.   For brevity, we omit the proof, since in \Cref{sect:fully-discrete-DLRA}, we prove a similar result in the low-rank setting.

\begin{proposition}\label{fullDG1lem}
For any $\Delta t>0$, there exists a unique solution $f^{\mathfrak{n}+1}_h$ of the fully-discrete, full-rank DG scheme \eqref{fullDG1} such that
\begin{itemize}
\item[(i)]\label{fullDG1lem-i} The solution $f^{\mathfrak{n}+1}_h$ is $L^2$ stable in the following sense:
\be\label{emabfullstab}
\bal
\|\varepsilon f_h^{\mathfrak{n}+1}\|_{L^2(\Omega)} \leq c^{\mathfrak{n}+1} \|\varepsilon f_{0}\|_{L^2(\Omega)} + \frac{1}{\chi_{\min}} (1-c^{\mathfrak{n}+1})\|\varepsilon\eta\|_{L^2(\Omega)},
\eal
\ee
where the parameter $c$ is given by
\be\label{ratior}
c = \frac{1}{1+\Delta t \chi_{\min}}.
\ee
\item[(ii)]\label{fullDG1lem-ii} The distance between $f_h^{\mathfrak{n}+1}$ and the discrete equilibrium $f^{\teq}_h$ is geometrically decreasing:
\be\label{Lstab}
\bal
\|\varepsilon (f_h^{\mathfrak{n}+1} - f_h^{\teq})\|_{L^2(\Omega)} \leq & c^{\mathfrak{n}+1}\|\varepsilon (f_h^{0} - f_h^\teq) \|_{L^2(\Omega)}.
\eal
\ee
where $f_h^\teq$ satisfies \Cref{feql2proj}.
\end{itemize}
\end{proposition}

\begin{remark}
From \eqref{Lstab}, it follows that $f_h^{\mathfrak{n}}$ converges to  
$f_h^\teq$ at a rate $O(\Delta t^{-\mathfrak{n}})$ as $\Delta t \to \infty$. 
\end{remark}

The main objective of this paper is to establish results analogous to \Cref{fullDG1lem} when the dynamical low-rank approximation is applied to the DG scheme. 
These are given in \Cref{sect:fully-discrete-DLRA}.

\section{Dynamical low-rank formulation}\label{sec-3}

In this section, we formulate low-rank approximations to \Cref{semiDG}. 

\subsection{Formulation of the matrix differential equation}

In order to apply the dynamical low-rank approximation, we first convert \Cref{semiDG} into an equivalent matrix differential equation via a basis expansion.  
Let $\{x_i(\mu)\}_{i=1}^m$ and $\{y_j(\varepsilon)\}_{j=1}^n$ be bases for the finite element spaces $V_{\mu,h}$ and $V_{\varepsilon,h}$, respectively.  
Here, $m=(k+1)N_\mu$ and $n=(k+1)N_\varepsilon$.  
We construct these bases using local Legendre polynomials on the local cells $K_\mathfrak{i}^{\mu}$ and $K_\mathfrak{j}^{\varepsilon}$ that are orthonormal with respect to the local inner products $L^{2}(K_\mathfrak{i}^{\mu})$ and $L^{2}(K_\mathfrak{j}^{\varepsilon})$, respectively.  
With this choice $\{x_i(\mu)\}_{i=1}^m$ forms on an orthonormal basis for $V_{\mu,h}$.  
However, $\{y_j(\varepsilon)\}_{j=1}^n$ does not form on an orthonormal basis for $V_{\varepsilon,h}$ due to the weight $\varepsilon^2$ in the inner product (cf.\eqref{eq:inner-product-def}).  
This fact has technical consequences for the remainder of the paper.

Given a function $w_h\in V_h$, its basis expansion can be written as
\be\label{DGtoMatDef}
w_h = \sum_{i=1}^m \sum_{j=1}^n W_{ij}(t) x_i(\mu) y_j(\varepsilon) = X^\top (\mu)W(t)Y(\varepsilon),
\ee
where $X:\Omega_\mu\to\R^m$ and $Y:\Omega_\varepsilon\to\R^n$ are defined by
\be\label{muepbasis}
X(\mu) = [x_1(\mu), \ldots, x_m(\mu)]^\top ~\text{ and }~ Y(\varepsilon)=[y_1(\varepsilon), \ldots, y_n(\varepsilon) ]^\top .
\ee
We call $W = [W_{ij}] \in \mathbb{R}^{m\times n}$ the \textit{coefficient matrix of} $w_h$ 
(\textit{with respect to the bases $\{x_i(\mu)\}_{i=1}^m$ and $\{y_j(\varepsilon)\}_{j=1}^n$)}.  For each fixed $i$, $W$ satisfies
\be
\sum_{j'=1}^n (y_j, y_{j'}; \varepsilon^2)_{\Omega_\varepsilon}W_{ij'} = (w_h, x_iy_{j}; \varepsilon^2)_{\Omega}, \quad j = 1, \ldots, n.
\ee

\begin{definition}\label{frobdef}
Given matrices $A,B \in \mathbb{R}^{m\times n}$ with entries $A_{ij}$ and $B_{ij}$, their {Frobenius inner product} is 
$
(A,B)_{\rm{F}}
=\text{tr}(A^\top B) 
=\sum_{i=1}^{m} \sum_{j=1}^{n} A_{ij}B_{ij}.
  $
The {Frobenius norm} of $A$ is 
$
\|A\|_{\rm{F}} = \sqrt{( A,A)_{\rm{F}}}.
$
\end{definition}

The next lemma relates weighted inner products of DG functions to weighted Frobenious inner products of the associated coefficient matrices. 
It follows from a direct calculation using (\ref{DGtoMatDef}).
\begin{lemma}\label{bilform}
Let $Z \in \mathbb{R}^{m\times n}$ and $W \in \mathbb{R}^{m\times n}$ be the coefficient matrices of $z_h\in V_h$ and $w_h \in V_h$, respectively, and let $\phi = \phi(\varepsilon)$ be a scalar function. 
Then
\be\label{L2Mat}
(\phi(\varepsilon)z_h,w_h;\varepsilon^{2})_{\Omega} = (I_m Z  A_{\phi}, W)_{\rm{F}} =(Z  A_{\phi}, W)_{\rm{F}},
\ee
where $I_m$ is the $m \times m$ identity matrix and the symmetric matrix
\be\label{dllp}
A_{\phi} = (\phi(\varepsilon)Y^\top (\varepsilon), Y(\varepsilon); \varepsilon^2)_{\Omega_\varepsilon}\in\R^{n\times n},
\ee
is block diagonal due to the locality of the basis.  
If further $\phi(\varepsilon)>0$, then $A_{\phi}$ is also positive-definite.
\end{lemma}

\begin{corollary}\label{biltomat}
Let $F \in \mathbb{R}^{m\times n}$ be the coefficient matrix of the DG solution $f_h \in V_h$ in \eqref{semiDG}, and $W \in \mathbb{R}^{m\times n}$ be the coefficient matrix of any function $w_h \in V_h$.
Then the semi-discrete DG scheme \eqref{semiDG} is equivalent to the following problem: Find $F(t)\in\mathbb{R}^{m\times n}$ such that 
\begin{subequations}\label{semiDGMat}
    \begin{align}
        &(\partial_t F(t) A_{\boldsymbol{1}}, W)_{\rm{F}} = \left( G(F), W \right)_{\rm{F}}, \quad \forall W\in\mathbb{R}^{m\times n},\\
        &F(0) = F_0,
    \end{align}
\end{subequations}
where $F_0 \in \mathbb{R}^{m\times n}$ is the coefficient matrix of $f_h(\mu,\varepsilon,0)$ obtained by solving \eqref{semidginit}. Here $A_{\boldsymbol{1}}$ is the symmetric positive-definite, block-diagonal matrix defined by \eqref{dllp} with $\phi = 1$,
and $G$ is the affine function defined by
\begin{equation} \label{Gformula1}
G(F) = L_{0} L_{\eta}^\top  - F A_{\EmAb}, \\
\end{equation}
where 
\be\label{lell}
L_{0} = (1, X)_{\Omega_\mu} \in \mathbb{R}^{m\times 1},
\quad 
L_{\eta}  =(\eta, Y;\varepsilon^2)_{\Omega_\varepsilon} \in \mathbb{R}^{n\times 1},
\ee
and the symmetric positive-definite, block-diagonal matrix $A_\chi$ is defined by \eqref{dllp} with $\phi=\chi$.
\end{corollary}
The variational problem \eqref{semiDGMat} immediately yields the following matrix-valued ODE:
\be\label{semidgmat+}
    \partial_t F = G(F)A_{\boldsymbol{1}}^{-1}.
\ee

\subsection{Weighted dynamical low-rank approximation}\label{subsec:weighted_dlra}

Let $\mathcal{M}_r \subset \mathbb{R}^{m\times n}$ be the manifold of rank-$r$ matrices ($r \leq \min\{m,n\}$).  
The Dynamical Low-Rank Approximation (DLRA) is traditionally formulated by evolving the matrix-valued ODE \eqref{semidgmat+} on $\mathcal{M}_r$ by a Galerkin projection of $\partial_t F$ onto the tangent space of $\mathcal{M}_r$ centered at $F$ (see e.g, \cite{lubichkoch2007DLRA}).  
This projection is on the space of $m\times n$ matrices and is traditionally orthogonal with respect to the standard Frobenius inner product in \Cref{frobdef}.  
However, such a formulation will not preserve the natural equivalence between the Galerkin equation of the DRLA and the matrix variational problem in \eqref{semiDGMat}.
In order to maintain this equivalence in the DLRA framework, we propose a modification to the standard DLRA approach that uses the weight $A_{\boldsymbol{1}}$ to characterize the tangent space.

\begin{definition}\label{Avar_inner_def}
For any $Z,W \in  \mathbb{R}^{m_1 \times n}$, $1\leq m_1 \leq m$, and any symmetric positive definite matrix $M \in  \mathbb{R}^{n\times n}$ with Cholesky factorization $M = C^\top C$, the $M$-weighted Frobenius inner product 
and its induced norm on $\mathbb{R}^{m_1 \times n}$ are given by 
\begin{equation}\label{Avardef}
    (Z,W)_{M} := (Z M,W)_{\rm{F}}=(ZC^\top ,WC^\top )_{\rm{F}} \quad\text{and}\quad
    \|W\|_{M}^2 := (W,W)_{M}.
\end{equation}
\end{definition}
\begin{remark}
The weighted Frobenius norm serves two purposes.  
The first is to introduce the matrix weight induced by the $\varepsilon^2$ integration weight in the definition of $\mathcal{A}$; see \eqref{bifs1}.  
The second is to introduce linear operations on the energy basis that, due to the transpose that appears in the rank-based representation of a matrix (e.g., the matrix $E^T$ in \eqref{FUSE} below), are often represented by left matrix multiplication. 
Thus for consistency, we reserve the usual vector norm on $\mathbb{R}^n$ for column vectors $x \in \mathbb{R}^{n \times 1}$ and use the Frobenius norm for row vectors $x^\top \in \mathbb{R}^{n \times 1}$, i.e., $\|x^\top\|_{\rm{F}}^{2} = tr(x x^\top) = \| x \|_{\rm{F}}^{2}$.  
\end{remark}

\begin{definition}\label{dlra_def}
Let $\hat{F}_0\in\mathcal{M}_r$ be given.  
The (weighted) dynamical low-rank approximation to \eqref{semidgmat+} is given by the solution $\hat{F}\in \mathcal{M}_r$ (where $\hat{F}$ approximates $F$) of the differential equation
\be\label{semidgam}
\partial_t \hat{F} = \argmin_{\delta \hat{F} \in \mathcal{T}_{\hat{F}} \mathcal{M}_r} J(\delta \hat{F}),
\quad \text{where} \quad 
J(\delta \hat{F})
= \|\delta \hat{F}-G(\hat{F})A_{\boldsymbol{1}}^{-1}\|_{A_{\boldsymbol{1}}},
\ee   
with initial condition $\hat{F}(0)=\hat{F}_0$.  
Here, $\mathcal{T}_{\hat{F}}\mathcal{M}_r$ is the tangent space of $\mathcal{M}_r$ at $\hat{F}$.
\end{definition}
\begin{remark}
    The initial condition $\hat{F}_0$ should be a rank-$r$ approximation to $F(0)$.  We delay the choice of $\hat{F}_0$ until the end of \Cref{subsec:weighted_dlra}.
\end{remark}

Like the usual DLRA \cite{lubichkoch2007DLRA}, \eqref{semidgam} can be rewritten into an equivalent system that updates the components of the low-rank decomposition of $\hat{F}$ in time; this equivalent system is often called the \textit{equations of motion}.  
Let $\hat{F}$ have the rank-$r$ decomposition
\begin{equation}\label{FUSE}
    \hat{F} = USE^\top, ~\text{ where }~U^\top U=E^\top A_{\boldsymbol{1}} E = I_r,
\end{equation}
with $U\in\mathbb{R}^{m\times r}$, $S\in\mathbb{R}^{r\times r}$, and $E\in\mathbb{R}^{n\times r}$ all full-rank matrices.%
\footnote{Unless otherwise stated, any matrices denoted with $U$ and $E$ satisfy $U^\top U = I_r$ and $E^\top A_{\boldsymbol{1}} E = I_r$, respectively.}
In terms of $U$, $S$, and $E$, the tangent space of $\mathcal{M}_r$ at $\hat{F}$ is (see e.g., \cite{lubichkoch2007DLRA}): 
\be\label{TFMr}
\mathcal{T}_{\hat{F}} \mathcal{M}_r = 
\{\delta USE^\top +U\delta S E^\top  + US\delta E^\top :  U^\top \delta U =0, \ E^\top A_{\boldsymbol{1}}\delta E=0 \},
\ee
where $\delta U \in \mathbb{R}^{m\times r}$, $\delta S \in \mathbb{R}^{r\times r}$, and $\delta E \in \mathbb{R}^{n\times r}$.
Due to the gauge conditions $U^\top \delta U = E^\top A_{\boldsymbol{1}} \delta E = 0$ in \eqref{TFMr},  any matrix $\delta \hat{F} \in \mathcal{T}_{\hat{F}} \mathcal{M}_r$ has the unique decomposition 
\be\label{deltaF}
\delta \hat{F} 
= \delta USE^\top +U\delta S E^\top  + US\delta E^\top
= P_U^\perp \delta \hat{F} A_{\boldsymbol{1}} P_E
    + P_U \delta \hat{F} A_{\boldsymbol{1}} P_E
    + P_U\delta \hat{F}^\top A_{\boldsymbol{1}} P^\perp_E,
\ee
where
\be
\bal
\delta U = P_U^\perp \delta \hat{F} A_{\boldsymbol{1}} ES^{-1},~~
\delta S =  U^\top  \delta \hat{F} A_{\boldsymbol{1}} E,~\text{ and }~\delta E =  P_E^\perp A_{\boldsymbol{1}} \delta \hat{F}^\top  U S^{-T}
\eal
\ee
with symmetric matrices
\be\label{projperp_u}
    P_U = UU^\top ~\text{ and }~P_U^\perp = I_m-P_U,
\ee
and
\be\label{projperp_e}
    P_E  = EE^\top ~\text{ and }~P_E^\perp=A_{\boldsymbol{1}}^{-1}-P_E .
\ee
The matrix $P_U$ is the orthogonal projection onto the column space of $U$ with respect to the standard inner product on $\mathbb{R}^m$ and $P_U^\perp$ is its orthogonal complement.  
The matrix $P_E A_{\boldsymbol{1}}$ is the orthogonal projection onto the column space of $E$ with respect to the inner product on $\mathbb{R}^n$ with weight $A_{\boldsymbol{1}}$.  
Moreover, for any $Z,W \in  \mathbb{R}^{\ell \times n}$, $1\leq \ell \leq m$, 
 \be
 \left(  Z A_{\boldsymbol{1}} P_E, W \right)_{ A_{\boldsymbol{1}}} = \left(Z, W A_{\boldsymbol{1}} P_E \right)_{A_{\boldsymbol{1}}}
 \ee
 where $(\cdot , \cdot)_{A_{\boldsymbol{1}}}$ is the Frobenius inner product defined in \Cref{Avar_inner_def}.

We now give several equivalent formulations of the weighted DLRA solution $\hat{F}$ in \Cref{dlra_def}.
\begin{proposition}\label{lrdiff}
The solution $\hat{F} = USE^\top  \in \mathcal{M}_r$ of \Cref{dlra_def}, with initial data $\hat{F}(0)= U^0S^0(E^0)^\top\in \mathcal{M}_r$ where $(U^0)^\top U^0=(E^0)^\top A_{\boldsymbol{1}}E^0 =I_r$,
satisfies the equivalent problems \cite{lubichkoch2007DLRA} 
\begin{enumerate}
\item[(i)] $\partial_t \hat{F} \in \mathcal{T}_{\hat{F}} \mathcal{M}_r$ is the solution of the Galerkin condition 
\begin{subequations}\label{Fglaerkin}
    \begin{align}
& \left ( \partial_t \hat{F}- G(\hat{F})A_{\boldsymbol{1}}^{-1}, \delta \hat{F}\right)_{A_{\boldsymbol{1}}}=0, \quad \forall 
 \delta \hat{F} \in \mathcal{T}_{\hat{F}} \mathcal{M}_r,\\
& \hat{F}(0)= U^0S^0(E^0)^\top.
    \end{align}
\end{subequations}
\item[(ii)]
The factors of $\hat F$ satisfy the equations of motion given by
\begin{subequations}\label{duse}
    \begin{align}
&\dot{U} = P_U^\perp G(\hat{F}) ES^{-1}, \quad \dot{S} = U^\top G(\hat{F})E, \quad \dot{E} = P_E^\perp G(\hat{F})^\top  US^{-T},\\
&U(0) = U^0, \quad S(0) =S^0, \quad E(0) = E^0,
    \end{align}
\end{subequations}
where $P_U^\perp$ and $P_E^\perp$ are defined in (\ref{projperp_u}) and (\ref{projperp_e}), respectively.
\item[(iii)] The matrices $\mathbf{K}=US \in \mathbb{R}^{m \times r}$, $\mathbf{L}=ES^\top \in \mathbb{R}^{n \times r} $, and $S$ satisfy the coupled ODE system
\begin{subequations}\label{semiDLRMat}
    \begin{align}
&\dot{\mathbf{K}} =  G(\mathbf{K}E^\top ) E, \quad
\dot{\mathbf{L}} = A_{\boldsymbol{1}}^{-1} G(U\mathbf{L}^\top )^\top  U,\quad
\dot{S} = U^\top G(USE^\top )E,\\
& \mathbf{K}(0) = U^0 S^0, \quad \mathbf{L}(0) = E^0 (S^0)^\top, \quad S(0) = S^0.
    \end{align}
\end{subequations}
\end{enumerate}
\end{proposition}
\begin{proof}We give a short sketch.
\begin{itemize}[leftmargin=*]
    \item \Cref{dlra_def} $\Leftrightarrow$ (i). The minimization problem \eqref{semidgam} is unchanged if $J$ is replaced by $\frac12 J^2$.  The minimization of this strongly convex quadratic functional over the linear subspace $\mathcal{T}_{\hat{F}}\mathcal{M}_r$ is equivalent to the Galerkin condition (i).
    \item (i) $\Rightarrow$ (ii). Since $\partial_t \hat{F} = \dot{U}SE^\top + U\dot{S}E^\top + US\dot{E}^\top$, the equations for $\dot U$, $\dot S$, and $\dot E$ in \eqref{duse} can be found from (\ref{Fglaerkin}a) by testing against
\be
\bal
\delta U SE^\top 
= P_U^\perp U_W S^{-\top}E^\top,
\quad
 U \delta S E^\top 
= U S_W E^\top, 
\quad
US (\delta E)^\top  
= US^{-\top}E_W^\top A_{\boldsymbol{1}} P_E^\perp,
\eal
\ee
respectively, where $U_W \in \mathbb{R}^{m\times r}$, $S_W \in \mathbb{R}^{r\times r}$, and $E_W \in \mathbb{R}^{n\times r}$ are arbitary.
By the arbitrariness of $U_W$, $S_W$, $E_W$, and the gauge condition $U^\top \dot{U} = E^\top A_{\boldsymbol{1}} \dot{E} = 0$, (\ref{Fglaerkin}a) reduces (\ref{duse}a).
\item (ii) $\Leftrightarrow$ (iii). Direct calculation: Take the derivative of $\mathbf{K}$ and $\mathbf{L}$ and use the product rule, \eqref{projperp_u}, and \eqref{projperp_e}.

\item (ii) $\Rightarrow$ (i). From the equations of motion \eqref{duse},  
\be
\label{eq:two-one}
\partial_t \hat{F} = \dot{U}SE^\top + U\dot{S}E^\top + US\dot{E}^\top = P_U^\perp G P_E
    + P_U G P_E
    + P_U GP^\perp_E.
\ee
Plugging \eqref{eq:two-one} and \eqref{deltaF} into (\ref{Fglaerkin}a), using \eqref{projperp_u} and \eqref{projperp_e}, verifies the result.
\end{itemize} 
\end{proof}

\begin{remark}
With the DLRA defined in \Cref{dlra_def}, the semi-discrete DG scheme in matrix formulation (\ref{semiDGMat}a) is identical to the Galerkin equation of the DRLA (\ref{Fglaerkin}a) when the coefficient matrix of the DG solution possesses a rank-$r$ decomposition and evolves tangentially to $\mathcal{M}_{r}$.
\end{remark}

\section{Fully discrete dynamical low-rank DG schemes}\label{sect:fully-discrete-DLRA}

In this section, we propose a fully discrete dynamical low-rank DG (DLR-DG) method. 
Similar to Proposition \ref{fullDG1lem} for the full-rank scheme, we investigate the well-posedness of the DLR-DG method and show the convergence of its solution to the equilibrium for a sufficiently large time step.

\subsection{The fully discrete DLR-DG schemes}

Applying a numerical integrator to the equations of motion in the form of Eq.~\eqref{duse} will produce an unstable method unless $\Delta t$ is of the same order as the smallest singular value of $S$ \cite{lubich2014projector}.  
Several DLRA temporal integrators have been developed with timestep restrictions that are much more reasonable \cite{lubich2014projector,ceruti2022,kieri2019projection}.  
Here we choose the unconventional integrator \cite{ceruti2022}, which is easily combined with the backward Euler method.

\subsubsection{A semi-implicit unconventional integrator}

The unconventional integrator of \cite{ceruti2022} can be viewed as an operator splitting method applied to the KLS system in Eq.~\eqref{semiDLRMat}, where the $K$ and $L$ equations are decoupled  and updated independently, followed by an update using the $S$ equation.  
We use backward (implicit) Euler for the underlying numerical integrator for all equations as collision operators generally induce timescales that cannot be efficiently advanced with an explicit method.  
Given $\Delta t>0$ and
the factored rank-$r$ matrix $\hat{F}^\mathfrak{n}= U^\mathfrak{n}S^\mathfrak{n}(E^\mathfrak{n})^\top $ with factors satisfying
\be\label{Fncond}
(U^\mathfrak{n})^\top U^\mathfrak{n}=I_{r}, \quad (E^\mathfrak{n})^\top A_{\boldsymbol{1}} E^\mathfrak{n} = I_{r},
\ee
one step of the method generates a new rank-$r$ matrix factorization 
\be\label{lrFnp1}
\hat{F}^{\mathfrak{n}+1}= U^{\mathfrak{n}+1}S^{\mathfrak{n}+1}(E^{\mathfrak{n}+1})^\top
\ee
with factors satisfying
\be\label{Fn1cond}
(U^{\mathfrak{n}+1})^\top U^{\mathfrak{n}+1}=I_{r}, \quad (E^{\mathfrak{n}+1})^\top A_{\boldsymbol{1}} E^{\mathfrak{n}+1} = I_{r}.
\ee
\Cref{lralg} precisely defines the semi-implicit unconventional integrator.
\begin{alg}\label{lralg}
A semi-implicit unconventional integrator (SIUI).
\begin{itemize}
\item{{Input:}} $U^\mathfrak{n}, S^{\mathfrak{n}}, E^{\mathfrak{n}}$, $\Delta t$; {output:} $U^{\mathfrak{n}+1},S^{\mathfrak{n}+1},E^{\mathfrak{n}+1}$.
    \item{\textbf{Step 1:}} Update $U^\mathfrak{n} \rightarrow U^{\mathfrak{n}+1}$ and $E^\mathfrak{n} \rightarrow E^{\mathfrak{n}+1}$ in parallel:
\begin{itemize}
    \item{\textbf{$K$-step}:} 
    \begin{itemize}
        \item Solve for $\mathbf{K}^{\mathfrak{n}+1}$ from the $m\times r$ matrix equation 
    \be\label{kstep}
    D_t \mathbf{K}^{\mathfrak{n}+1} 
    = G(\mathbf{K}^{\mathfrak{n}+1}(E^\mathfrak{n})^\top )E^\mathfrak{n}, \quad \mathbf{K}^\mathfrak{n}=U^\mathfrak{n}S^\mathfrak{n}.
    \ee
    \item Perform a QR factorization $\mathbf{K}^{\mathfrak{n}+1}=U^{\mathfrak{n}+1}R_{\mathbf{K}}$.
    \item Compute the $r \times r$ matrix $M^{\mathfrak{n}+1}=(U^{\mathfrak{n}+1})^\top  U^\mathfrak{n}$.
    \end{itemize}
    \item{\textbf{$L$-step}:} 
    \begin{itemize}
        \item Solve for $\mathbf{L}^{\mathfrak{n}+1}$ from the $n\times r$ matrix equation
    \be\label{lstep}
    D_t \mathbf{L}^{\mathfrak{n}+1}
    = A_{\boldsymbol{1}}^{-1}G(U^\mathfrak{n}(\mathbf{L}^{\mathfrak{n}+1})^\top )^\top U^\mathfrak{n}, \quad \mathbf{L}^\mathfrak{n}=E^\mathfrak{n}(S^\mathfrak{n})^\top.
    \ee
    \item Perform a generalized QR factorization (\Cref{wQR}) $\mathbf{L}^{\mathfrak{n}+1}= E^{\mathfrak{n}+1}R_{\mathbf{L}}$.
    \item  Compute the $r \times r$ matrix $N^{\mathfrak{n}+1}=(E^{\mathfrak{n}+1})^\top  A_{\boldsymbol{1}} E^\mathfrak{n}$.
    \end{itemize}
   
\end{itemize}
   \item{\textbf{Step 2:}} Update $S^\mathfrak{n}\rightarrow S^{\mathfrak{n}+1}$:
   \begin{itemize}
       \item{\textbf{$S$-step}:} 
       \begin{itemize}
           \item Project $S^\mathfrak{n}$ to the new bases
           \be\label{boldsn}
           S^{\mathfrak{n},*}=M^{\mathfrak{n}+1}S^\mathfrak{n}(N^{\mathfrak{n}+1})^\top.
           \ee
           \item Solve for $S^{\mathfrak{n}+1}$ from the $r\times r$ matrix equation
    \be\label{sstep}
    \frac{{S}^{\mathfrak{n}+1}-S^{\mathfrak{n},*}}{\Delta t} = (U^{\mathfrak{n}+1})^\top G(U^{\mathfrak{n}+1}S^{\mathfrak{n}+1} (E^{\mathfrak{n}+1})^\top )E^{\mathfrak{n}+1}.
    \ee
       \end{itemize}

   \end{itemize}
\end{itemize}
\end{alg}
\begin{remark}\label{remark_uc}~
The following remarks apply to \Cref{lralg}.
\begin{enumerate}
\item The choice of bases $U^{\mathfrak{n}+1}$ and $E^{\mathfrak{n}+1}$ used in the $S$-step is not unique.  
For any unitary matrices $V_U$, $V_E\in\mathbb{R}^{r\times r}$, 
the matrices $U^{\mathfrak{n}+1}V_U$ and $E^{\mathfrak{n}+1}V_E$ could replace $U^{\mathfrak{n}+1}$ and $E^{\mathfrak{n}+1}$, respectively, without changing $\hat{F}^{\mathfrak{n}+1}$.
\item The algorithm is semi-implicit since it uses explicit evaluation of the bases $U^\mathfrak{n}$ and $E^\mathfrak{n}$ in Eqs.~\eqref{kstep} and \eqref{lstep}, respectively, but makes implicit updates for $\mathbf{K}^{\mathfrak{n}+1}$, $\mathbf{L}^{\mathfrak{n}+1}$, and $S^{\mathfrak{n}+1}$.
\item 
$S^{\mathfrak{n},*}$ in \eqref{boldsn} is the projection of $S^\mathfrak{n}$ under the new bases $U^{\mathfrak{n}+1}$ and $E^{\mathfrak{n}+1}$.
Thus, $\|S^{\mathfrak{n},*}\|_{\rm{F}}\leq \|S^{\mathfrak{n}}\|_{\rm{F}}$.
For small $\Delta t$, the projection error is small \cite{ceruti2022}. 
For sufficiently large $\Delta t$, the projection error does not affect the SIUI solution's convergence to an equilibrium.
\item The matrix $R_{\mathbf{L}}$ in the $L$-step can also be computed by a regular QR factorization $\mathbf{L}^{\mathfrak{n}+1}= \tilde E^{\mathfrak{n}+1}\tilde R_{\mathbf{L}}$, followed by the weighted Gram--Schmidt decomposition $\tilde E^{\mathfrak{n}+1} = E^{\mathfrak{n}+1} \bar R_{\mathbf{L}}$, and then setting $R_{\mathbf{L}}=\bar R_{\mathbf{L}} \tilde R_{\mathbf{L}}$. \label{remark_ucd}
\item 
If $L_0$, defined in \eqref{lell}, is in the span of the columns of $U^\mathfrak{n}$, 
then $U^\teq = L_0/\|L_0\|=U^{\mathfrak{n}} z$ for some vector $z \in \mathbb{R}^{r \times 1}$. 
In this case, \eqref{kstep} reduces to 
\be\label{Kn1decomp0}
\bal
\mathbf{K}^{\mathfrak{n}+1}  
=  {U}^{\mathfrak{n}} \bar{R},
\eal
\quad
\text{where $\bar{R} = (S^\mathfrak{n}+\Delta t \|L_0\|zL_{\eta}^\top E^{\mathfrak{n}} )\left(I_r+\Delta t (E^\mathfrak{n})^\top A_\chi E^\mathfrak{n}\right)^{-1} \in \mathbb{R}^{r\times r}$}.
\ee
Thus,  \textbf{$K$-step} 
can be omitted, and we can set $U^{\mathfrak{n}+1}=U^\mathfrak{n}$.  
(See also \Cref{specialcase}, following \Cref{fsnproj+}.)
\end{enumerate}
\end{remark}

\subsubsection{DG formulation of the SIUI}

Given a low-rank approximation $\hat{f}_h^\mathfrak{n}$ with coefficient matrix $\hat{F}^{\mathfrak{n}} = U^\mathfrak{n}S^\mathfrak{n}E^\mathfrak{n}$, define the following subspaces of $V_h$ (which depend on $\hat{f}_h^\mathfrak{n}$): 
\begin{subequations}\label{subspaces}
\begin{align}
V_0^{\mathfrak{n}}=& \left\{v \ | \ v (\mu,\varepsilon)
=X^\top(\mu) U^{\mathfrak{n}}S(E^{\mathfrak{n}})^\top  Y(\varepsilon), \quad  \forall S \in \mathbb{R}^{r\times r} \right\},\\
V_1^\mathfrak{n}=& \left\{v \ | \ v (\mu,\varepsilon) 
=  X^\top(\mu) \mathbf{K}(E^\mathfrak{n})^\top  Y(\varepsilon), \quad   \forall \mathbf{K} \in \mathbb{R}^{m\times r} \right\}, \\
V_2^\mathfrak{n}= & \left\{v \ | \ v (\mu,\varepsilon)
=  X^\top(\mu) U^\mathfrak{n}\mathbf{L}^\top  Y(\varepsilon), \quad  \forall \mathbf{L} \in \mathbb{R}^{n\times r} \right\}.
\end{align}
\end{subequations}
It is easy to check that $\hat{f}_h^\mathfrak{n}= X^\top U^\mathfrak{n}S^\mathfrak{n}(E^\mathfrak{n})^\top Y\in V_0^{\mathfrak{n}} \cap V_1^{\mathfrak{n}}\cap V_2^{\mathfrak{n}}$, but $\hat{f}_h^\mathfrak{n} \not \in V_0^{\mathfrak{n}+1}$. 
However,  
\be\label{fbolds}
f_{S}^{\mathfrak{n},*} := X^\top U^{\mathfrak{n}+1} S^{\mathfrak{n},*}(E^{\mathfrak{n}+1})^\top Y \in V_0^{\mathfrak{n}+1},
\ee
where $S^{\mathfrak{n},*}$ is given in \eqref{sstep}. 
Moreover, $f_S^{\mathfrak{n},*}$ is the $L^2$ projection of $\hat f_h^\mathfrak{n}$ onto $V_0^{\mathfrak{n}+1}$:
\be\label{fboldseq}
(f_S^{\mathfrak{n},*}, w_h;\varepsilon^2)_{\Omega} = (\hat f_h^\mathfrak{n}, w_h;\varepsilon^2)_{\Omega}, \quad \forall w_h \in V_0^{\mathfrak{n}+1}.
\ee

The following lemma establishes an equivalent DG formulation for \eqref{kstep}-\eqref{sstep}.  

\begin{lemma}\label{dlrdg}
The matrices $\mathbf{K}^{\mathfrak{n}+1},\mathbf{L}^{\mathfrak{n}+1},S^{\mathfrak{n}+1}$ are solutions to \eqref{kstep}, \eqref{lstep}, \eqref{sstep}, respectively, iff $f_\mathbf{K}^{\mathfrak{n}+1}:=X^\top(\mu) \mathbf{K}^{\mathfrak{n}+1}(E^\mathfrak{n})^\top  Y(\varepsilon) \in V_1^{\mathfrak{n}}$, $f_\mathbf{L}^{\mathfrak{n}+1}:=X^\top(\mu) U^\mathfrak{n}(\mathbf{L}^{\mathfrak{n}+1})^\top  Y(\varepsilon)\in V_2^{\mathfrak{n}}$, and $f_S^{\mathfrak{n}+1} := X^\top(\mu) U^{\mathfrak{n}+1}S^{\mathfrak{n}+1}(E^{\mathfrak{n}+1})^\top  Y(\varepsilon) \in V_0^{\mathfrak{n}+1}$ solve the following DLR-DG scheme
\begin{subequations}\label{fullDG1K}
\begin{align}
    \left(D_t f_\mathbf{K}^{\mathfrak{n}+1}, w_1;\varepsilon^2 \right)_{\Omega} = & \mathcal{A}(f_\mathbf{K}^{\mathfrak{n}+1}, w_1), \quad \forall w_1 \in V_1^{\mathfrak{n}},
    \label{fullDG1K-K}\\
    \left(D_t f_\mathbf{L}^{\mathfrak{n}+1}, w_2;\varepsilon^2 \right)_{\Omega} = & \mathcal{A}(f_\mathbf{L}^{\mathfrak{n}+1}, w_2), \quad \forall w_2 \in V_2^{\mathfrak{n}},
    \label{fullDG1K-L}\\
    \left(D_t f_S^{\mathfrak{n}+1}, w_0;\varepsilon^2 \right)_{\Omega} = & \mathcal{A}(f_S^{\mathfrak{n}+1}, w_0), \quad \forall w_0 \in V_0^{\mathfrak{n}+1},
    \label{fullDG1K-S}
\end{align}
\end{subequations}
where 
$f_{\bK}^{\mathfrak{n}}=
f_\mathbf{L}^{\mathfrak{n}} = 
f_S^{\mathfrak{n}} = \hat f_h^\mathfrak{n}.
$
\end{lemma}
\begin{proof}
We only prove the equivalence between \eqref{sstep} and (\ref{fullDG1K}c); the others can be proved similarly.
Suppose $f_S^{\mathfrak{n}+1}$ solves (\ref{fullDG1K}c). 
Then, by \Cref{bilform}, $S^{\mathfrak{n}+1}$ solves
\be\label{halfmatodeK}
\bal
\left(U^{\mathfrak{n}+1}
D_t S^{n+1} (E^{\mathfrak{n}+1})^\top  A_{\boldsymbol{1}}, U^{\mathfrak{n}+1}W_0(E^{\mathfrak{n}+1})^\top  \right)_{\rm{F}}
=  \left(G(U^{\mathfrak{n}+1}S^{\mathfrak{n}+1} (E^{\mathfrak{n}+1})^\top ), U^{\mathfrak{n}+1}W_0(E^{\mathfrak{n}+1})^\top \right)_{\rm{F}},
\eal
\ee
for all $W_0\in \mathbb{R}^{r\times r}$.
The matrix form of \eqref{fboldseq}:
\be\label{lrprojmat}
\left(U^{\mathfrak{n}} {S}^\mathfrak{n} (E^{\mathfrak{n}})^\top  A_{\boldsymbol{1}}, U^{\mathfrak{n}+1}W_0(E^{\mathfrak{n}+1})^\top  \right)_{\rm{F}} = \left(U^{\mathfrak{n}+1} S^{\mathfrak{n},*}(E^{\mathfrak{n}+1})^\top  A_{\boldsymbol{1}}, U^{\mathfrak{n}+1}W_0(E^{\mathfrak{n}+1})^\top  \right)_{\rm{F}} 
\ee
can be used to replace $S^{\mathfrak{n}}$ by $S^{\mathfrak{n},*}$ in \eqref{halfmatodeK}. Then 
applying \Cref{FrobProp} and \eqref{Fn1cond} gives 
\be\label{halfmatodeK2}
\bal
\left(
\frac{S^{\mathfrak{n}+1}-S^{\mathfrak{n},*}}{\Delta t}, W_0 \right)_{\rm{F}}
= \left((U^{\mathfrak{n}+1})^\top G(U^{\mathfrak{n}+1}S^{\mathfrak{n}+1} (E^{\mathfrak{n}+1})^\top )E^{\mathfrak{n}+1}, W_0\right)_{\rm{F}}.
\eal
\ee
Since $W_0$ is arbitrary, (\ref{halfmatodeK2}) is equivalent to (\ref{sstep}).
\end{proof}

\subsection{Well-posedness}

We now obtain an analog of \Cref{fullDG1lem}~(i) for the SIUI listed in \Cref{lralg} -- namely that the DLR-DG scheme is uniquely solvable and uniformly stable.  

\begin{lemma}\label{KLSexistDG}
Given the low-rank representation $\hat{f}_{h}^\mathfrak{n}$, from which $f_\mathbf{K}^{\mathfrak{n}}$, $f_\mathbf{L}^{\mathfrak{n}}$, and $ f_S^{\mathfrak{n}}$ can be computed, the first order fully discrete DG scheme (\ref{fullDG1K}) admits a unique solution $(f_\mathbf{K}^{\mathfrak{n}+1}, f_\mathbf{L}^{\mathfrak{n}+1}, f_S^{\mathfrak{n}+1}) \in V_1^{\mathfrak{n}} \times  V_2^{\mathfrak{n}} \times V_0^{\mathfrak{n}+1}$ for any $\Delta t >0$. Equivalently, Algorithm \ref{lralg} admits a unique matrix solution $(\mathbf{K}^{\mathfrak{n}+1}, \mathbf{L}^{\mathfrak{n}+1}, S^{\mathfrak{n}+1})$.
\end{lemma}
\begin{proof}
We only prove the existence and uniqueness for $f_S^{\mathfrak{n}+1}$; the corresponding results for $f_\mathbf{K}^{\mathfrak{n}+1}$ and $f_\mathbf{L}^{\mathfrak{n}+1}$ can be proved in a similar way.  
Since (\ref{fullDG1K}c) is a linear system in a finite dimensional space where the domain and codomain have the same dimension, existence is equivalent to uniqueness. Let $\delta f_S^{\mathfrak{n}+1} \in V_0^{n+1}$ be the difference between two possible solutions to (\ref{fullDG1K}c).  
Then
\be
\left(\delta f_S^{\mathfrak{n}+1}, w_0;\varepsilon^2 \right)_{\Omega}= - \Delta t ( \chi(\varepsilon) \delta f_S^{\mathfrak{n}+1},w_0;\varepsilon^2)_{\Omega} \quad \forall w_0 \in V_0^{\mathfrak{n}+1}.
\ee
If $w_0=\delta f_S^{\mathfrak{n}+1}$, then
$
\|\varepsilon \delta f_S^{\mathfrak{n}+1}\|_{L^2(\Omega)}^2+\Delta t \|\varepsilon \sqrt{\chi(\varepsilon)}\delta f_S^{\mathfrak{n}+1} \|_{L^2(\Omega)}^2=0,
$
which implies $\delta f_S^{\mathfrak{n}+1}=0$. 
Therefore, the DG scheme (\ref{fullDG1K}c) admits a unique solution.
The uniqueness of $(f_\mathbf{K}^{\mathfrak{n}+1}, f_\mathbf{L}^{\mathfrak{n}+1}, f_S^{\mathfrak{n}+1})$ and the equivalence established by  Lemma \ref{dlrdg} imply that Algorithm \ref{lralg} admits the unique matrix solution $(\mathbf{K}^{\mathfrak{n}+1}, \mathbf{L}^{\mathfrak{n}+1}, S^{\mathfrak{n}+1})$.
\end{proof}

\begin{definition}
\label{lrdgn1}
We define the DG approximation $\hat{f}_h^{\mathfrak{n}+1} = f_S^{\mathfrak{n}+1}$ as the DLR-DG solution, and the subspace $V_0^{\mathfrak{n}+1}$ as the DLR-DG space.
\end{definition}

The $L^2$ stability of the DLR-DG solution $\hat{f}_h^{\mathfrak{n}+1}$ is established by the following lemma.
\begin{lemma}\label{lrdgbdd}
Suppose that $
\|\varepsilon \hat f_h^0\|_{L^2(\Omega)} 
\leq \|\varepsilon f_h^0\|_{L^2(\Omega)}.
$
Then the solution of the DG scheme \eqref{fullDG1K} is stable in the following sense
\be\label{emabfullstabK}
\|\varepsilon \hat f_h^{\mathfrak{n}+1}\|_{L^2(\Omega)} \leq c^{\mathfrak{n}+1} \|\varepsilon f_{0}\|_{L^2(\Omega)} + \frac{1}{\chi_{\min}} (1-c^{\mathfrak{n}+1})\|\varepsilon \eta\|_{L^2(\Omega)},
\ee
where $c$ is given in \eqref{ratior}.
\end{lemma}
\begin{proof}
Setting $w_0 = f_S^{\mathfrak{n}+1} \equiv \hat f_h^{\mathfrak{n}+1}$ (see \Cref{lrdgn1})  in \eqref{fullDG1K}c gives
\be\label{emab1ststab}
\bal
\left((1+\Delta t \chi) \hat f_h^{\mathfrak{n}+1}, \hat f_h^{\mathfrak{n}+1};\varepsilon^2 \right)_{\Omega} 
=  (\hat f_h^{\mathfrak{n}}, \hat f_h^{\mathfrak{n}+1};\varepsilon^2)_{\Omega}
+ \Delta t (\eta,\hat f_h^{\mathfrak{n}+1};\varepsilon^2)_{\Omega},
\eal
\ee
which, with Cauchy-Schwartz, leads to
\be\label{emab1ststab1}
\|\varepsilon \hat f_h^{\mathfrak{n}+1}\|_{L^2(\Omega)} \leq c \|\varepsilon \hat f_h^{\mathfrak{n}}\|_{L^2(\Omega)} + c\Delta t  \|\varepsilon \eta\|_{L^2(\Omega)},
\ee
where $c$ is given by \eqref{ratior}. 
Applying \eqref{emab1ststab1} recursively gives
\be\label{emab1ststab2}
\bal
\|\varepsilon \hat f_h^{\mathfrak{n}+1}\|_{L^2(\Omega)} \leq & c^{\mathfrak{n}+1} \|\varepsilon \hat f_h^{0}\|_{L^2(\Omega)} + \Delta t  \|\varepsilon \eta\|_{L^2(\Omega)}\sum_{i=1}^{\mathfrak{n}+1} c^i\\
\leq & c^{\mathfrak{n}+1} \|\varepsilon \hat f_h^{0}\|_{L^2(\Omega)} + \frac{1}{\chi_{\min}} (1-c^{\mathfrak{n}+1})\|\varepsilon \eta\|_{L^2(\Omega)}.
\eal
\ee
Thus the estimate \eqref{emabfullstabK} follows from \eqref{emab1ststab2}, the assumption on the initial data, and the fact that $\|\varepsilon  f_h^{0}\|_{L^2(\Omega)} 
\leq 
\|\varepsilon f_{0}\|_{L^2(\Omega)}$.
\end{proof}

\subsection{Convergence to the equilibrium distribution}

The convergence result in \Cref{fullDG1lem}~(ii) follows from the fact that the discrete equilibrium $f_h^\teq$ is in the trial space of the fully discrete full-rank DG scheme \eqref{fullDG1}.  
However, for the DLR-DG scheme, the space of trial functions may not contain the discrete equilibrium.
In this subsection, we provide additional conditions to ensure convergence of the DLR-DG solution $\hat{f}_h^\mathfrak{n}$ to $f_h^\teq$. 
We first evaluate the error between the equilibrium solution $f_h^\teq$ and its projection in the DLR-DG space and then investigate the convergence of $\hat{f}_h^\mathfrak{n}$ to this projection.

The equilibrium solution of the steady state equation \eqref{feql2proj} has the form $f_h^{\teq} = X^\top(\mu) F^{\teq} Y(\varepsilon)$, where $G(F^{\teq})=0$ and $G$ is given in \eqref{Gformula1}.  Equivalently, 
\be\label{Feqmat}
    F^{\teq}=L_{0}(L_\eta)^\top (A_\EmAb)^{-1}
    \quad 
    \text{and}
    \quad
    G = (F^{\teq} - F)(A_\EmAb)^{-1}
\ee
where the vectors $L_{0}$, $L_\eta$, and the matrix $A_\EmAb$ are given in \Cref{biltomat}.
The matrix $F^{\teq}$ in \eqref{Feqmat} is a rank-$1$ matrix that can be decomposed as
\be
\label{Feq}
F^{\teq}  = U^{\teq} S^\teq (E^\teq)^\top  ,  
\ee
where $U^{\teq} \in \mathbb{R}^{m \times 1}$, $E^\teq \in \mathbb{R}^{n\times 1}$, and $S^\teq \in \mathbb{R}^{1\times 1}$ are given by
\begin{equation}\label{USEeq}
U^{\teq} =  \frac{L_{0}}{\|L_{0}\|}, \quad
E^\teq =  \frac{(A_\EmAb)^{-1} L_\eta}{\|(L_\eta)^\top (A_\EmAb)^{-1}\|_{A_{\boldsymbol{1}}} } ,\quad
S^\teq = \|L_{0}\| \|(L_\eta)^\top (A_\EmAb)^{-1}\|_{A_{\boldsymbol{1}}}  .
\end{equation}
The vectors $U^{\teq}$ and $E^\teq$ satisfy the orthogonality conditions
$(U^{\teq})^\top U^{\teq}=1$ and ${(E^\teq)^\top  A_{\boldsymbol{1}} E^\teq=1}$, and the matrix $S^\teq$ satisfies the following estimate.
\begin{lemma}\label{seqbddlem}
The scalar $S^\teq$ is uniformly bounded in the following sense:  
\be\label{seqbdd}
|S^{\teq}| = \|\varepsilon  f_h^\teq\|_{L^2(\Omega)} \leq \chi_{\min}^{-1/2}\|\varepsilon \eta\|_{L^2(\Omega)}.
\ee
\end{lemma}
\begin{proof}
Setting $w_h = f_h^\teq$ in \eqref{feql2proj}, it is easy to show that
\be\label{feqbdd}
\chi_{\min}^{1/2} \|\varepsilon  f_h^\teq\|_{L^2(\Omega)} \leq \|\varepsilon \chi^{1/2} f_h^\teq\|_{L^2(\Omega)} \leq \|\varepsilon \eta\|_{L^2(\Omega)}.
\ee
A direct calculation using \Cref{FrobProp} gives
$
 \|\varepsilon  f_h^\teq\|_{L^2(\Omega)}^2 
=  \left( F^{\teq}, F^{\teq} \right)_{A_{\boldsymbol{1}}} 
= |S^\teq|^2
$
which, when substituted into \eqref{feqbdd}, recovers the estimate \eqref{seqbdd}.
\end{proof}

Let $f_{S}^{\teq,{\mathfrak{n}}}$ be the projection of $f_{h}^{\teq}$ onto the DLR-DG space $V_0^{\mathfrak{n}}$ that is orthogonal with respect to the inner product $(\cdot,\cdot;\varepsilon^2)_{\Omega}$, defined in \eqref{eq:inner-product-def}.  Then $f_{S}^{\teq,{\mathfrak{n}}}$ has an expansion of the form
\be
\label{feqsnproj}
f_{S}^{\teq,{\mathfrak{n}}} 
= X^\top(\mu)U^{\mathfrak{n}}{S}^{\teq,\mathfrak{n}}(E^{\mathfrak{n}})^\top  Y(\varepsilon)
 =  X^\top(\mu)P_{U^{\mathfrak{n}}} F^\teq A_{\boldsymbol{1}} P_{E^{\mathfrak{n}}}  Y(\varepsilon)
\in V_0^{\mathfrak{n}},
\ee
where the orthogonality condition $(f_{S}^{\teq,{\mathfrak{n}}}, w_h;\varepsilon^2 )_{\Omega} = (f_h^{\teq}, w_h;\varepsilon^2 )_{\Omega} \ \forall w_h \in V_0^{\mathfrak{n}}$ implies that  
$
{S}^{\teq,\mathfrak{n}} 
= 
(U^{\mathfrak{n}})^\top F^\teq A_{\boldsymbol{1}} E^{\mathfrak{n}}.
$

\subsubsection{Projection error of the equilibrium in the DLR-DG space}

The projection of $U^\teq$ onto the columns of $U^\mathfrak{n}$ is $P_{U^\mathfrak{n}}U^\teq:= U^\mathfrak{n}(U^\mathfrak{n})^\top U^\teq$, and the projection error is
\be\label{projectionErrorU}
\|U^\teq - P_{U^\mathfrak{n}}U^\teq\|^2 
= 1- \|P_{U^\mathfrak{n}}U^\teq\|^2= 1- \|(U^{\mathfrak{n}})^\top  U^{\teq}\|^2 \in [0,1].
\ee
Similarly, the (weighted) projection of $E^\teq$ onto the space spanned by the columns of $E^\mathfrak{n}$ is $P_{E^\mathfrak{n}}A_{\boldsymbol{1}} E^\teq := E^\mathfrak{n}(E^\mathfrak{n})^\top A_{\boldsymbol{1}} E^\teq$, and the (weighted) projection error is
\be\label{projectionErrorE}
\|(E^\teq)^\top -  (E^\teq)^\top A_{\boldsymbol{1}} P_{E^\mathfrak{n}} \|_{A_{\boldsymbol{1}}}^2
= 1- \|(E^\teq)^\top A_{\boldsymbol{1}} P_{E^\mathfrak{n}}\|_{A_{\boldsymbol{1}}}^2
= 1- \|(E^{\mathfrak{n}})^\top A_{\boldsymbol{1}}E^\teq\|^2\in [0,1].
\ee

Lemma~\ref{fsnproj} and Lemma~\ref{fsnprojU} below provide upper bounds for the projection errors in \eqref{projectionErrorE} and \eqref{projectionErrorU}, respectively. 
Their proofs can be found in \Cref{pfsnproj} and \Cref{pfsnprojU}, respectively.

\begin{lemma}\label{fsnproj} 
Assume that for some constant $\beta \in (0, 1]$,
\be\label{UUEEasumpt}
\|P_{U^\mathfrak{n}}U^\teq\| \geq \beta.  
\ee
Then, for any $\delta > 0$ and any $\Delta t \geq \Delta t_1=\frac{\sqrt{r}}{\beta \delta \chi_{\min}}$,
\be\label{basebdd}
\bal
1-\|(E^\teq)^\top A_{\boldsymbol{1}} P_{E^{\mathfrak{n}+1}}\|_{A_{\boldsymbol{1}}}^2
\leq & \frac{\delta^2}{|S^\teq|^2}\|\varepsilon(\hat{f}_h^{\mathfrak{n}}- f^\teq_h) \|_{L^2(\Omega)}^2.
\eal
\ee
Moreover, if $\hat{f}_h^\mathfrak{n}=f_h^\teq$, then for any $\Delta t>0$, 
\be\label{UUEEasumpt1}
\|(E^\teq)^\top A_{\boldsymbol{1}} P_{E^{\mathfrak{n}+1}}\|_{A_{\boldsymbol{1}}}
=1.
\ee
\end{lemma}

Define the symmetric matrix
\begin{equation}
\label{PE-chi}
    P_{E}^\chi = E(E^TA_\chi E)^{-1}E^T.
\end{equation}
Then $P_E^\chi A_{\chi}$ is the orthogonal projection onto the column space of $E$ with respect to the inner product on $\mathbb{R}^n$ with weight $A_{\chi}$.

\begin{lemma}\label{fsnprojU} 
Assume there exists a constant $\alpha>0$ such that 
\be\label{UUEEasumptU}
\| (E^\teq)^\top A_\EmAb P_{E^{\mathfrak{n}}}^\chi \|_{A_{\boldsymbol{1}}} \geq\alpha.
\ee

Then for any $\delta > 0$ and any $\Delta t \geq \Delta t_2=\frac{{r}^{1/2}\chi_{\max}^{1/2}}{\alpha \delta \chi_{\min}^{3/2}}$, 
\be\label{basebddU}
\bal
1-\|P_{U^{\mathfrak{n}+1}}U^\teq\| 
\leq & \frac{\delta^2}{|S^\teq|^2} \|\varepsilon(\hat{f}_h^{\mathfrak{n}}- f^\teq_h) \|_{L^2(\Omega)}^2. 
\eal
\ee
Moreover, if $\hat{f}_h^\mathfrak{n}=f_h^\teq$, then for any $\Delta t>0$, 
\be\label{UUEEasumpt1U}
\|P_{U^{\mathfrak{n}+1}}U^\teq\|
=1.
\ee
\end{lemma}

\Cref{fsnproj} and \Cref{fsnprojU} can be used to bound the projection error of the equilibrium with respect to $U^{\mathfrak{n}+1}$ and $E^{\mathfrak{n}+1}$. 

\begin{lemma}\label{fsnproj+} 
Assume there exist constants $\beta \in (0,1]$ and $\alpha>0$ such that $\|P_{U^\mathfrak{n}}U^\teq\|\geq \beta$ and $\| (E^\teq)^\top A_\EmAb P_{E^{\mathfrak{n}}}^\chi \|_{A_{\boldsymbol{1}}} \geq \alpha$.  
Then for any $\delta > 0$, 
there exists
\be\label{tstari}
    \Delta t_0=\frac{\sqrt{2}}{\delta}\max\left\{\frac{{r}^{1/2}}{\beta \chi_{\min}} , \frac{{r}^{1/2}\chi_{\max}^{1/2}}{\alpha \chi_{\min}^{3/2}}\right\}
\ee
such that when $\Delta t \geq \Delta t_0$, 
\be\label{Eerror}
    \|\varepsilon( f_{S}^{\teq, \mathfrak{n}+1}-f_h^{\teq})\|_{L^2(\Omega)}  \leq  \delta \|\varepsilon(\hat{f}_h^{\mathfrak{n}}- f^\teq_h) \|_{L^2(\Omega)}.  
\ee
Moreover, if $\hat{f}_h^{\mathfrak{n}} = f^\teq_h$, it follows that for any $\Delta t>0$,
\be\label{projeqS}
f_{S}^{\teq, \mathfrak{n}+1} = f^\teq_h.
\ee
\end{lemma}
\begin{proof}
By (\ref{Feq}), \Cref{bilform}, and Lemma \ref{FrobProp}, 
\be\label{htosnproj}
\bal
& \|\varepsilon( f_h^{\teq}-f_{S}^{\teq, \mathfrak{n}+1})\|_{L^2(\Omega)}^2 =  \| F^{\teq}  - P_{U^{\mathfrak{n}+1}} F^\teq A_{\boldsymbol{1}} P_{E^{\mathfrak{n}+1}} \|_{A_{\boldsymbol{1}}}^2\\
& = |S^\teq|^2\left[\, \left(1-\|P_{U^{\mathfrak{n}+1}}U^\teq\|^2\right) + \|P_{U^{\mathfrak{n}+1}}U^\teq\|^2 \left( 1 - \|(E^\teq)^\top A_{\boldsymbol{1}} P_{E^{\mathfrak{n}+1}}\|_{A_{\boldsymbol{1}}}^2\right) \,\right].
\eal
\ee
By Lemma \ref{fsnproj} and Lemma \ref{fsnprojU} (with $\delta$ being replaced by $\delta/\sqrt{2}$), it follows that when $\Delta t \geq \Delta t_0$, where $\Delta t_0$ is given by Eq.~\eqref{tstari},
the following estimates hold
\begin{subequations}\label{basebdd0}
\begin{align}
1-\|P_{U^{\mathfrak{n}+1}}U^\teq\|^2 &\leq  \frac{\delta^2}{2|S^\teq|^2} \|\varepsilon(\hat{f}_h^{\mathfrak{n}}- f^\teq_h) \|_{L^2(\Omega)}^2, \\ 
1-\|(E^\teq)^\top A_{\boldsymbol{1}} P_{E^{\mathfrak{n}+1}}\|_{A_{\boldsymbol{1}}}^2 &\leq  \frac{\delta^2}{2|S^\teq|^2} \|\varepsilon(\hat{f}_h^{\mathfrak{n}}- f^\teq_h) \|_{L^2(\Omega)}^2,
\end{align}
\end{subequations}
which, when substituted into (\ref{htosnproj}), yields \eqref{Eerror}. Then (\ref{projeqS}) follows from (\ref{htosnproj}), using (\ref{UUEEasumpt1}) and (\ref{UUEEasumpt1U}).
\end{proof}

\begin{remark}\label{specialcase} 
If $\|P_{U^\mathfrak{n}}U^\teq\|=\|(U^{\mathfrak{n}})^\top  U^{\teq}\|=1$, then $U^{\eq} = U^{\mathfrak{n}} z$ for some vector $z \in \mathbb{R}^{r \times 1}$, and
\eqref{Kn1decomp0} implies that the \textbf{$K$-step} can be omitted for \Cref{lralg} by simply taking $U^{\mathfrak{n}+1}=U^\mathfrak{n}$. Then for any $\delta>0$, there exists $\Delta t_0 = \frac{{r}^{1/2}}{\delta \chi_{\min}}$, such that when $\Delta t \geq \Delta t_0$, \eqref{Eerror} holds.
\end{remark}

\subsubsection{Convergence of the DLR-DG solution to the equilibrium}

We estimate the convergence of the DLR-DG solution $\hat{f}_h^{\mathfrak{n}+1}$ to the equilibrium $f_h^\teq$.
We first provide a one-step estimate.

\begin{theorem}\label{lem:onestep}
Suppose the assumptions in \Cref{fsnproj+} hold. For any $\delta>0$, let $\Delta t_0$ be given in \eqref{tstari}.  Then for any $\Delta t \geq \Delta t_0$, 
\be\label{LstabS++}
\|\varepsilon (\hat f_h^{\mathfrak{n}+1} - f_h^{\teq})\|_{L^2(\Omega)} \leq (c+\delta_\chi) \|\varepsilon (\hat f^{\mathfrak{n}}_h-f_h^{\teq})\|_{L^2(\Omega)},
\ee
where $c$ is given by \eqref{ratior} and $\delta_\chi = \left(1+\frac{\chi_{\max}}{\chi_{\min}}\right) \delta$. Moreover, if $\hat{f}_h^{\mathfrak{n}} = f^\teq_h$, then for any $\Delta t>0$, 
\be\label{DLReq}
\hat f_h^{\mathfrak{n}+1} = f^\teq_h.
\ee

\end{theorem}
\begin{proof}
Since $\eta(\varepsilon) = \chi(\varepsilon) f^{\teq}(\varepsilon)$, the DG scheme (\ref{fullDG1K}c) can be written using (\ref{fboldseq}) and (\ref{feql2proj}) as
\be\label{fullDG1r2S}
\left((1+\Delta t \chi) f_S^{\mathfrak{n}+1}, w_h; \varepsilon^2\right)_{\Omega} =\left (\Delta t\chi  f^{\teq}_h + f_S^{\mathfrak{n},*},w_h;\varepsilon^2\right)_{\Omega}
\quad \forall w_h \in V_0^{\mathfrak{n}+1}.
\ee
Subtracting $((1+\Delta t \chi) f_{S}^{\teq, \mathfrak{n}+1}, w_h;\varepsilon^2)_{\Omega}$ from (\ref{fullDG1r2S}) yields
\be\label{fullDG1r3S}
\left((1+\Delta t \chi)  (f_S^{\mathfrak{n}+1} - f_{S}^{\teq, \mathfrak{n}+1}), w_h ;\varepsilon^2\right)_{\Omega} =\left (\Delta t\chi (f^{\teq}_h- f_{S}^{\teq, \mathfrak{n}+1}) + (f_S^{\mathfrak{n},*}-f_{S}^{\teq, \mathfrak{n}+1}),w_h;\varepsilon^2\right)_{\Omega}.
\ee
Setting $w_h = f_S^{\mathfrak{n}+1} - f_{S}^{\teq, \mathfrak{n}+1} \in V_0^{\mathfrak{n}+1}$ in (\ref{fullDG1r3S}) and applying the Cauchy–Schwarz inequality gives
\be\label{fullDG1r4S}
\|\varepsilon (f_S^{\mathfrak{n}+1} - f_{S}^{\teq, \mathfrak{n}+1})\|_{L^2(\Omega)} \leq c \|\varepsilon (f_S^{\mathfrak{n},*}- f_{S}^{\teq, \mathfrak{n}+1})\|_{L^2(\Omega)} +c\Delta t  \|\varepsilon\chi( f_{S}^{\teq, \mathfrak{n}+1}- f^{\teq}_h)\|_{L^2(\Omega)},
\ee
where $c$ is given by (\ref{ratior}).
By the triangle inequality and (\ref{fullDG1r4S}),
\be\label{LstabS0} 
\bal
\|\varepsilon (f_S^{\mathfrak{n}+1} - f_h^{\teq})\|_{L^2(\Omega)} 
&\leq \|\varepsilon (f_S^{\mathfrak{n}+1} - f_{S}^{\teq, \mathfrak{n}+1})\|_{L^2(\Omega)} + \|\varepsilon (f_{S}^{\teq, \mathfrak{n}+1}-f_h^{\teq})\|_{L^2(\Omega)} \\
& \leq  c \|\varepsilon (f_S^{\mathfrak{n},*}- f_{S}^{\teq, \mathfrak{n}+1})\|_{L^2(\Omega)}  +(1+c\Delta t \chi_{\max}) \|\varepsilon( f_{S}^{\teq, \mathfrak{n}+1}-f_h^{\teq})\|_{L^2(\Omega)}.
\eal
\ee
Additionally, $f_S^{\mathfrak{n},*}$ and $f_{S}^{\teq, \mathfrak{n}+1}$ are both $L^2$ projections of $\hat{f}_h^\mathfrak{n}$ and $f_h^\teq$ onto $V_0^{\mathfrak{n}+1}$; therefore
\be\label{ctmapping}
    \|\varepsilon (f_S^{\mathfrak{n},*}- f_{S}^{\teq, \mathfrak{n}+1})\|_{L^2(\Omega)} \leq \|\varepsilon (\hat{f}_h^{\mathfrak{n}} - f_h^{\teq})\|_{L^2(\Omega)}
\ee
By \eqref{LstabS0}, \eqref{ctmapping},  and \Cref{lrdgn1}, the stability estimate follows:
\be\label{LstabS+}
\bal
\|\varepsilon (\hat f_h^{\mathfrak{n}+1} - f_h^{\teq})\|_{L^2(\Omega)} \leq & c \|\varepsilon (\hat f_h^{\mathfrak{n}} - f_h^{\teq})\|_{L^2(\Omega)}  + \left(1+\frac{\chi_{\max}}{\chi_{\min}}\right) \|\varepsilon( f_{S}^{\teq, \mathfrak{n}+1} - f_h^{\teq})\|_{L^2(\Omega)}.
\eal
\ee
If $\Delta t \geq \Delta t_0$, then  \eqref{Eerror} holds and, when substituted into \eqref{LstabS+}, gives \eqref{LstabS++}. 
If $\hat{f}_h^{\mathfrak{n}} = f^\teq_h$, then \eqref{projeqS} 
and \eqref{LstabS+} imply  \eqref{DLReq}.  
\end{proof}

Unlike the full-rank case, the one-step estimate in \Cref{lem:onestep} cannot be trivially extended to a multi-step estimate.  
This is because of the disconnect between the conclusion of \Cref{fsnproj} and the hypothesis of \Cref{fsnprojU}, which bound the projection with respect to the $A_{\boldsymbol{1}}$- and $A_\chi$-inner products respectively.  
In order to bootstrap the one-step estimate further, we require the following lemma which controls $\| (E^\teq)^\top A_\EmAb P_{E^{\mathfrak{n}}}^\chi \|_{A_{\boldsymbol{1}}}$ by $\|(E^\teq)^\top A_{\boldsymbol{1}} P_{E^{\mathfrak{n}}}\|_{A_{\boldsymbol{1}}}$, where $P_{E^{\mathfrak{n}}}^\chi$ is defined in \eqref{PE-chi}.
This estimate depends on $\frac{\chi_{\max}}{\chi_{\min}}$, the weighted condition number of $A_\chi$. 
\begin{lemma}\label{normrel}
For any $\alpha \in (0,1)$, there exists $\gamma^*\in (\alpha,1)$, dependent only on $\frac{\chi_{\max}}{\chi_{\min}}$ and $\alpha$, such that if $E\in\mathbb{R}^{n\times r}$ with $E^TA_{\boldsymbol{1}}E=I_r$ and $\|(E^\teq)^\top A_{\boldsymbol{1}} P_{E}\|_{A_{\boldsymbol{1}}} \geq \gamma^*$, then
$\| (E^\teq)^\top A_\EmAb P_{E}^\chi \|_{A_{\boldsymbol{1}}} \geq \alpha$. 
\end{lemma}
\begin{proof}
Decompose $E^\teq$ as
\be\label{Eqdecomp}
E^\teq = E_1^\teq + E_2^\teq,
\ee
where $E_1^\teq=P_{E}A_{\boldsymbol{1}} E^\teq$ is the orthogonal projection of $E^\teq$ onto the column space of $E$ and $E_2^\teq= P_{E}^\perp A_{\boldsymbol{1}} E^\teq$ is the orthogonal complement
satisfying $ \|(E_1^\teq)^\top\|_{A_{\boldsymbol{1}}}^2 +\|(E_2^\teq)^\top \|_{A_{\boldsymbol{1}}}^2=1$.
Since $P_E^\chi A_\chi$ is also a projection onto the column space of $E$,  
\be\label{chiproj}
P_E^\chi A_\chi E_1^\teq  = E_1^\teq.
\ee
Suppose $\|(E_1^\teq)^\top\|_{A_{\boldsymbol{1}}} =: \gamma\in \left(\alpha,1\right]$. By \Cref{projachieig}, 
\be\label{Approj}
\| (E_2^\teq)^\top A_\EmAb P_{E}^\chi \|_{A_{\boldsymbol{1}}}^2 \leq \frac{\chi_{\max}}{\chi_{\min}} \| (E_2^\teq)^\top \|_{A_{\boldsymbol{1}}}^2=(1-\gamma^2)\frac{\chi_{\max}}{\chi_{\min}}.
\ee
By \eqref{Eqdecomp}, \eqref{chiproj}, \eqref{Approj}, H\"older's and Young's inequalities, for any $\tau(\gamma) \in (0,1)$, 
\be\label{Pchiproj}
\bal
\| (E^\teq)^\top A_\EmAb P_{E}^\chi \|_{A_{\boldsymbol{1}}}^2 & = \| (E_1^\teq)^\top A_\EmAb P_{E}^\chi + (E_2^\teq)^\top A_\EmAb P_{E}^\chi \|_{A_{\boldsymbol{1}}}^2 = \| (E_1^\teq)^\top + (E_2^\teq)^\top A_\EmAb P_{E}^\chi \|_{A_{\boldsymbol{1}}}^2 \\
&= \|(E_1^\teq)^\top \|_{A_{\boldsymbol{1}}}^2 + \| (E_2^\teq)^\top A_\EmAb P_{E}^\chi \|_{A_{\boldsymbol{1}}}^2+ 2\left((E_1^\teq)^\top,(E_2^\teq)^\top A_\EmAb P_{E}^\chi \right)_{A_{\boldsymbol{1}}} \\
&\geq \|(E_1^\teq)^\top \|_{A_{\boldsymbol{1}}}^2 (1-\tau(\gamma)) + \| (E_2^\teq)^\top A_\EmAb P_{E}^\chi \|_{A_{\boldsymbol{1}}}^2 \left(1 - \frac{1}{\tau(\gamma)} \right) \\
&\geq   \gamma^2(1-\tau(\gamma)) + (1-\gamma^2)\frac{\chi_{\max}}{\chi_{\min}} \left(1-\frac{1}{\tau(\gamma)} \right)=:g(\gamma).
\eal
\ee
Let $\tau(\gamma)=\tfrac{1}{2}(1-\tfrac{\alpha^2}{\gamma^2})$.  Then for every $\gamma\in(\alpha,1]$, $\tau$ satisfies $0<\tau(\gamma)<1$ and $1-\tau(\gamma) = \tfrac12 + \tfrac{\alpha^2}{\gamma^2} >
\tfrac{\alpha^2}{\gamma^2}$. 
Since $g$ is continuous at $\gamma=1$ and $g(1) = 1 - \tau(1) >\alpha^2$, there exists $\gamma^*\in(\alpha,1)$, dependent on $\alpha$ and $\frac{\chi_{\max}}{\chi_{\min}}$, such that for any $\gamma^*\leq\gamma\leq 1$, $g(\gamma)\geq \alpha^2$.  Therefore by \eqref{Pchiproj} the result follows.
\end{proof}

We now have the following multi-step estimate.
\begin{theorem}\label{LstablemS}
Assume there exist constants $\beta \in (0,1)$ and $\alpha \in (0,1)$ such that  $\|P_{U^0}U^\teq\| \geq \beta$ and $\| (E^\teq)^\top A_{\chi} P_{E^{0}}^\chi \|_{A_{\chi}} \geq \alpha$.
Let $\gamma^* \in (\alpha,1)$ be given in \Cref{normrel}.
Then for any 
\be\label{deltacond}
0<\delta < \min\left\{ (1-c)\left(1+\frac{\chi_{\max}}{\chi_{\min}}\right)^{-1}, \frac{\sqrt{2(1-\max\{\gamma^*, \beta\}^2)} \|\varepsilon f_h^\teq\|_{L^2(\Omega)}}{\|\varepsilon (\hat f^{0}_h-f_h^{\teq})\|_{L^2(\Omega)}} \right\},
\ee
and $\Delta t_0$ given in \Cref{lem:onestep}, when $\Delta t \geq \Delta t_0$, 
\be\label{LstabS}
\|\varepsilon (\hat f_h^{\mathfrak{n}+1} - f_h^{\teq})\|_{L^2(\Omega)} \leq (c+\delta_\chi)^{\mathfrak{n}+1} \|\varepsilon (\hat f^{0}_h-f_h^{\teq})\|_{L^2(\Omega)}\quad\forall\mathfrak{n} \geq 1,
\ee
where $c$ is given by \eqref{ratior} and $\delta_\chi = \left(1+\frac{\chi_{\max}}{\chi_{\min}}\right) \delta$.  Moreover, if $\hat{f}_h^0 = f^\teq_h$, then for any $\Delta t>0$,
$
\hat{f}_h^{\mathfrak{n}+1} = f^\teq_h.
$
\end{theorem}
\begin{proof}
We prove the result by the method of induction. For $\mathfrak{n}=0$, \eqref{LstabS} follows from the one-step result in \Cref{lem:onestep}; see \eqref{LstabS++}. We assume that \eqref{LstabS} holds for some $\mathfrak{n} \geq 1$, that is
\be\label{LstabSn}
\|\varepsilon (\hat f_h^{\mathfrak{n}} - f_h^{\teq})\|_{L^2(\Omega)} \leq (c+\delta_\chi)^{\mathfrak{n}} \|\varepsilon (\hat f^{0}_h-f_h^{\teq})\|_{L^2(\Omega)}.
\ee
By \eqref{deltacond}, $(c+\delta_\chi) < 1$; thus $\|\varepsilon (\hat f_h^{\mathfrak{n}} - f_h^{\teq})\|_{L^2(\Omega)} \leq \|\varepsilon (\hat f^{0}_h-f_h^{\teq})\|_{L^2(\Omega)}$.
Then for $\mathfrak{n}+1$, the bounds in \eqref{basebdd0}, the fact that $|S^\teq|=\|\varepsilon f_h^\teq\|_{L^2(\Omega)}$ (see \Cref{seqbddlem}), and the definition of $\delta$ in \eqref{deltacond} imply that
\begin{subequations}\label{regularcond}
\begin{align}
\|P_{U^{\mathfrak{n}+1}}U^\teq\|^2 &\geq  1-\frac{\delta^2}{2\|\varepsilon f_h^\teq\|_{L^2(\Omega)}^2} \|\varepsilon(\hat{f}_h^{\mathfrak{n}}- f^\teq_h) \|_{L^2(\Omega)}^2 \geq \beta^2, \\ 
\|(E^\teq)^\top A_{\boldsymbol{1}} P_{E^{\mathfrak{n}+1}}\|_{A_{\boldsymbol{1}}}^2 &\geq  1-\frac{\delta^2}{2\|\varepsilon f_h^\teq\|_{L^2(\Omega)}^2} \|\varepsilon(\hat{f}_h^{\mathfrak{n}}- f^\teq_h) \|_{L^2(\Omega)}^2 \geq (\gamma^*)^2.
\end{align}
\end{subequations}
By \Cref{normrel}, (\ref{regularcond}b) implies $\|(E^\teq)^\top A_{\chi} P_{E^{\mathfrak{n}+1}}^\chi\|_{A_{\boldsymbol{1}}} \geq \alpha$. Therefore, the one-step estimate \eqref{LstabS++} holds.
The estimate (\ref{LstabS}) then follows from \eqref{LstabS++} and \eqref{LstabSn}.  Finally, if $\hat{f}_h^0 = f^\teq_h$, by \eqref{DLReq}, $\hat f_h^{\mathfrak{n}+1} = \hat f_h^{\mathfrak{n}} = \ldots = \hat f_h^{0} = f_h^\teq$.
\end{proof}
\section{Numerical Results}\label{sec-5}

In this section, we present numerical examples to validate our theoretical findings.  
For all the numerical tests in this section, we construct initial data $\hat{F}(0) = U^0 S^0 (E^0)^\top \in \mathcal{M}_r$ for \Cref{lralg} by applying the generalized singular value decomposition (GSVD) \cite{abdi2007} (\Cref{wsvd}) to $F(0)$, followed by truncation. 

\begin{example}\label{ex1}
In this example, we test the performance of the dynamical low-rank DG scheme in \eqref{fullDG1K}, or \Cref{lralg}, by comparing with the full-rank DG scheme in \eqref{fullDG1}.  
We let $\varepsilon_{\max}=1$, and set the opacity $\chi(\varepsilon)=4+\frac{\varepsilon^2}{2}$ and the emissivity $\eta(\varepsilon) = f^\teq(\varepsilon)\chi(\varepsilon)$, where 
\be\label{exeq}
f^\teq(\varepsilon) = \frac{1}{\varepsilon^2+1}.
\ee
is the rank-1 equilibrium distribution.
With initial data
$
f(\mu,\varepsilon,0) = \frac{1}{\varepsilon^2+1}+\frac{1}{\mu^2+\varepsilon^2+1/2},
$
the exact solution to \eqref{tp1} is
\be\label{exactsol}
f(\mu,\varepsilon,t) = \frac{1}{\varepsilon^2+1}+\frac{1}{\mu^2+\varepsilon^2+1/2}e^{-\chi(\varepsilon)t}.
\ee
We use $\cQ_2$ polynomials for all the tests in this example.  

To establish a baseline, we first test the spatial and temporal accuracy of the full-rank DG scheme in \eqref{fullDG1} with $N = N_\mu= N_\varepsilon$ cells in each direction. 
Errors at $t=1$ are shown in \Cref{acc_a}.  The convergence rate of the full rank DG scheme \eqref{fullDG1} is first-order in time (as expected with backward Euler time stepping) and third-order in phase-space (as expected with $\cQ_2$ polynomials) until saturation due to the temporal error. 
Errors at $t=10$ are shown in \Cref{acc_b}.  In this case, the phase-space convergence rate is still third-order for sufficiently small $\Delta t$, but the temporal accuracy is super linear due to the fact that the solution is very the near time-independent equilibrium distribution. Thus the error follows the bound in \eqref{Lstab}, which decreases geometrically.

\begin{figure}
\centering
\subfigure[Errors at $t=1$ .]{\includegraphics[width=0.49\textwidth]{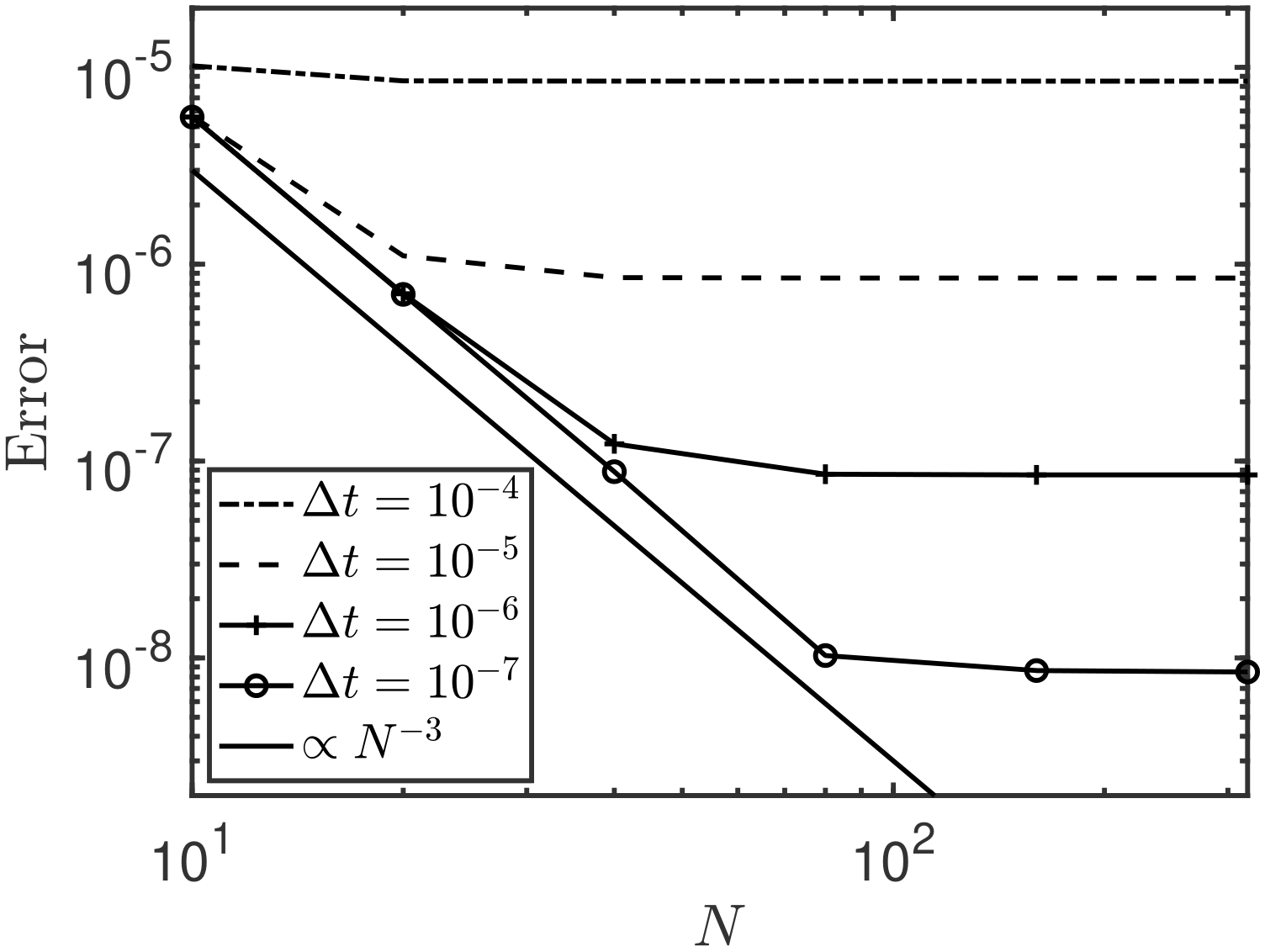}\label{acc_a}}
\subfigure[Errors at $t=10$.]{\includegraphics[width=0.495\textwidth]{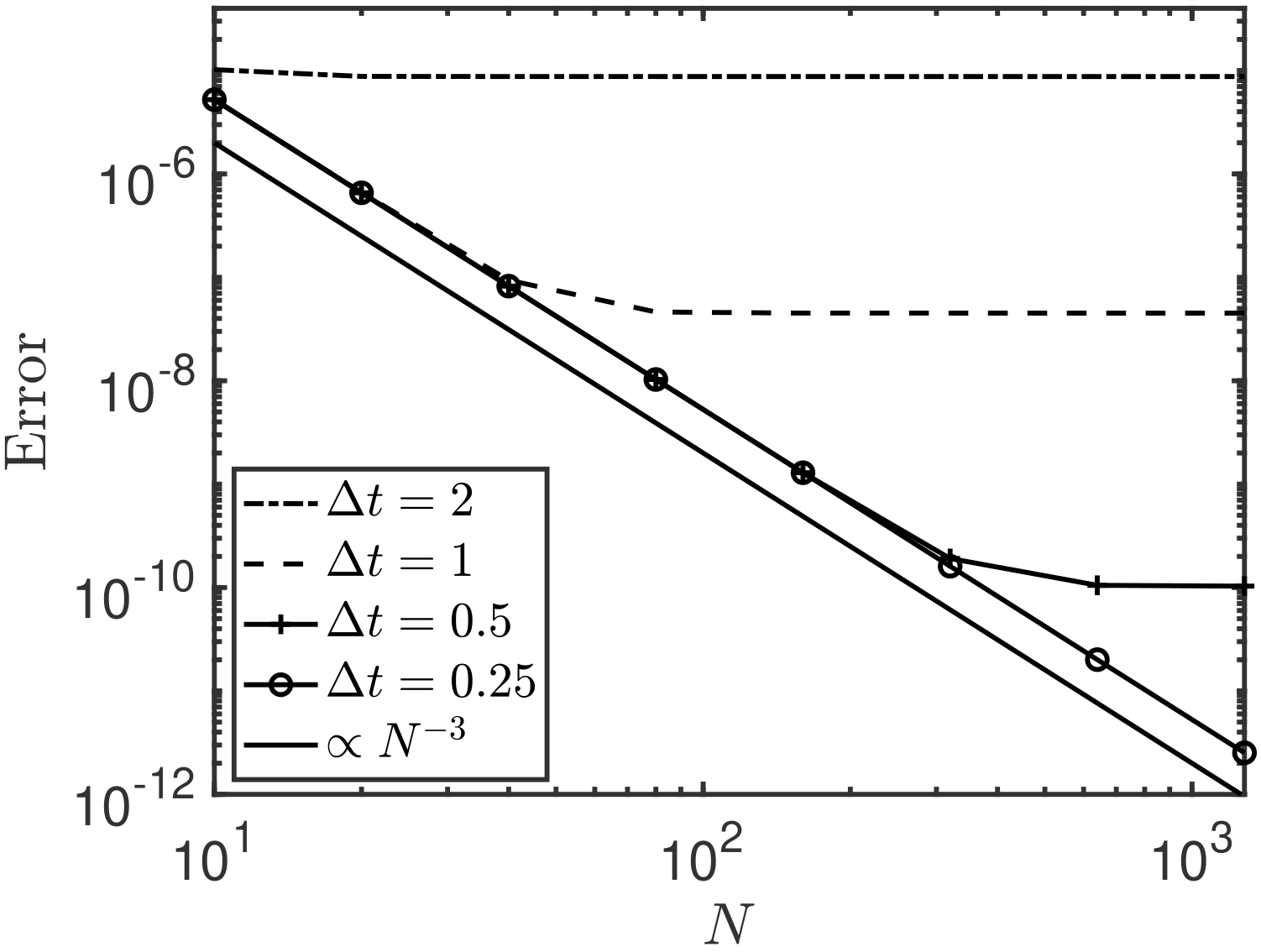}\label{acc_b}}
\caption{Error, $\|\varepsilon(f-f_h^\mathfrak{n})\|_{L^2(\Omega)}$, for the full-rank DG scheme in \eqref{fullDG1} versus number of elements, $N=N_\mu=N_\varepsilon$, for two different time step sizes.  The scheme uses $\cQ_2$ polynomials in phase-space and backward Euler time stepping.
In each panel, the solid lines without symbols are reference lines proportional to $N^{-3}$.
} \label{Accuracy}
\end{figure}

Second, we show the evolution of the rank of the coefficient matrix $F^n$ for the full-rank DG scheme \eqref{fullDG1}, using a mesh with $N_\mu= N_\varepsilon=160$. The numerical rank is calculated with the Matlab function \texttt{rank}($F^n$,$10^{-12}$), which returns the total number of singular values of $F^n$ that are larger than $10^{-12}$. 
The results with different time steps are plotted in \Cref{fullrank}.  
We observe that the numerical rank of the coefficient matrix decreases from $r=9$ at the initial condition to $r=1$ as the solution approaches equilibrium.

\begin{figure}
    \centering
    \subfigure[]{\includegraphics[width=0.49\textwidth]{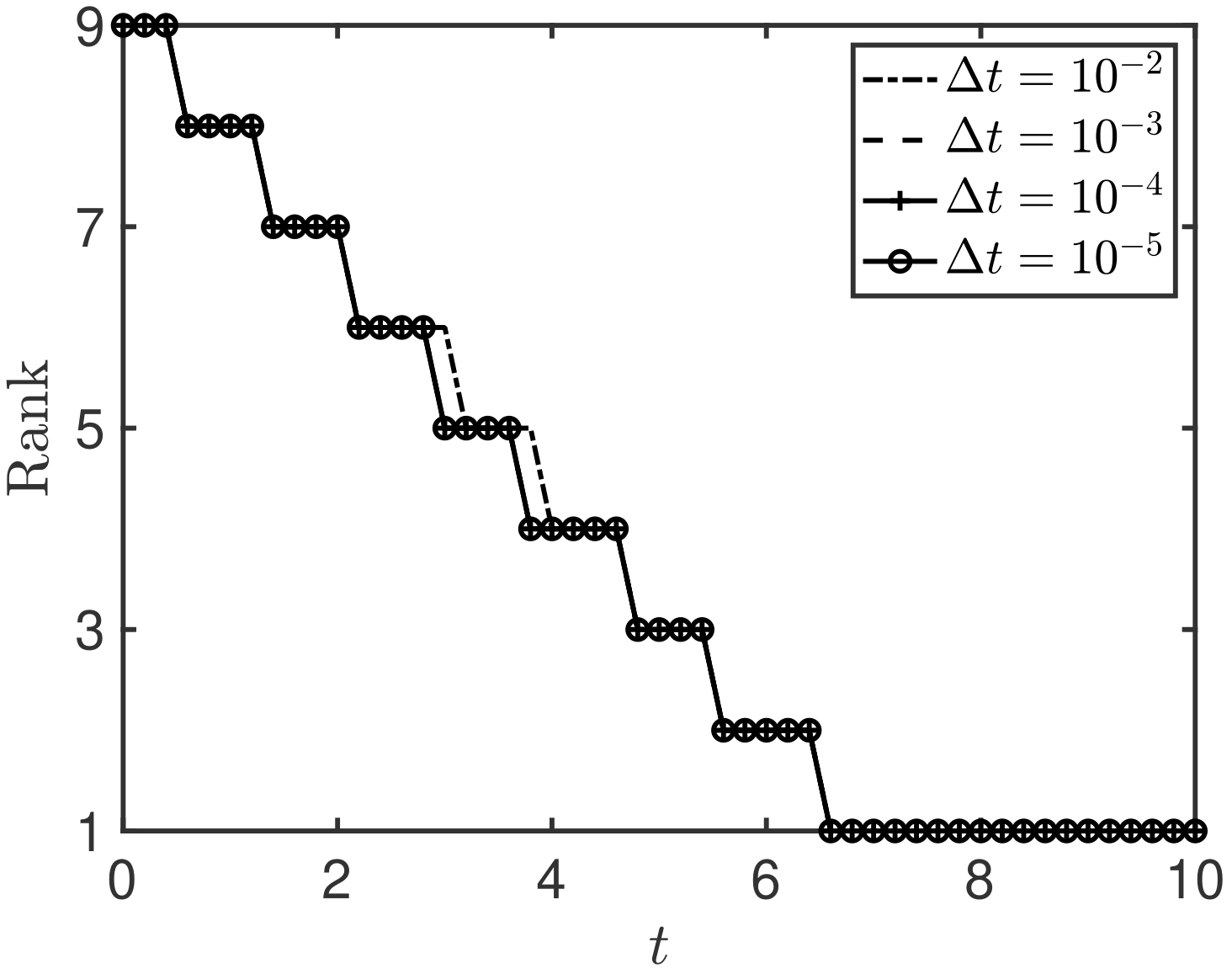}\label{fullrank}}
\subfigure[]{\includegraphics[width=0.495\textwidth]{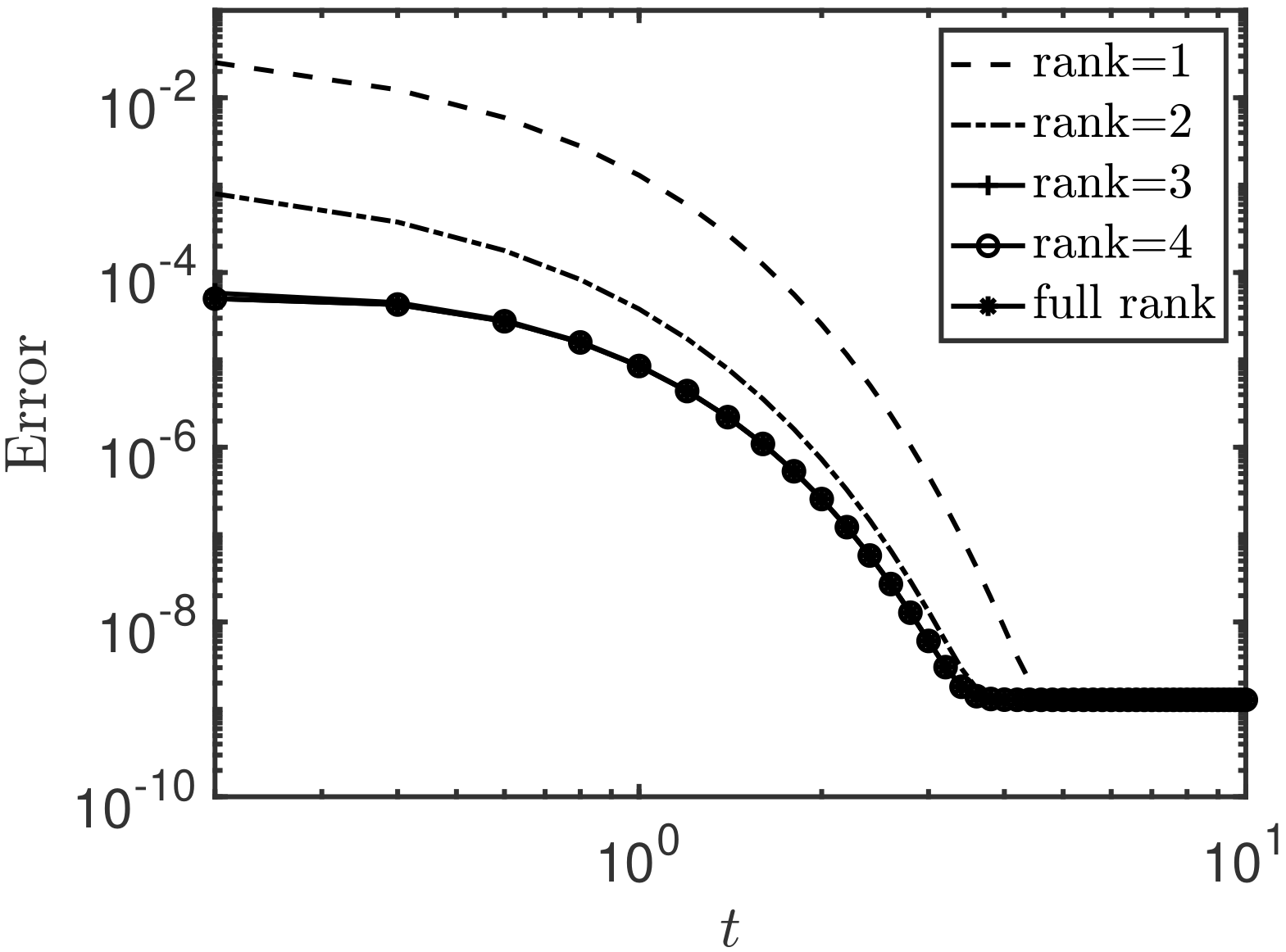}\label{lowrank}}
    \caption{(a) Evolution of the numerical rank for the coefficient matrix $F^n$ of the full-rank DG scheme, plotted vs. time using various time step sizes. (b) Weighted $L^2$ errors of the DLR-DG method (using $r=1$, $2$, $3$, and $4$) and the full rank DG scheme (\ref{fullDG1}) relative to the exact solution versus time, using $\Delta t = 10^{-4}$ and $N_{\mu}\times N_{\varepsilon}=160\times 160$.} 
\end{figure}

Third, we solve \Cref{tp1} using both the DLR-DG scheme in \Cref{lralg} and the full-rank DG scheme \eqref{fullDG1}. 
The purpose of this test is to compare the DLR-DG solution with the solution of the full-rank DG scheme \eqref{fullDG1} as the rank $r$ in \Cref{lralg} increases. The $L^2$ errors of the numerical solutions are plotted in Figure \ref{lowrank}. 
These errors decrease as the rank $r$ increases.  In particular, Algorithm \ref{lralg} with $r=3$ and time step $\Delta t = 10^{-4}$ produces numerical solutions that are practically identical to that of the full-rank scheme (\ref{fullDG1}).  
All low-rank solutions eventually give accurate equilibrium approximations. 

Fourth, we test the convergence of the dynamical low-rank DG solution and the full-rank DG solution to the equilibrium with one time step $\Delta t = T$, and two time steps $\Delta t = T/2$ for some final time $T$.  
The $L^2$ error of the numerical solution as a function of $T$ is plotted in \Cref{lowrankBE}.  
The results show that both algorithms converge up to discretization error, and the convergence rates of both algorithms to the equilibrium is equal to the total number of the time steps (i.e., $\propto T^{-1}$ for one step and $\propto T^{-2}$ for two steps), which is consistent with the theoretical results of \Cref{LstablemS}, regarding the low-rank scheme, and \Cref{fullDG1lem}~(ii), regarding the full-rank scheme.  
The $L^2$ error saturates for large $T$, when it becomes dominated by the projection error of the equilibrium (around $10^{-12}$).  

 \begin{figure}
 \centering
 \subfigure{\includegraphics[width=0.50\textwidth]{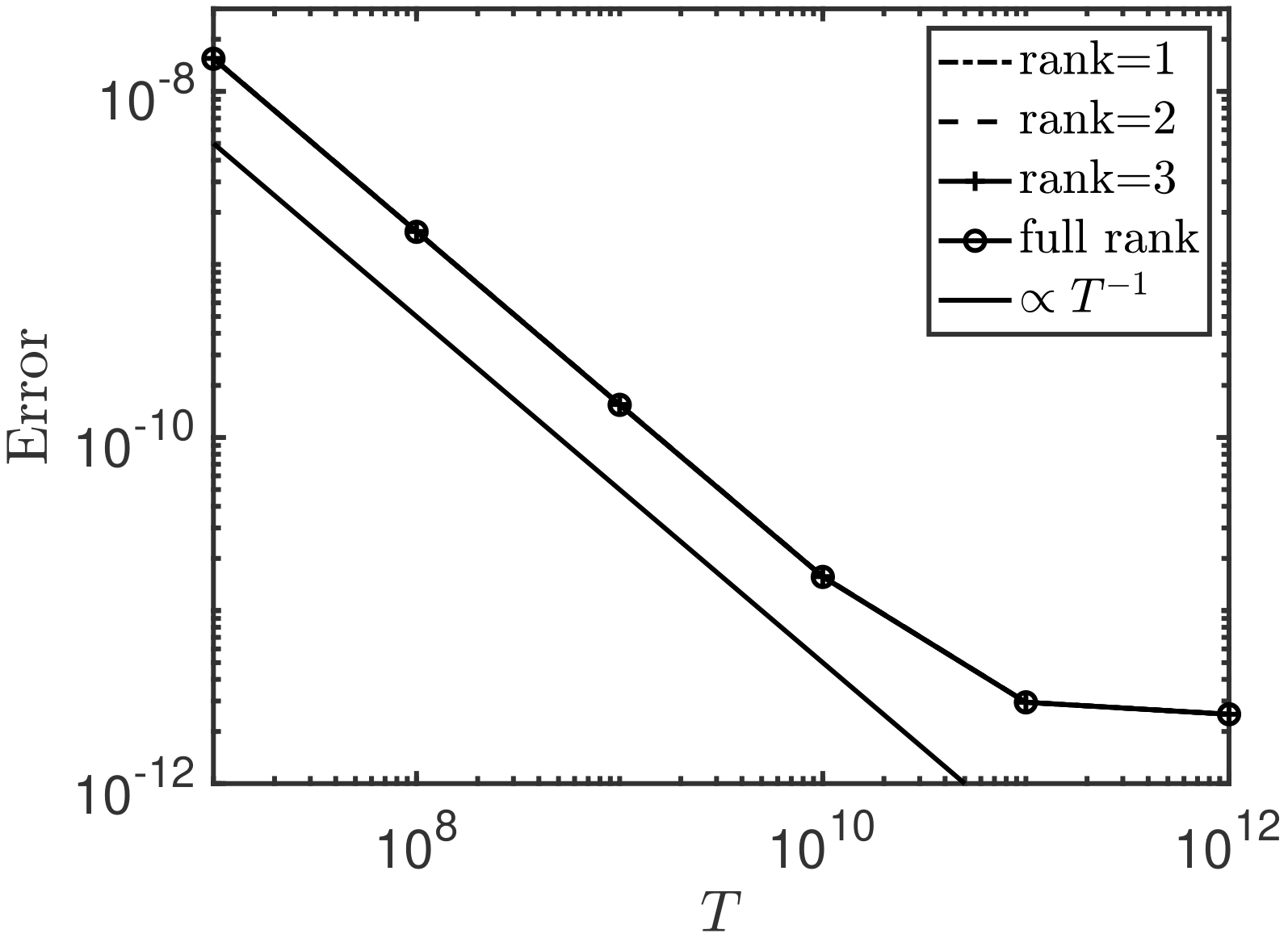}\label{lowrankBE_a}}
  \subfigure{\includegraphics[width=0.48\textwidth]{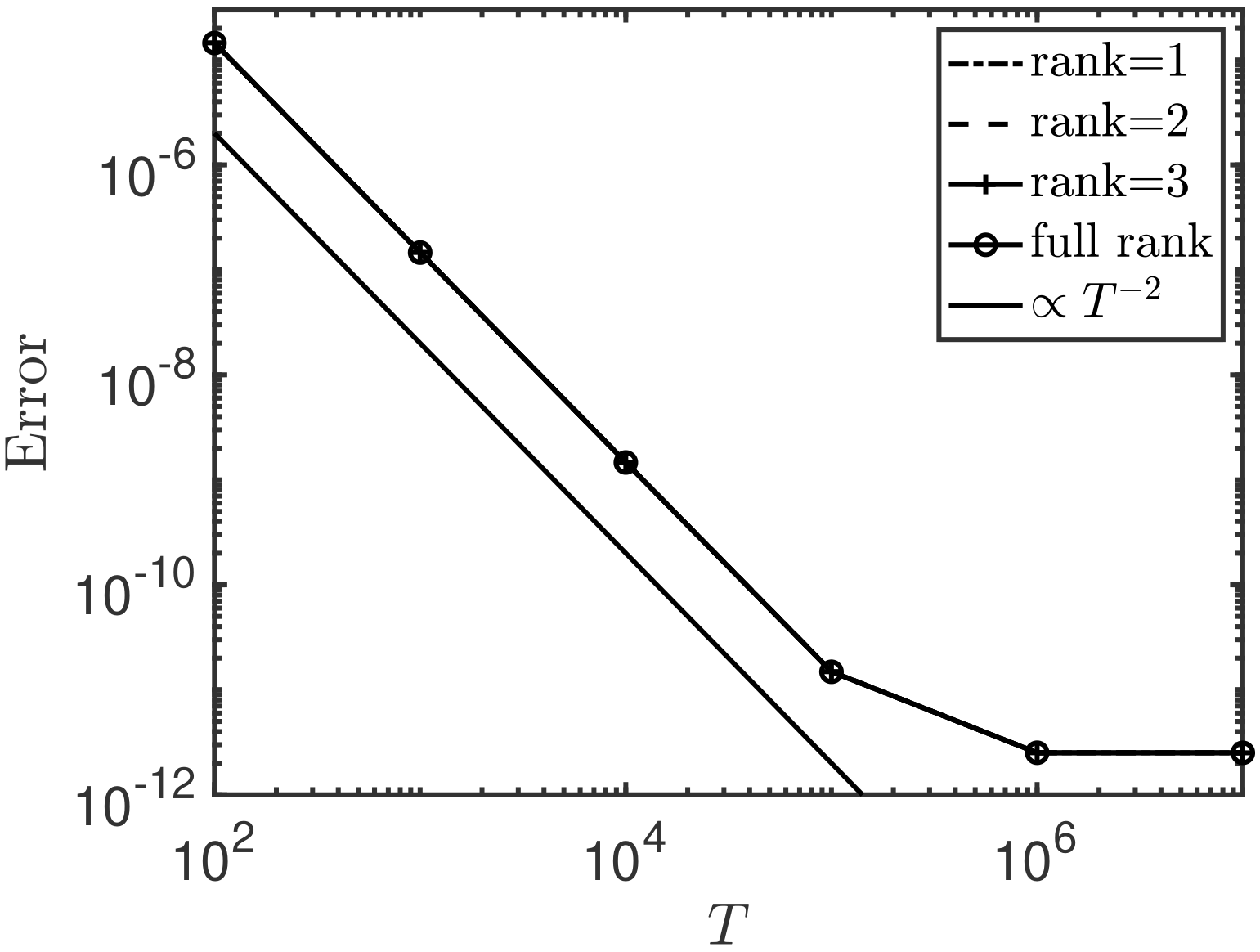}\label{lowrankBE_b}}
 \caption{Weighted $L^2$ errors of the dynamical low-rank DG method (with $r=1$, $2$, and $3$) and the full-rank DG scheme relative to the exact solution versus final time $T$, computed with one time step (with $\Delta t = T$; left panel) and two time steps (with $\Delta t= T/2$; right panel).} \label{lowrankBE}
 \end{figure}

Finally, we test the convergence of the dynamical low-rank DG solution and the full-rank DG solution to the equilibrium after $\mathfrak n$ steps, using
two differnt time step sizes: $\Delta t = 2$ and  $\Delta t = 10$.  
We show the $L^2$ error between the numerical solution and the exact solution versus $\mathfrak n$ in \Cref{lowrankRT}.  
The results show that both algorithms converge with convergence rates equal to the decay rate $c=\frac{1}{1+\Delta t\chi_{\min}}=\frac{1}{1+4\Delta t}$ (i.e., $\propto 9^{-1}$ for $\Delta t = 2$ and $\propto 41^{-1}$ for $\Delta t = 10$), which is consistent with the theoretical results of \Cref{LstablemS}, regarding the low-rank scheme, and \Cref{fullDG1lem}~(ii), regarding full-rank scheme.  

\begin{figure}
 \centering
 \subfigure{\includegraphics[width=0.49\textwidth]{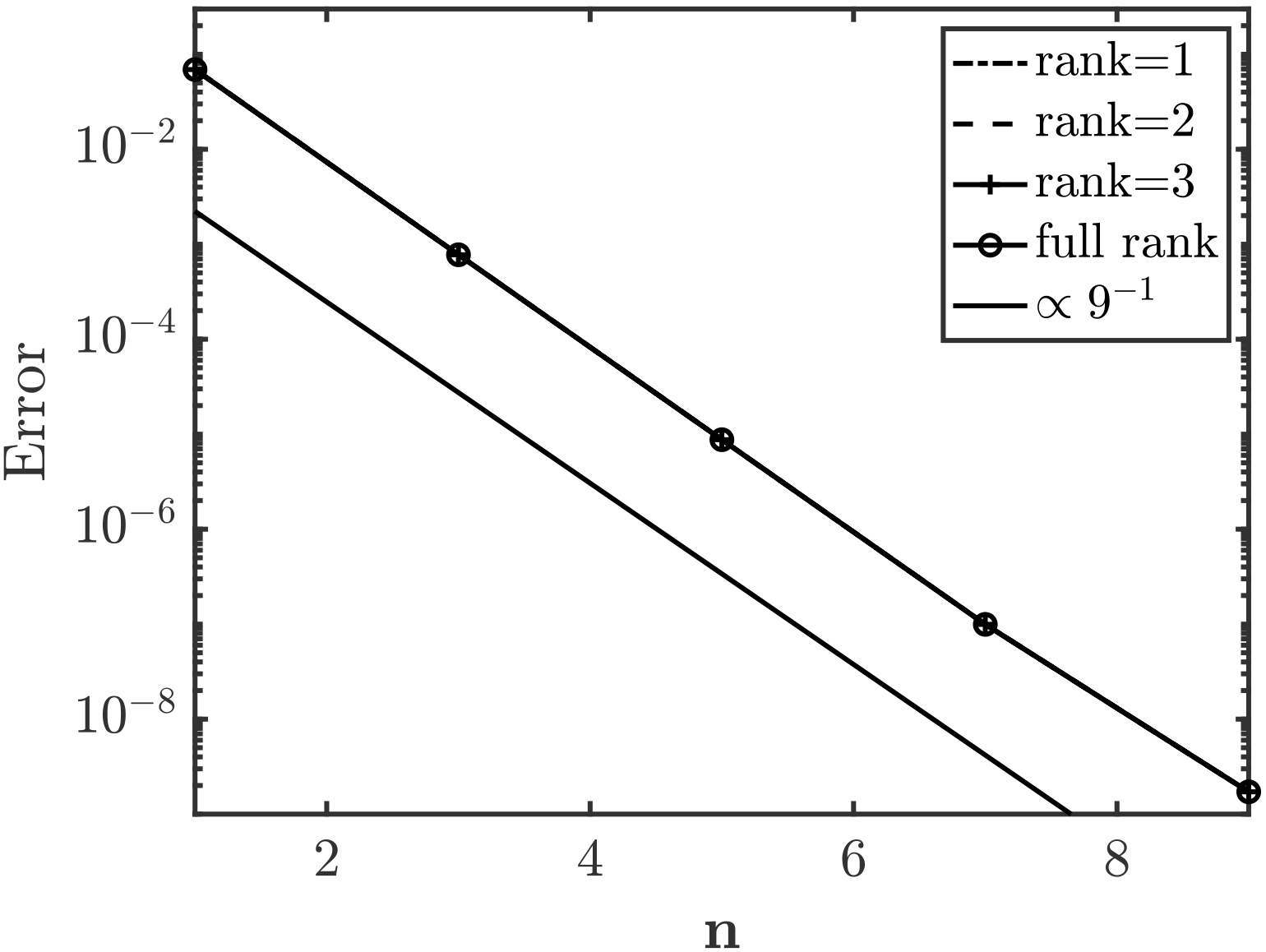}\label{lowrankRT_a}}
  \subfigure{\includegraphics[width=0.49\textwidth]{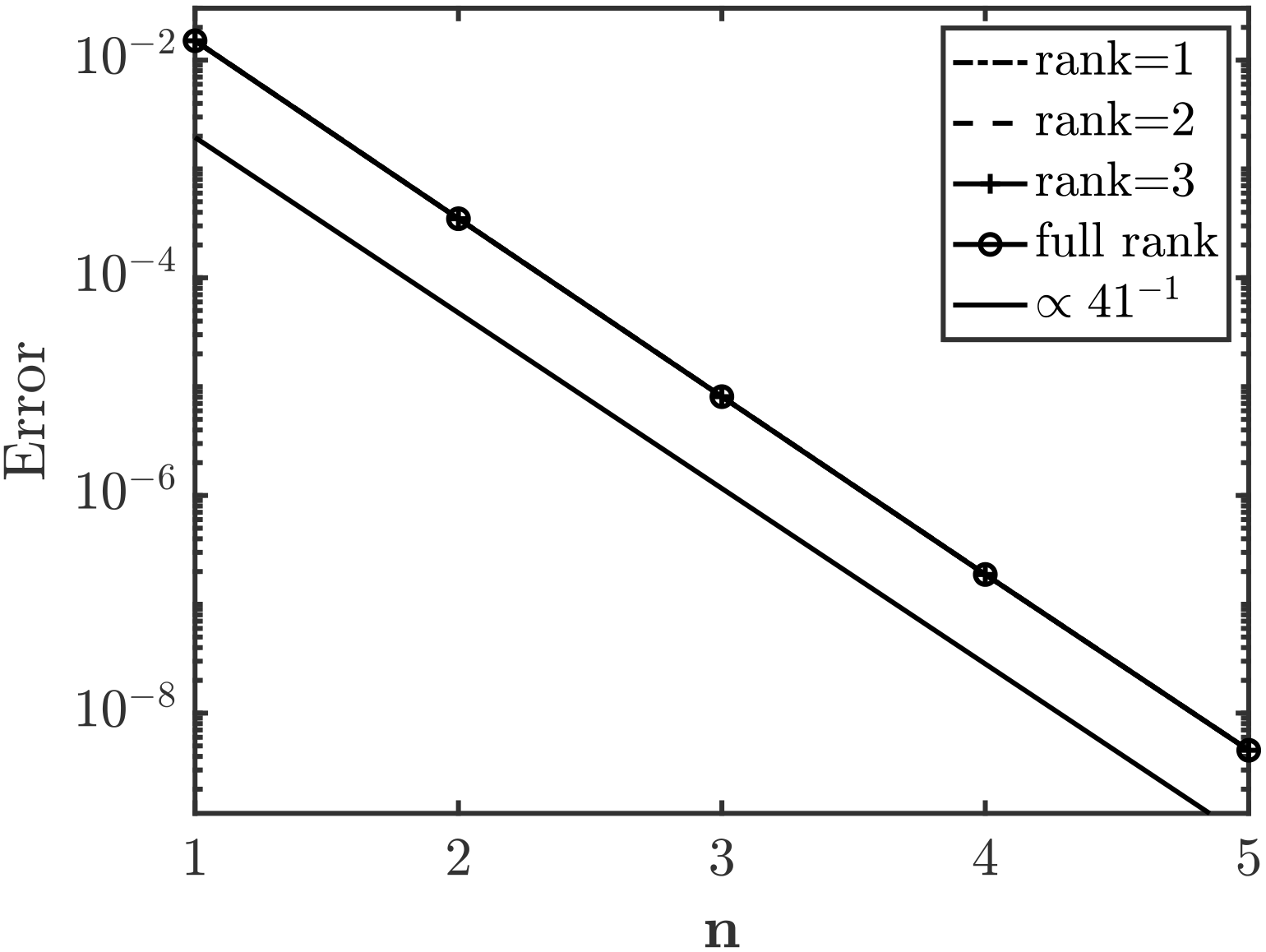}\label{lowrankRT_b}}
 \caption{Weighted $L^2$ errors of the dynamical low-rank DG method (with $r=1$, $2$, and $3$) and the full-rank DG scheme relative to the exact solution versus the total number of steps $\mathfrak n$, computed with ($\Delta t = 2$; left panel) and ($\Delta t= 10$; right panel).  In both panels, we compare the numerical results with the predicted decay rate $c=\frac{1}{1+\Delta t\chi_{\min}}=\frac{1}{1+4\Delta t}$.} \label{lowrankRT}
 \end{figure}

\end{example}

\begin{example}\label{ex2}
The purpose of this example is to demonstrate how the condition given in \Cref{lem:onestep} affects the convergence of the DLR-DG solution to the equilibrium.
We solve \eqref{tp1} with the same parameters as in Example \ref{ex1}, but with $\cQ_{1}$ polynomials and different initial conditions.  
The equilibrium is given in \eqref{exeq} and is independent of the initial data.  

To construct different initial conditions, we first prepare some basis functions.
\begin{enumerate}
\item[(i)] Let $\mathbf{K}^0=[U^\teq, \widetilde U^0]$, and $\mathbf{L}^0=[E^\teq, \widetilde E^0]$, where $U^\teq$, $E^\teq$ are given in (\ref{USEeq}), and  $\widetilde U^0$ and $\widetilde E^0$ are rank-$2$ matrices, computed from \Cref{lralg} using the initial data from Example \ref{ex1}.
\item[(ii)] Perform a QR factorization to obtain $\mathbf{K}^0=\hat U^0 R^0_U$, where $\hat U^0=[U^\teq, \hat U^0_2, \hat U^0_3]$.
\item[(iii)] Perform an $A_{\boldsymbol{1}}$-weighted Gram--Schmidt decomposition to obtain $\mathbf{L}^0=\hat E^0  R^0_E$, where $\hat E^0=[E^\teq, \hat E^0_2, \hat E^0_3]$.
\item[(iv)] Generate $\check{U}= \frac{U^\teq+\hat U^0_2}{\|U^\teq+\hat U^0_2\|}$ and $\check E = \frac{E^\teq+\hat E^0_2}{\|E^\teq+\hat E^0_2\|_{A_{\boldsymbol{1}}}}$. 
\item[(v)] Perform an $A_\EmAb$-weighted Gram--Schmidt decomposition to obtain $\mathbf{L}^0=\tilde E^0 \tilde R^0_E$, where $\tilde E^0 = [\tilde E^\teq, \tilde E^0_2, \tilde E^0_3]$.  Then perform an $A_{\boldsymbol{1}}$-weighted Gram--Schmidt decomposition to obtain $[\tilde E^0_2, \tilde E^0_3]=[\bar E^0_2, \bar E^0_3]  \bar R^0_{E}$.  Here, we expect $\|(\tilde E^0_j)^\top A_\EmAb E^\teq\|$ to be close to zero for $j=2,3$.
\end{enumerate}

\noindent\textbf{Test Case 5.2-1.} 
We use these different matrices to construct the various initial conditions given in the second and third row of \Cref{initcond}, with $S^0=1$.
We solve \eqref{tp1} with rank-$1$ initial conditions given in \Cref{initcond}, using \Cref{lralg} with $r=1$, $\cQ_1$ polynomials, and a mesh size of $N_\mu=N_\varepsilon=160$.  
We show the one time step ($\Delta t = T$) convergence of the dynamical low-rank DG solution to the equilibrium in \Cref{UEcondfig_a}. 
In \Cref{initcond}, we show the initial basis $U^{0}$ and $E^{0}$, the values in \eqref{UUEEasumpt} and \eqref{UUEEasumptU}, whether the assumptions of \Cref{lem:onestep} are satisfied (\cmark) or not (\xmark), and whether the scheme converges to the equilibrium (C) or not (NC).  
For Cases (a)-(c) in \Cref{initcond}, the conditions for convergence in \Cref{lem:onestep} are not satisfied, and the corresponding solution in \Cref{UEcondfig_a} does not converge to the equilibrium.
Case (d) is a special case that is addressed in \Cref{specialcase}. Specifically,
$\|(E^\teq)^\top A_\EmAb P_{E^0}^\chi\|_{A_{\boldsymbol{1}}}$ is zero (to algorithmic precision) and hence does not satisfy the associated condition in \Cref{fsnproj+}.  However, because $\left\|P_{U^0} U^\teq\right\|=1$, the corresponding numerical solution in \Cref{UEcondfig_a} still converges to the equilibrium with a first-order convergence rate.
Cases (e)-(f) satisfy the conditions of \Cref{lem:onestep}, and, as expected, the corresponding numerical solutions in \Cref{UEcondfig_a} converge to the equilibrium solution with a first-order convergence rate.  
All these results indicate that the conditions given in \Cref{lem:onestep}, or \Cref{specialcase}, are sufficient to determine the convergence of the one-step DLR-DG solution to the equilibrium.

\begin{table}[!htbp]\tabcolsep0.03in
\begin{tabular}[c]{||c|c|c|c|c|c|c||}
\hline
 & Case (a) & Case (b) & Case (c) & Case (d) & Case (e) & Case (f) \\
 \hline
$U^0$ & $\hat U^0_2$  & $\hat U^0_2$ & $\check{U}$ & $U^\teq$ & $\check{U}$ & $U^\teq$ \\
 \hline
$E^0$ & $\bar E^0_2$   & $E^\teq$ & $\bar E^0_2$ & $\bar E^0_2$ & $\check{E}$ & $E^\teq$ \\
\hline
$\left\|P_{U^0} U^\teq\right\|$ &  0 &  0 & $\sqrt{2}/2$ & 1 & $\sqrt{2}/2$ & 1 \\
\hline
$\| (E^\teq)^\top A_\EmAb P_{E^{0}}^\chi \|_{A_{\boldsymbol{1}}}$ &  1.7699e-15 & 1 & 1.7699e-15 & 1.7699e-15 & 0.7119 & 1 \\
\hline
Theorem \ref{LstablemS} &  \xmark & \xmark & \xmark & \cmark (\Cref{specialcase}) & \cmark & \cmark \\
\hline
Figure \ref{UEcondfig}(a) &  NC & NC & NC & C & C & C \\
\hline
\end{tabular}

\caption{Initial bases for Algorithm~\ref{lralg} used in Test Case~5.2-1, the values for the conditions in \eqref{UUEEasumpt} and \eqref{UUEEasumptU}, whether the conditions of \Cref{lem:onestep} are satisfied (\cmark) or not (\xmark), and the observed numerical behavior: convergence (C) or no convergence (NC).}
\label{initcond}
\end{table}

\begin{figure}
 \centering
 \subfigure[Example \ref{ex2} Test Case 5.2-1.]{\includegraphics[width=0.49\textwidth]{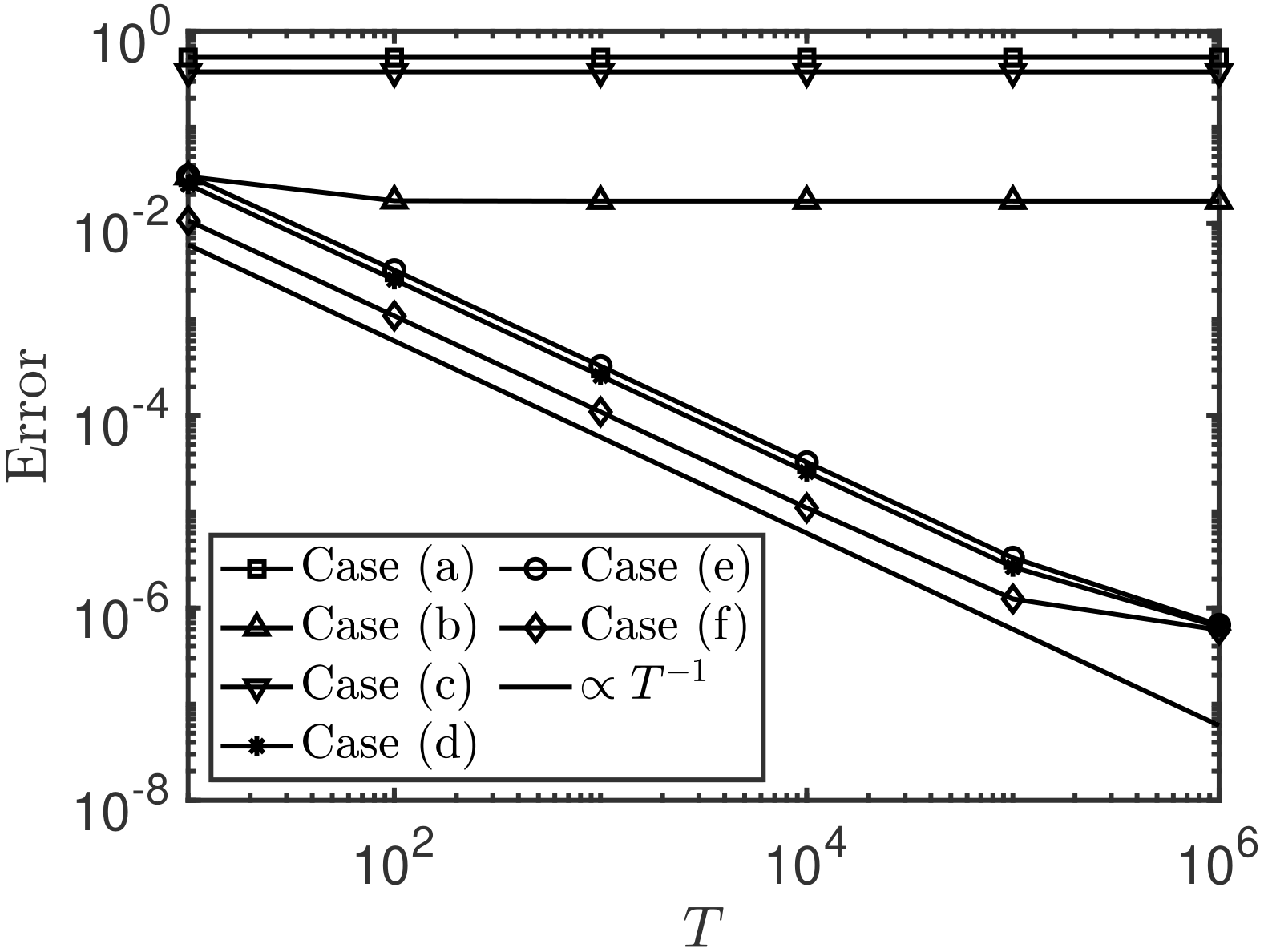}\label{UEcondfig_a}}
 \subfigure[Example \ref{ex2} Test Case 5.2-2.]{\includegraphics[width=0.49\textwidth]{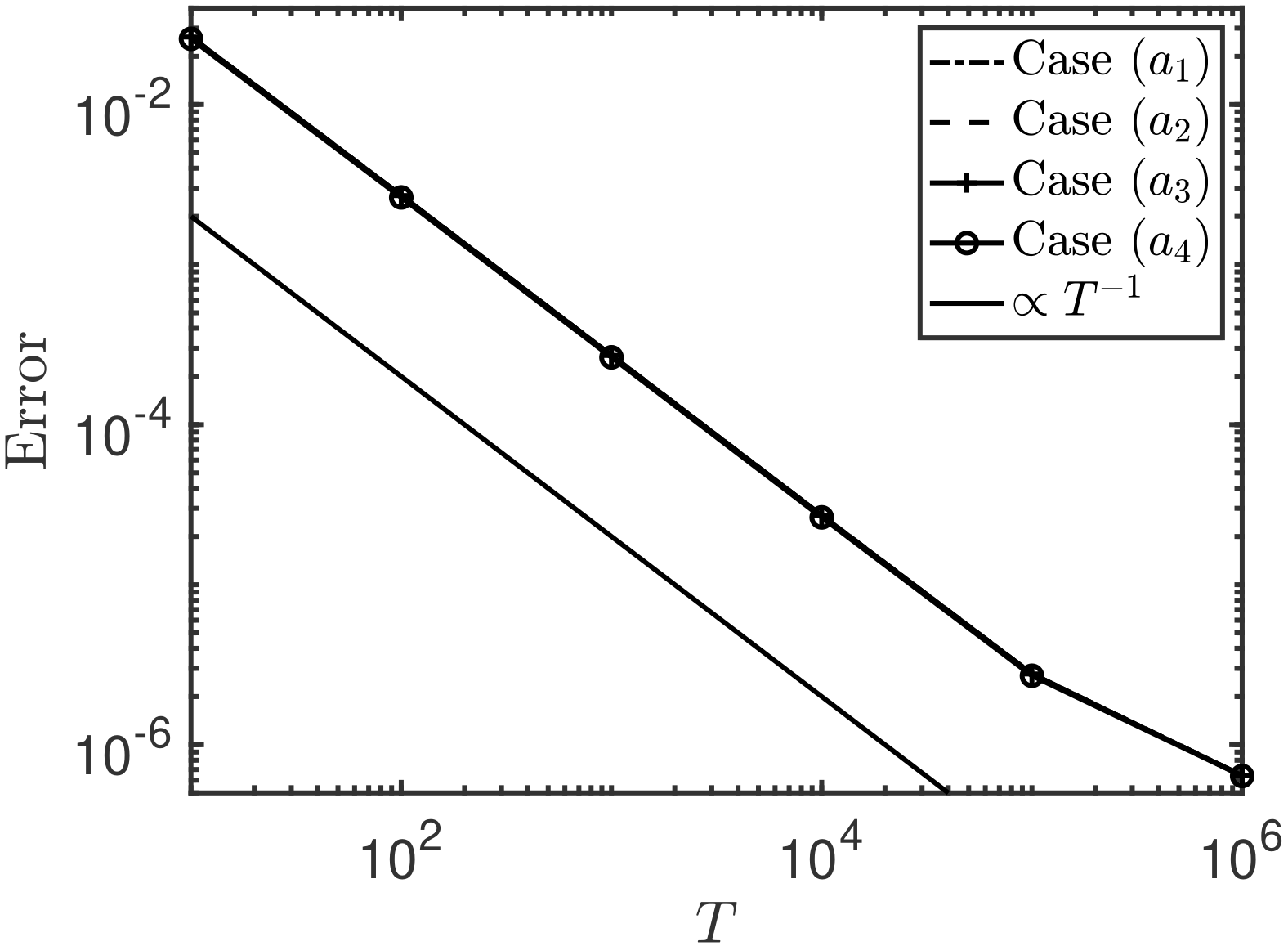}\label{UEcondfig_b}}
 \caption{The weighted $L^2$ errors between the dynamical low-rank DG solution and the equilibrium solution based on $\cQ_1$ polynomials and $N_\mu=N_\varepsilon=160$. }
 \label{UEcondfig}
 \end{figure}

\noindent\textbf{Test Case 5.2-2.} 
Though the conditions for convergence in \Cref{lem:onestep} are not satisfied for Case (a)-(c), the initial bases can be manually adjusted to yield a convergent algorithm.  
In the following, we take Case (a) in \Cref{initcond} as an example and show how to modify \Cref{lralg} such that the solution converges to the equilibrium.  
(Similar modifications can be applied to Case (b) and Case (c).)  
To achieve convergence, we increase the rank to $r=2$ and append the basis such that the conditions in \Cref{lem:onestep} or \Cref{specialcase} are satisfied.  
Let $x$ and $y$ be scalar parameters (not both zero), and define the functions
\begin{equation}
    U_\perp(x,y) := \frac{x U^\teq+y \hat U^0_3}{\|x U^\teq+y \hat U^0_3\|}
    \qquad \text{and} \qquad
    E_\perp(x,y) := \frac{xE^\teq+y\bar E^0_3}{\|xE^\teq+y\bar E^0_3\|_{A_{\boldsymbol{1}}}}.
\end{equation}
Then $\{\hat U^0_2, U_\perp(x,y)\}$ and $\{\bar E^0_2, E_\perp(x,y)\}$ are orthonormal and ${A_{\boldsymbol{1}}}$-orthonormal bases, respectively. 

\begin{table}[!htbp]\tabcolsep0.03in
\begin{tabular}[c]{||c|c|c|c|c||}
\hline
& Case ($a_1$) & Case ($a_2$) & Case ($a_3$) & Case ($a_4$) \\
 \hline
$U^0$ & $[\hat U^0_2,\ U_\perp(1,1)]$  & $[\hat U^0_2,\ U_\perp(1,0)]$ & $[\hat U^0_2,\ U_\perp(0.1,10)]$ & $[\hat U^0_2,\ U_{\text{rand}}]$ \\
 \hline
$E^0$ & $[\bar E^0_2,\ E_\perp(1,1)]$   & $[\bar E^0_2,\ E_\perp(0,1)]$ & $[\bar E^0_2,\ U_\perp(0.1,10)]$ & $[\bar E^0_2,\ E_{\text{rand}}]$ \\
 \hline
$\left\|P_{U^0} U^\teq\right\|$ & 0.7071 & 1 & 0.01 & 0.1601 \\
\hline
$\| (E^\teq)^\top A_\EmAb P_{E^{0}}^\chi \|_{A_{\boldsymbol{1}}}$ & 0.7046 & 3.0119e-13 & 0.01 & 0.0869 \\
\hline
\end{tabular}
\caption{Modified bases and the corresponding values for the condition in (\ref{UUEEasumpt}) and \eqref{UUEEasumptU}. }
\label{initmcond}
\end{table}

We use $U_\perp$ and $E_\perp$ to generate different initial bases for \Cref{lralg}; these are listed as Cases ($a_1$-$a_3$) in \Cref{initmcond}.  Case ($a_4$) is different: we randomly generate the basis functions by calling \texttt{randn}($(k+1)N$,1) in Matlab and apply the QR decomposition followed by an $A_{\boldsymbol{1}}$-weighted Gram--Schmidt decomposition to obtain the random basis functions $U_{\text{rand}}$ and $E_{\text{rand}}$, respectively. We set $S^0=\text{diag}(1,0)$, so that the initial matrix $F^0$ is unchanged after the basis enrichment.

The results are shown in \Cref{UEcondfig_b}, from which we see that after adding an additional component to the original bases, all of the initial conditions satisfy the conditions in \Cref{lem:onestep}, and consequently converge to the equilibrium.  
In addition, we also repeated Case ($a_4$) for more than $1000$ times with different random basis functions, and observed that all converged to the equilibrium.  
This is not surprising as the probability of drawing a random vector that is orthogonal to the equilibrium is very small.

\end{example}

\section{Conclusions}

In this paper, we have proposed a semi-implicit dynamical low-rank, discontinuous Galerkin (DLR-DG) method for a space homogeneous kinetic equation with a relaxation operator that models the emission and absorption of particles by a background medium.  
We have derived a weighted dynamical low-rank approximation (DLRA) that is consistent with the matrix differential equation of the DG scheme.  
A semi-implicit unconventional integrator (SIUI) is used to integrate the DLRA, and we show that the solution is identical to the solution of a DLR-DG scheme in a DLR-DG space.
We have shown the well-posedness of the fully discrete DLR-DG scheme and identified a sufficient condition on the time step size, together with conditions on the DLR-DG basis, such that the distance between the DLR-DG solution and the equilibrium solution decays geometrically with the number of time steps.  
Numerical results show that the DLR-DG solution is comparable to the full-rank DG solution and converges to the equilibrium solution when the bases satisfy the conditions of the theory.

In future work, it would be interesting to apply the proposed DLR-DG method to more general kinetic equations, e.g., that model scattering with a background.  
Then, in addition to the properties stated in \Cref{fullDG1lem} for the kinetic equation modeling  emission and absorption, the conservation of particles in the scattering process should be captured.  
It may be challenging for the proposed DLR-DG scheme to conserve particles, but extensions inspired by ideas proposed in \cite{einkemmer2022, baumann2023} may be fruitful.  
We will investigate this in future works.  

\section*{Acknowledgments}
Research at Oak Ridge National Laboratory is supported under contract DE-AC05-00OR22725 from the U.S. Department of Energy to UT-Battelle, LLC.
This work was supported by the U.S. Department of Energy, Office of Science, Office of Advanced Scientific Computing Research via the Applied Mathematics Program and the Scientific Discovery through Advanced Computing (SciDAC) program.  
This research was supported by Exascale Computing Project (17-SC-20-SC), a collaborative effort of the U.S. Department of Energy Office of Science and the National Nuclear Security Administration.

\appendix

\section{Some useful matrix results}

From Lemma \ref{FrobProp} to Lemma \ref{exteig}, we assume that $m$, $n$, and $r$ are some positive integers satisfying $r \leq \min\{m,n\}$.
If $A\in \mathbb{R}^{n\times n}$ is a symmetric and positive definite matrix, then Cholesky decomposition implies that there exists a nonsingular matrix $C \in \mathbb{R}^{n\times n}$ such that 
\be\label{AtoC}
A=C^\top C.
\ee

\begin{lemma}\label{FrobProp}
For any matrices $A\in \bbR^{m\times r}$, $B\in \bbR^{n\times r}$, and $D \in \bbR^{m\times n}$, 
\be\label{Fpropermat}
\left(AB^\top ,D\right)_{\rm{F}} =\left(B^\top ,A^\top D\right)_{\rm{F}} = \left(A,DB\right)_{\rm{F}}.
\ee
\end{lemma}

\begin{lemma}\label{eigtoray}
Let $a \in \mathbb{R}$ and $b \in \mathbb{R}$ be constants satisfying $0 \leq a \leq b$. Suppose $D \in \mathbb{R}^{n\times n}$ is a symmetric positive semi-definite matrix with eigenvalues $\{\lambda_i\}_{i=1}^n$ satisfying $a \leq \lambda_1 \leq \dots \leq \lambda_n \leq b$. 
Then for any nonzero $Z\in \mathbb{R}^{m \times n}$, 
\be\label{rayquo+}
a\leq \frac{(Z D, Z)_{\rm{F}}}{(Z, Z)_{\rm{F}}} \leq b.
\ee
\end{lemma}
\begin{proof}
For any nonzero $z\in \mathbb{R}^{n \times 1}$, the Rayleigh quotient satisfies
\be\label{rayquo}
a \leq \frac{(Dz,z)}{(z,z)} =  \frac{(z^\top D^\top, z^\top)}{(z^\top, z^\top)} \leq b.
\ee
 Set $Z^\top =[z_1,\ldots, z_m]$ where each $z_j\in \mathbb{R}^{n\times 1}$.  Then
\be\label{Rquotient1A-}
\frac{(ZD,Z)_{\rm{F}}}{(Z,Z)_{\rm{F}}} = \frac{\sum_{j=1}^m (Dz_j,z_j)}{\sum_{j=1}^m (z_j,z_j)},
\ee
which gives \eqref{rayquo+} by applying \eqref{rayquo} to each term in the sum of the numerator.
\end{proof}

\begin{lemma}\label{exteig}
Let $A \in \mathbb{R}^{n\times n}$ be a symmetric positive definite matrix. Suppose $D\in\mathbb{R}^{n\times n}$ with eigenvalues $\{\lambda_i\}_{i=1}^n$ satisfying $a \leq \lambda_1 \leq \dots \leq \lambda_n \leq b$. Then for any $Z \in \mathbb{R}^{m\times n}$,
\be\label{exteigformula}
0\leq (ZD^\top A,ZD^\top)_{\rm{F}} \leq b^2 (ZA,Z)_{\rm{F}}.
\ee
\end{lemma}
\begin{proof}
Let $(\lambda,q)$ be an eigenpair of the matrix $D$ so that
$Dq = \lambda q$.
If $q' = Cq$ for $C$ given in \eqref{AtoC}, then
$
CDC^{-1}q' = \lambda q'
$,
which implies that $\lambda$ is also the eigenvalue of the matrix $CDC^{-1}$
and  that the symmetric positive-definite matrix $(CDC^{-1})^T(CDC^{-1})$ has an eigenvalue $\lambda^2 \in [0,b^2]$. Let $Z' \in \mathbb{R}^{m \times n}$ be any matrix. By Lemma \ref{eigtoray}, we have 
\be\label{zcdccon}
0\leq{(Z'(CDC^{-1})^T(CDC^{-1}), Z')_{\rm{F}}} \leq b^2{(Z', Z')_{\rm{F}}},
\ee
which can be reformulated as \eqref{exteigformula}
by taking $Z'=ZC^T$.
\end{proof}

\section{Some useful algorithms}

Motivated by \cite{abdi2007}, we introduce the generalized singular value decomposition (GSVD) and the generalized QR factorization (GQR).  Let $\mathbb{S}^n_{++}$ be the set of
$n \times n$ symmetric positive definite matrices.
\begin{alg}\label{Aeigdecom}
(Matrix square root)
{Input:} $A_{\boldsymbol{1}}\in \mathbb{S}^n_{++}$. Output: $A_{\boldsymbol{1}}^{\pm \frac{1}{2}}\in \mathbb{S}^n_{++}$.
\begin{itemize}
    \item Apply the eigen-decomposition (\texttt{svd} in MATLAB) to $A_{\boldsymbol{1}}$ and obtain 
    \be
    A_{\boldsymbol{1}} = \Phi \Lambda \Phi^\top,
    \ee
    where  $\Phi$ satisfies $\Phi^\top \Phi = I_n$ and $\Lambda = \text{diag}(\lambda_1, \ldots, \lambda_n)$ with $\lambda_i>0$.
    \item Compute $\Lambda^{\pm \frac{1}{2}} = \text{diag}(\lambda_1^{\pm\frac{1}{2}}, \ldots, \lambda_n^{\pm\frac{1}{2}})$.
    \item Compute the symmetric matrix $A_{\boldsymbol{1}}^{\pm\frac{1}{2}} = \Phi \Lambda^{\pm\frac{1}{2}} \Phi^\top$. 
\end{itemize}
\end{alg}
Algorithm \ref{Aeigdecom} gives 
\be
A_{\boldsymbol{1}} = A_{\boldsymbol{1}}^{\frac{1}{2}} (A_{\boldsymbol{1}}^{\frac{1}{2}})^\top = A_{\boldsymbol{1}}^{\frac{1}{2}} A_{\boldsymbol{1}}^{\frac{1}{2}}, \quad A_{\boldsymbol{1}}^{\frac{1}{2}} A_{\boldsymbol{1}}^{-\frac{1}{2}} = I_n.
\ee

\begin{alg}\label{wsvd} (GSVD)
Input: $F \in \mathbb{R}^{m\times n}$, $A_{\boldsymbol{1}}^{\pm\frac{1}{2}}\in \mathbb{S}^n_{++}$, $r \leq\min\{m,n\}$. Output: $U \in \mathbb{R}^{m\times r}$, $S\in \mathbb{R}^{r\times r}$ and $E\in \mathbb{R}^{n\times r}$. 
\begin{itemize}
    \item Apply the SVD decomposition to $ F A_{\boldsymbol{1}}^{\frac{1}{2}}$ and obtain
    \be\label{Fsvd}
     F A_{\boldsymbol{1}}^{\frac{1}{2}} = US\hat{E}^\top,
    \ee
    where $U$ satisfies 
    \be\label{UUTI}
    U^\top U = I_r,
    \ee
    and $\hat{E}$ satisfies $\hat{E}^\top \hat{E} = I_r$.
    \item Compute $E = A_{\boldsymbol{1}}^{-\frac{1}{2}}\hat{E}$.
\end{itemize}
\end{alg}
\Cref{wsvd} gives the GSVD
\be
F=USE^T,
\ee
where $U$ satisfies \eqref{UUTI} and $E$ satisfies
\be\label{newEbasis}
E^\top A_{\boldsymbol{1}} E =  \hat{E}^\top A_{\boldsymbol{1}}^{-\frac{1}{2}} A_{\boldsymbol{1}} A_{\boldsymbol{1}}^{-\frac{1}{2}} \hat{E} = I_r.
\ee

\begin{alg}\label{wQR} (GQR)
Input: $\mathbf{L} \in \mathbb{R}^{n\times r}$, $A_{\boldsymbol{1}}^{\pm\frac{1}{2}}\in \mathbb{S}^n_{++}$. Output: $E\in \mathbb{R}^{n \times r}$.
\begin{itemize}
    \item Apply the QR decomposition to $A_{\boldsymbol{1}}^{\frac{1}{2}} \mathbf{L}$ and obtain
    \be
    A_{\boldsymbol{1}}^{\frac{1}{2}} \mathbf{L} = \hat{E} R,
    \ee
    where $\hat{E}$ satisfies $\hat{E}^\top \hat{E} = I_r$.
    \item Compute $E = A_{\boldsymbol{1}}^{-\frac{1}{2}} \hat{E}$.
\end{itemize}
\end{alg}
\Cref{wQR} gives the generalized QR factorization 
\be
\mathbf{L} = ER,
\ee
where $E$ also satisfies \eqref{newEbasis}.

\section{Technical Proofs}

In this section, we present the proofs to some lemmas. 
For any nonzero function $w_h=X^\top(\mu)WY(\varepsilon) \in V_h$ for some nonzero $W \in \mathbb{R}^{m \times n}$, let
\be\label{Rquotient}
R_{\chi}(w_h) = \frac{(\chi w_h,w_h;\varepsilon^2 )_{\Omega}}{(w_h,w_h;\varepsilon^2)_{\Omega}} = \frac{(WA_\EmAb,W)_{\rm{F}}}{(WA_{\boldsymbol{1}},W)_{\rm{F}}} = \frac{\|W\|_{A_{\chi}}^2}{ \|W\|_{A_{\boldsymbol{1}}}^2 }.
\ee

\begin{lemma}\label{Rquolem}
Let $m_1$ be an integer satisfying $1 \leq m_1 \leq m$. 
Then for any nonzero matrix $Z \in \mathbb{R}^{m_1\times n}$,
\begin{align} 
\label{Ralquo}
\chi_{\min}  \leq & \frac{\|Z\|_{A_\chi}^2}{\|Z\|^2_{A_{\boldsymbol{1}}}} \leq \chi_{\max},
 \end{align}
Therefore if $\lambda$ is an eigenvalue of the matrix $(A_\EmAb)^{-1}A_{\boldsymbol{1}}$ or $A_{\boldsymbol{1}}(A_\EmAb)^{-1}$, then
\be\label{AAEAeig}
\chi_{\max}^{-1} \leq \lambda \leq \chi_{\min}^{-1}.
\ee
\end{lemma}
\begin{proof}
A consequence of \Cref{remchi} is that
\be\label{quobdd}
\chi_{\min} \leq R_\chi(w_h) \leq \chi_{\max}.
\ee
The inequality in \eqref{Ralquo} follows from setting $w_h=X^\top(\mu)WY(\varepsilon)$ in \eqref{quobdd}, where  $W^\top=[Z^\top, Z_1^\top]$ and $Z_1=0 \in \mathbb{R}^{(m-m_1)\times n}$.  Inverting \eqref{Ralquo}, gives
\begin{equation}
\label{Rinverse}
   \chi_{\max}^{-1}  \leq  \frac{\|Z\|^2_{A_{\boldsymbol{1}}}}{\|Z\|_{A_\chi}^2} \leq \chi_{\min}^{-1} \quad \text{for all nonzero $Z \in \mathbb{R}^{m_1\times n}$}
\end{equation}
The inequalities in  \eqref{AAEAeig} follow immediately by setting $Z^\top$ in \eqref{Rinverse} to be an eigenvector of  $(A_\EmAb)^{-1}A_{\boldsymbol{1}}$.
\end{proof}

\subsection{Proof of \Cref{fsnproj}}\label{pfsnproj}

We will first need a rather technical lemma.

\begin{lemma}\label{LstepQR}
Let $E^\teq \in \mathbb{R}^{n\times 1}$, $B_{\mathbf{L}}=[b_1 \ldots, b_r] \in \mathbb{R}^{n\times r}$, and $l=[l_1,\ldots,l_r]$ be a nonzero vector, where $b_i \in \mathbb{R}^{n\times 1}$ and $l_i \in \mathbb{R}$ for $i=1,\dots,r$.
Assume that the matrix
\be\label{Ln1eqb}
\mathbf{L}^{\mathfrak{n}+1}
= \left[ l_1 E^{\teq} + \frac{1}{\Delta t} {b}_1, \cdots, l_r E^{\teq} + \frac{1}{\Delta t} {b}_r\right ] \in \mathbb{R}^{n\times r},
\ee
has a decomposition $\mathbf{L}^{\mathfrak{n}+1} = E^{\mathfrak{n}+1} S_{\mathbf{L}}^{\mathfrak{n}+1}$ with $E^{\mathfrak{n}+1} = \left[ E^{\mathfrak{n}+1}_1, \cdots, E^{\mathfrak{n}+1}_r  \right ]$ satisfying \eqref{Fn1cond}. 
Then
\be\label{EeqEnbdd}
1-\|(E^\teq)^\top A_{\boldsymbol{1}} P_{E^{\mathfrak{n}+1}}\|_{A_{\boldsymbol{1}}}^2 = 1-\|(E^{\mathfrak{n}+1})^\top  A_{\boldsymbol{1}} E^\teq\|^2  \leq \frac{\| B_{\mathbf{L}}^\top \|_{A_{\boldsymbol{1}}}^2}{\Delta t^2 \|l\|_\infty^2}.
\ee
\end{lemma}
\begin{proof}
As long as \eqref{Fn1cond} holds, \eqref{EeqEnbdd} is independent of the choice of basis for the span of $\mathbf{L}^{\mathfrak{n}+1}$.  Hence without loss of generality, we assume a weighted Gram-Schmidt decomposition:
\be\label{Ein1}
E^{\mathfrak{n}+1}_i = \frac{ \mathbf{L}^{\mathfrak{n}+1}_i - \sum_{j=1}^{i-1} \left( (\mathbf{L}^{\mathfrak{n}+1}_i)^\top A_{\boldsymbol{1}} E^{\mathfrak{n}+1}_j \right)E^{\mathfrak{n}+1}_j }{\sqrt{ (\mathbf{L}^{\mathfrak{n}+1}_i)^\top A_{\boldsymbol{1}} \mathbf{L}^{\mathfrak{n}+1}_i - \sum_{j=1}^{i-1} \left( (\mathbf{L}^{\mathfrak{n}+1}_i)^\top A_{\boldsymbol{1}} E^{\mathfrak{n}+1}_j \right)^2 } },
\ee
where $\mathbf{L}^{\mathfrak{n}+1}_i = l_i E^{\teq} + \frac{1}{\Delta t} {b}_i$.
Then
\be\label{EAE2}
\bal
((E^\teq)^\top  A_{\boldsymbol{1}} E^{\mathfrak{n}+1}_i)^2 
= & \frac{ \left( (E^\teq)^\top  A_{\boldsymbol{1}}\mathbf{L}^{\mathfrak{n}+1}_i - \sum_{j=1}^{i-1} \left( (\mathbf{L}^{\mathfrak{n}+1}_i)^\top A_{\boldsymbol{1}} E^{\mathfrak{n}+1}_j \right) (E^\teq)^\top  A_{\boldsymbol{1}} E^{\mathfrak{n}+1}_j \right)^2 }{ (\mathbf{L}^{\mathfrak{n}+1}_i)^\top A_{\boldsymbol{1}} \mathbf{L}^{\mathfrak{n}+1}_i - \sum_{j=1}^{i-1} \left( (\mathbf{L}^{\mathfrak{n}+1}_i)^\top A_{\boldsymbol{1}} E^{\mathfrak{n}+1}_j \right)^2  }\\
= &  \frac{ \left( l_i \xi_i^2 + \frac{1}{\Delta t} \alpha_i \right)^2 }{ l_i^2 \xi_i^2 + \frac{1}{\Delta t} (2l_i\alpha_i + \frac{1}{\Delta t}\gamma_i^2) },
\eal
\ee
where
\be\label{notagx}
\bal
\alpha_i &= (E^\teq)^\top  A_{\boldsymbol{1}} b_i - \sum_{j=1}^{i-1}  ((E^\teq)^\top  A_{\boldsymbol{1}} E^{\mathfrak{n}+1}_j)( (E^{\mathfrak{n}+1}_j)^\top  A_{\boldsymbol{1}} b_i ),\\
\gamma_i &= \left(b_i^\top  A_{\boldsymbol{1}} b_i  - \sum_{j=1}^{i-1} (b_i^\top A_{\boldsymbol{1}} E^{\mathfrak{n}+1}_j)^2\right)^\frac{1}{2}\\
\xi_{i} &= \left( (E^\teq)^\top  A_{\boldsymbol{1}} E^\teq-\sum_{j=1}^{i-1} ((E^\teq)^\top  A_{\boldsymbol{1}} E^{\mathfrak{n}+1}_j)^2 \right)^{\frac{1}{2}} = \left( 1-\sum_{j=1}^{i-1} ((E^\teq)^\top  A_{\boldsymbol{1}} E^{\mathfrak{n}+1}_j)^2 \right)^{\frac{1}{2}}
\eal
\ee
are all non-negative.  We extend the orthonormal basis $E^{\mathfrak{n}+1}_j$ from $1\leq j \leq r$ to $1 \leq j \leq n$. Then $E^\teq$ and $b_i$ in \eqref{Ln1eqb} can be expressed in terms of the basis functions $\{E_j^{\mathfrak{n}+1}\}_{j=1}^n$ as
\be\label{EbtoEj}
\bal
E^\teq = \sum_{j=1}^n ((E^\teq)^\top A_{\boldsymbol{1}} E_j^{\mathfrak{n}+1})E_j^{\mathfrak{n}+1}, \quad b_i = \sum_{j=1}^n (b_i^\top A_{\boldsymbol{1}} E_j^{\mathfrak{n}+1})E_j^{\mathfrak{n}+1},
\eal
\ee
which implies
\be\label{gamxi}
\gamma_i^2= \sum_{j=i}^{n} (b_i^\top A_{\boldsymbol{1}} E^{\mathfrak{n}+1}_j)^2,\qquad \xi_i^2 = \sum_{j=i}^{n} ((E^\teq)^\top A_{\boldsymbol{1}} E^{\mathfrak{n}+1}_j)^2,
\ee
and
\be\label{alphas-}
\bal
\alpha_i= & \sum_{j=i}^{n}  ((E^\teq)^\top  A_{\boldsymbol{1}} E^{\mathfrak{n}+1}_j)( (E^{\mathfrak{n}+1}_j)^\top  A_{\boldsymbol{1}} b_i ).
\eal
\ee
Therefore, we have
\be\label{alphas}
\bal
|\alpha_i| \leq  \left(\sum_{j=i}^{n} \left((E^\teq)^\top  A_{\boldsymbol{1}} E^{\mathfrak{n}+1}_j \right)^2 \right)^{\frac{1}{2}} \left( \sum_{j=i}^{n}\left ( (E^{\mathfrak{n}+1}_j)^\top  A_{\boldsymbol{1}} b_i \right)^2) \right)^{\frac{1}{2}}
= \xi_i \gamma_i.
\eal
\ee
By \eqref{EbtoEj} and \eqref{gamxi}, it follows that
\be\label{gammatoB}
\gamma_i^2 \leq b_i^\top  A_{\boldsymbol{1}} b_i  \leq (B_{\mathbf{L}}^\top A_{\boldsymbol{1}}, B_{\mathbf{L}}^\top )_{\rm{F}} = \| B_{\mathbf{L}}^\top \|_{A_{\boldsymbol{1}}}^2.
\ee
Meanwhile, the direct calculation gives
\be\label{En1Eq}
\|(E^{\mathfrak{n}+1})^\top  A_{\boldsymbol{1}} E^\teq\|^2 = \sum_{j=1}^r ((E^\teq)^\top  A_{\boldsymbol{1}} E^{\mathfrak{n}+1}_j)^2.
\ee
Choose $i$ such that $1\leq i\leq r$ and $|l_i| = \|l\|_\infty := \max_{1\leq j\leq r} |l_j|$.  We consider the following cases.
\begin{enumerate}[leftmargin=0pt, label = \bf Case \arabic*\phantom{.a}]
    \item If $l_i \xi_i^2 + \frac{1}{\Delta t} \alpha_i  = 0$, that is $\xi^2_i = - \frac{\alpha_i}{\Delta t l_i}=\frac{|\alpha_i|}{\Delta t |l_i|}$, then 
    by \eqref{alphas},
\be\label{xibdd}
0\leq \xi_i \leq \frac{\gamma_i}{\Delta t |l_i|},
\ee
which together with \eqref{gammatoB} and \eqref{En1Eq} implies that
\be
1-\|(E^{\mathfrak{n}+1})^\top  A_{\boldsymbol{1}} E^\teq\|^2
\leq 1-\sum_{j=1}^{i-1}((E^\teq)^\top  A_{\boldsymbol{1}} E^{\mathfrak{n}+1}_j)^2 
= \xi_i^2 \leq \frac{\gamma_i^2}{\Delta t^2 l_i^2} \leq \frac{\| B_{\mathbf{L}}^\top \|_{A_{\boldsymbol{1}}}^2}{\Delta t^2 \|l\|_\infty^2}.
\ee
\item Now we consider $l_i \xi_i^2 + \frac{1}{\Delta t} \alpha_i  \not = 0$.
\begin{enumerate}[leftmargin=0pt,label= \bf Case \arabic{enumi}.\alph*]
    \item If $\xi_i = 0$, then
\be
1-\|(E^{\mathfrak{n}+1})^\top  A_{\boldsymbol{1}} E^\teq\|^2
    \leq 
    1- \sum_{j=1}^{i-1} ((E^\teq)^\top  A_{\boldsymbol{1}} E^{\mathfrak{n}+1}_j)^2
    = \xi_i^2 =0.
\ee
Therefore, the inequality \eqref{EeqEnbdd} holds.
\item If $\xi_i \not = 0$, we consider two cases:

\begin{enumerate}[leftmargin=0pt,label= \bf Case \arabic{enumi}.\alph{enumii}.\roman*]
\item 
If $\gamma_i=0$, then by \eqref{alphas}, $\alpha_i=0$. 
By \eqref{EAE2}, $((E^\teq)^\top  A_{\boldsymbol{1}} E^{\mathfrak{n}+1}_i)^2 = \xi_i^2$, which implies
\be\label{gamma0eq}
1-\|(E^{\mathfrak{n}+1})^\top  A_{\boldsymbol{1}} E^\teq\|^2
    \leq 1-\sum_{j=1}^{i}((E^\teq)^\top  A_{\boldsymbol{1}} E^{\mathfrak{n}+1}_j)^2 \\
    =\xi_i^2 - ((E^\teq)^\top  A_{\boldsymbol{1}} E^{\mathfrak{n}+1}_i)^2 = 0.
\ee
Therefore, the inequality \eqref{EeqEnbdd} still holds.\\

\item If $\gamma_i\neq 0$,
by \eqref{alphas} there exists a parameter $\tau \in [-1,1]$ such that
\be \label{alpha_equiv}
\alpha_i = \tau \gamma_i \xi_i.
\ee
Substituting \eqref{alpha_equiv} into \eqref{EAE2} and rewriting yield
\be\label{EAE20}
\bal
((E^\teq)^\top  A_{\boldsymbol{1}} E^{\mathfrak{n}+1}_i)^2 = &  \xi_i^2\frac{ \left( l_i \xi_i + \frac{1}{\Delta t} \tau \gamma_i  \right)^2 }{ \left( l_i \xi_i + \frac{1}{\Delta t} \tau \gamma_i  \right)^2 + \frac{1}{(\Delta t)^2}(1-\tau^2)\gamma_i^2 }\\
= & \xi_i^2 - g(\tau)
=  1-\sum_{j=1}^{i-1} ((E^\teq)^\top  A_{\boldsymbol{1}} E^{\mathfrak{n}+1}_j)^2 - g(\tau),
\eal
\ee
where we have used \eqref{notagx} for the third equality and $g:[-1,1]\to\R$ is a non-negative and differentiable function given by
\be\label{g_tau_def}
g(\tau) = \frac{ \xi_i^2 \frac{1}{(\Delta t)^2}(1-\tau^2)\gamma_i^2 }{ \left( l_i \xi_i + \frac{1}{\Delta t} \tau \gamma_i  \right)^2 + \frac{1}{(\Delta t)^2}(1-\tau^2)\gamma_i^2 }
\ee
We wish to maximize $g$ on $[-1,1]$.  Since $g(-1)=g(1)=0$, we solve for the critical points $\tau^*$ satisfying
\be\label{gprime}
g'(\tau)=\frac{-  \frac{2}{(\Delta t)^2}\xi_i^2\gamma_i^2 \left( l_i \xi_i + \frac{1}{\Delta t} \tau \gamma_i  \right) \left( \tau l_i \xi_i + \frac{1}{\Delta t} \gamma_i \right) }{ \left(\left( l_i \xi_i + \frac{1}{\Delta t} \tau \gamma_i  \right)^2 + \frac{1}{(\Delta t)^2}(1-\tau^2)\gamma_i^2 \right)^2}=0.
\ee
Since $\xi_i \not = 0$ and $\left(l_i \xi_i^2 + \frac{1}{\Delta t} \alpha_i\right)^2 = \xi_i^2 \left( l_i \xi_i + \frac{1}{\Delta t} \tau \gamma_i  \right)^2 \not =0$, 
it follows that $l_i \xi_i + \frac{1}{\Delta t} \tau \gamma_i \not = 0$. Therefore, the only critical point for \eqref{gprime} is $\tau^* = -\frac{\gamma_i}{\Delta t l_i \xi_i} \in (-1,1)$. 
Plugging in $\tau^*$ into \eqref{g_tau_def} yields
\be\label{gtaubdd}
g(\tau) \leq g\left(-\frac{\gamma_i}{\Delta t l_i \xi_i}\right) = \frac{\gamma_i^2}{\Delta t^2 l_i^2}.
\ee
Therefore \eqref{EAE20}, \eqref{gtaubdd}, and \eqref{gammatoB} imply
\be
1-\|(E^{\mathfrak{n}+1})^\top  A_{\boldsymbol{1}} E^\teq\|^2\leq 1-\sum_{j=1}^{i}((E^\teq)^\top  A_{\boldsymbol{1}} E^{\mathfrak{n}+1}_j)^2 = g(\tau) \leq \frac{\gamma_i^2}{\Delta t^2 l_i^2} \leq \frac{\| B_{\mathbf{L}}^\top \|_{A_{\boldsymbol{1}}}^2}{\Delta t^2 \|l\|_\infty^2}.
\ee
\end{enumerate}

\end{enumerate}

\end{enumerate}

\end{proof}

Next, we present the proof of \Cref{fsnproj}.

\begin{proof}[Proof of \Cref{fsnproj}]
Start with \eqref{fullDG1K-L}, which is equivalent to finding $\mathbf{L}^{\mathfrak{n}+1} \in \mathbb{R}^{n\times r}$ such that for any $\mathbf{L}_W\in \mathbb{R}^{n\times r}$, 
\be\label{halfmatodeKb}
\left(U^\mathfrak{n} \left(  D_t \mathbf{L}^{\mathfrak{n}+1}\right)^\top  A_{\boldsymbol{1}}, U^\mathfrak{n}\mathbf{L}_W^\top  \right)_{\rm{F}}
= \left(G(U^\mathfrak{n}(\mathbf{L}^{\mathfrak{n}+1})^\top ), U^\mathfrak{n}\mathbf{L}_W^\top \right)_{\rm{F}}.
\ee
Set $\mathbf{L}_W = (A_\EmAb)^{-1}\mathbf{L}'_W$, where $\mathbf{L}'_W$ is arbitrary, into \eqref{halfmatodeKb}.  Then use \eqref{Feqmat} for $G$, and apply Lemma \ref{FrobProp}:
\be\label{halfmatodeKb1}
\bal
\left(\left(  D_t \mathbf{L}^{\mathfrak{n}+1}\right)^\top  A_{\boldsymbol{1}}(A_\EmAb)^{-1}, (\mathbf{L}'_W)^\top  \right)_{\rm{F}}
= & \left((U^\mathfrak{n} )^\top F^{\teq}  - (\mathbf{L}^{\mathfrak{n}+1})^\top , (\mathbf{L}'_W)^\top \right)_{\rm{F}},
\eal
\ee
Since $\mathbf{L}'_W$ is arbitrary, it follows that
\bes
\bal
& \left(I_{n} + \frac{1}{\Delta t} (A_\EmAb)^{-1}A_{\boldsymbol{1}} \right)\mathbf{L}^{\mathfrak{n}+1}
=  (F^\teq)^\top   U^\mathfrak{n}  + \frac{1}{\Delta t}(A_\EmAb)^{-1}A_{\boldsymbol{1}} \mathbf{L}^\mathfrak{n}\\
& =  \left(I_{n} + \frac{1}{\Delta t} (A_\EmAb)^{-1}A_{\boldsymbol{1}} \right) (F^\teq)^\top    U^\mathfrak{n} + \frac{1}{\Delta t} (A_\EmAb)^{-1}A_{\boldsymbol{1}} \left( \mathbf{L}^\mathfrak{n} -  (F^\teq)^\top  U^\mathfrak{n}  \right),
\eal
\ees
which gives
\be\label{halfmatodeKb3}
\bal
\mathbf{L}^{\mathfrak{n}+1} 
&= (F^\teq)^\top  U^\mathfrak{n} + \frac{1}{\Delta t}B_{\mathbf{L}},
\eal
\ee
where
\begin{equation}\label{BBB}
\begin{aligned}
    B_{\mathbf{L}} & :=  D_{\mathbf{L}} \left(\mathbf{L}^\mathfrak{n} - (F^\teq)^\top  U^\mathfrak{n} \right) \in\mathbb{R}^{n\times r},\\
    D_{\mathbf{L}} &:= \left(I + \frac{1}{\Delta t} (A_\EmAb)^{-1}A_{\boldsymbol{1}} \right)^{-1} (A_\EmAb)^{-1}A_{\boldsymbol{1}}\in\mathbb{R}^{n\times n}.
\end{aligned}
\end{equation}
Because $U^\mathfrak{n}$ is orthogonal and $U^\mathfrak{n}(\mathbf{L}^\mathfrak{n})^\top=U^\mathfrak{n}S^\mathfrak{n}(E^\mathfrak{n})^\top = \hat{F}^\mathfrak{n}$, it follows that
\be\label{Lbdd}
\|(\mathbf{L}^\mathfrak{n})^\top  - (U^\mathfrak{n})^\top  F^{\teq} \|_{A_{\boldsymbol{1}}} 
=\|\hat{F}^\mathfrak{n} - U^\mathfrak{n}(U^\mathfrak{n})^\top  F^\teq \|_{A_{\boldsymbol{1}}}
\leq \|\hat{F}^\mathfrak{n} -  F^\teq\|_{A_{\boldsymbol{1}}} = \|\varepsilon(\hat{f}_h^{\mathfrak{n}}- f^\teq_h) \|_{L^2(\Omega)}.
\ee
Any eigenvalue $\lambda_{D_{\mathbf{L}}}$ of the matrix $D_{\mathbf{L}}$ can be expressed in terms of the corresponding eigenvalue $\lambda$ of $(A_\EmAb)^{-1}A_{\boldsymbol{1}}$ as follows
\be
\lambda_{D_{\mathbf{L}}} = \frac{\lambda}{1+ \frac{1}{\Delta t} \lambda} = \frac{1}{\frac{1}{\lambda}+\frac{1}{\Delta t}}.
\ee
Therefore, according to \eqref{AAEAeig}, $\lambda_{D_{\mathbf{L}}}$ satisfies
\be\label{LDbdd}
0<\frac{1}{\chi_{\max} + \frac{1}{\Delta t} } \leq \lambda_{D_{\mathbf{L}}} \leq  \frac{1}{\chi_{\min} +\frac{1}{\Delta t} }<\frac{1}{\chi_{\min} }.
\ee
Together, \eqref{Lbdd}, \eqref{LDbdd}, and \Cref{exteig} imply that
\be\label{babbdd}
\bal
& \|B_{\mathbf{L}}^\top \|_{A_{\boldsymbol{1}}}^2 
\leq  \frac{1}{\chi_{\min}^2 } \|\varepsilon(\hat{f}_h^{\mathfrak{n}}- f^\teq_h) \|_{L^2(\Omega)}^2.
\eal
\ee
Let $S^\teq   (U^\teq)^\top   U^\mathfrak{n}=[l_1, \ldots, l_r] = l \in \mathbb{R}^{1 \times r}$ for scalars $l_i$ ($i=1,\ldots,r$).
Using \eqref{Feq}, \eqref{halfmatodeKb3} becomes 
\be\label{halfmatodeKb4}
\bal
\mathbf{L}^{\mathfrak{n}+1}
= &  E^{\teq} S^\teq   (U^\teq)^\top   U^\mathfrak{n} + \frac{1}{\Delta t} B_{\mathbf{L}} =   \left[ l_1 E^{\teq} + \frac{1}{\Delta t} {b}_1, \dots, l_r E^{\teq} + \frac{1}{\Delta t} {b}_r\right ], 
\eal
\ee
where $b_i\in\mathbb{R}^{n}$ ($i=1,\ldots,r$) are the column vectors of $B_{\mathbf{L}}$.
Combining \Cref{LstepQR} and the bound in \eqref{babbdd} gives
\be\label{babbddE}
1-\|(E^{\mathfrak{n}+1})^\top A_{\boldsymbol{1}} E^\teq \|^2  \leq \frac{\|B_{\mathbf{L}}^\top\|_{A_{\boldsymbol{1}}}^2}{\Delta t^2 \|l\|_{\infty}^2} \leq \frac{\|\varepsilon(\hat{f}_h^{\mathfrak{n}}- f^\teq_h) \|_{L^2(\Omega)}^2}{\Delta t^2 \|l\|_{\infty}^2 \chi_{\min}^2}.
\ee
For $l$, it holds
\be\label{inf_to_l2_est}
\|l\|_{\infty} 
= |S^\teq| \, \|(U^\mathfrak{n})^\top  U^\teq\|_{\infty} \geq \frac{|S^\teq| \, \|(U^\mathfrak{n})^\top  U^\teq\|}{\sqrt{r}} 
= \frac{|S^\teq|\|P_{U^\mathfrak{n}}U^\teq\|}{\sqrt{r}} \geq \frac{\beta|S^\teq|}{\sqrt{r}},
\ee
where the first inequality follows from the norm equivalence, and the last inequality follows from the assumption in \eqref{UUEEasumpt}.
Thus, if $\Delta t \geq \frac{\sqrt{r}}{\beta \delta \chi_{\min}}$, the estimate \eqref{basebdd} holds.  

The equality \eqref{UUEEasumpt1} follows from \eqref{babbddE} when $\hat{f}_h^\mathfrak{n}=f_h^\teq$.
\end{proof}

\subsection{Proof of \Cref{fsnprojU}}\label{pfsnprojU}
Similar to \Cref{LstepQR}, we prepare the following result.
\begin{lemma}\label{KstepQR}
Let $U^\teq \in \mathbb{R}^{n\times 1}$, $B_{\mathbf{K}}=[b_1 \ldots, b_r] \in \mathbb{R}^{m\times r}$, and $l=[l_1,\ldots,l_r]$ be a nonzero vector, where $b_i \in \mathbb{R}^{m\times 1}$ and $l_i \in \mathbb{R}$ for $i=1,\dots,r$.
Assume that the matrix
\be\label{Kn1eqb}
\mathbf{K}^{\mathfrak{n}+1}
= \left[ l_1 U^{\teq} + \frac{1}{\Delta t} {b}_1, \cdots, l_r U^{\teq} + \frac{1}{\Delta t} {b}_r\right ] \in \mathbb{R}^{m\times r},
\ee
has a decomposition $\mathbf{K}^{\mathfrak{n}+1} = U^{\mathfrak{n}+1} S_{\mathbf{K}}^{\mathfrak{n}+1}$ with $U^{\mathfrak{n}+1} = \left[ U^{\mathfrak{n}+1}_1, \cdots, U^{\mathfrak{n}+1}_r  \right ]$ satisfying \eqref{Fn1cond}. 
Then
\be\label{UeqUnbdd}
1- \|P_{U^{\mathfrak{n}+1}}U^\teq\|^2= 1- \|(U^{\mathfrak{n}+1})^\top  U^{\teq}\|^2 \leq \frac{\| B_{\mathbf{K}} \|_{\rm{F}}^2}{\Delta t^2 \|l\|_\infty^2}.
\ee
\end{lemma}

For any $E \in \mathbb{R}^{n\times r}$ satisfying $E^\top A_{\boldsymbol{1}} E=I_r$, 
because the term $E^\top A_\EmAb E$ will appear frequently, we introduce the symmetric matrix
\be\label{BEAE}
B = E^\top A_\EmAb E \in \mathbb{R}^{r\times r},
\ee
for which we have the following results.
\begin{lemma}
\label{lem:Beigen}
    Let $(\lambda_B, q_B)$ be an eigenpair of $B$ in \eqref{BEAE}.  Then 
    \be\label{lambdaB}
\chi_{\min} \leq \lambda_B = \frac{(E q_B)^\top A_\EmAb E q_B}{(E q_B)^\top A_{\boldsymbol{1}} E q_B}=\frac{\|(Eq_B)^\top\|_{A_\chi}^2}{\|(Eq_B)^\top\|^2_{A_{\boldsymbol{1}}}}\leq \chi_{\max}.
\ee 
\end{lemma}
\begin{proof}
If $(\lambda_B, q_B)$ is an eigenpair of $B$, then
\be\label{Beigen}
E^\top A_\EmAb E q_B = B q_{B} = \lambda_B q_B = \lambda_B E^\top A_{\boldsymbol{1}} Eq_B.
\ee
 Left-multiplying \eqref{Beigen} by $q_B^\top$ and applying \eqref{Ralquo} with $Z = (E q_B)^\top$ gives \eqref{lambdaB}.
\end{proof}

\begin{lemma}\label{projachieig}
Let $E \in \mathbb{R}^{n\times r}$ satisfy $E^\top A_{\boldsymbol{1}} E=I_r$, and recall the definition of $P_{E}^\chi$ from \eqref{PE-chi}.
Then for any $Z \in \mathbb{R}^{\ell \times n}$, $1\leq \ell \leq m$, 
\be\label{eqn:projachieig:0}
\|ZA_\chi P_{E}^{\chi}\|_{A_{\boldsymbol{1}}} \leq \sqrt{\frac{\chi_{\max}}{\chi_{\min}}} \|Z\|_{A_{\boldsymbol{1}}}.
\ee
\end{lemma}
\begin{proof}
Recall that $A_\EmAb$ is symmetric and positive definite, and thus can be decomposed as $A_\EmAb = C_\EmAb^\top C_\EmAb$ where $C_\EmAb$ is nonsingular.
Let $D_\chi = C_\chi E B^{-1}$, where $B$ is given in \eqref{BEAE}, and compute
\be\label{eqn:projachieig:1}
\|ZA_\chi P_{E}^{\chi}\|_{A_{\boldsymbol{1}}}= \|ZA_\chi E B^{-1}\|_{\rm{F}} \leq \|ZC_\chi^\top \|_{\rm{F}} \| C_\chi E B^{-1}\| = \|Z\|_{A_\chi}\|D_\chi\|.
\ee
Since $\|D_\chi\|^2$ is the largest eigenvalue of $D_\chi^\top D_\chi$ and 
\be
D_\chi^\top D_\chi = (B^{-1})^\top E^\top C_\chi^\top C_\chi EB^{-1} = B^{-1}BB^{-1} = B^{-1},
\ee
then \eqref{lambdaB} implies $\|D_{\chi}\| \leq \chi_{\min}^{-1/2}$, which along with \eqref{eqn:projachieig:1} and \eqref{Ralquo} yields \eqref{eqn:projachieig:0}.
\end{proof}
Next, we present the proof of \Cref{fsnprojU}.
\begin{proof}[Proof of \Cref{fsnprojU}]
The proof follows along the same lines as the proof of \Cref{fsnproj}. 
\eqref{fullDG1K-K} is equivalent to finding $\mathbf{K}^{\mathfrak{n}+1} \in \mathbb{R}^{m\times r}$ such that for any $\mathbf{K}_W\in \mathbb{R}^{m\times r}$,
\be\label{halfmatodeKa}
\left(D_t \mathbf{K}^{\mathfrak{n}+1}  (E^\mathfrak{n})^\top  A_{\boldsymbol{1}}, \mathbf{K}_W(E^\mathfrak{n})^\top  \right)_{\rm{F}}
=  \left(G(\mathbf{K}^{\mathfrak{n}+1}(E^\mathfrak{n})^\top ), \mathbf{K}_W(E^\mathfrak{n})^\top \right)_{\rm{F}}
\ee
Applying \eqref{Feqmat} and \Cref{FrobProp} to \eqref{halfmatodeKa} gives
\be
\left(D_t \mathbf{K}^{\mathfrak{n}+1}, \mathbf{K}_W \right)_{\rm{F}}
=  \left(F^\teq A_\EmAb E^\mathfrak{n}- \mathbf{K}^{\mathfrak{n}+1}(E^\mathfrak{n})^\top A_\EmAb E^\mathfrak{n}, \mathbf{K}_W\right)_{\rm{F}}.
\ee
Let $\mathbf{K}_W = \mathbf{K}'_W (B^\mathfrak{n})^{-1}$ for any $\mathbf{K}'_W \in \mathbf{R}^{m\times r}$, where $B^\mathfrak{n} = (E^\mathfrak{n})^\top A_\EmAb E^\mathfrak{n}$.  
Then
\be
\left((D_t \mathbf{K}^{\mathfrak{n}+1} (B^\mathfrak{n})^{-1}, \mathbf{K}'_W \right)_{\rm{F}}
=  \left(F^\teq A_\EmAb E^\mathfrak{n}(B^\mathfrak{n})^{-1}- \mathbf{K}^{\mathfrak{n}+1}, \mathbf{K}'_W\right)_{\rm{F}}.
\ee
Since $\mathbf{K}'_W$ in arbitrary, it follows that
\be\label{Keqdecomp}
\mathbf{K}^{\mathfrak{n}+1} =  F^\teq  A_\EmAb E^\mathfrak{n}(B^\mathfrak{n})^{-1} + \frac{1}{\Delta t} B_{\mathbf{K}},
\ee
where 
\begin{equation}\label{BBB+}
\begin{aligned}
B_\mathbf{K} 
    & :=[b_1, \ldots, b_r] 
    = \left( \mathbf{K}^\mathfrak{n}- F^\teq 
 A_\EmAb E^\mathfrak{n}(B^\mathfrak{n})^{-1}\right)D_{\mathbf{K}}\in \mathbb{R}^{m\times r},\\
    D_{\mathbf{K}} &:= (B^\mathfrak{n})^{-1} \left(I_r + \frac{(B^\mathfrak{n})^{-1}}{\Delta t} \right)^{-1}\in\mathbb{R}^{r\times r}.
\end{aligned}
\end{equation}
Since $\mathbf{K}^\mathfrak{n}=U^\mathfrak{n}S^\mathfrak{n}$, we can write
\be\label{KUSEeq}
\bal
& \mathbf{K}^\mathfrak{n}- F^\teq A_\EmAb E^\mathfrak{n}(B^\mathfrak{n})^{-1}  
= \left(F^\mathfrak{n} - F^\teq \right)  A_\EmAb E^\mathfrak{n}(B^\mathfrak{n})^{-1}.
\eal
\ee
By \eqref{KUSEeq}, \Cref{eigtoray}, \Cref{projachieig}, and \Cref{bilform}, 
\be\label{Kbdd}
\bal
 \|\mathbf{K}^\mathfrak{n}- F^\teq A_\EmAb E^\mathfrak{n}(B^\mathfrak{n})^{-1} \|_{\rm{F}}
&= \|\left(F^\mathfrak{n} - F^\teq \right)  A_\EmAb P_{E^{\mathfrak{n}}}^\chi\|_{A_{\boldsymbol{1}}} \\
&\leq \sqrt{\frac{\chi_{\max}}{\chi_{\min}}}\left\|F^\mathfrak{n} - F^\teq \right\|_{A_{\boldsymbol{1}}} 
 =  \sqrt{\frac{\chi_{\max}}{\chi_{\min}}} \|\varepsilon(\hat{f}_h^{\mathfrak{n}}- f^\teq_h) \|_{L^2(\Omega)}.
\eal
\ee  
Any eigenvalue $\lambda_{D_{\mathbf{K}}}$ of $D_{\mathbf{K}}$ satisfies
\be\label{LDbddK}
0<\frac{1}{\chi_{\max} + \frac{1}{\Delta t} } \leq \lambda_{D_{\mathbf{K}}} \leq  \frac{1}{\chi_{\min} +\frac{1}{\Delta t} }<\frac{1}{\chi_{\min} }.
\ee
Then, by \eqref{Kbdd}, \eqref{LDbddK}, and \Cref{exteig},
\be\label{BK_bound}
\|B_\mathbf{K}\|_{\rm{F}} \leq \frac{1}{\chi_{\min} }\sqrt{\frac{\chi_{\max}}{\chi_{\min}} }\|\varepsilon(\hat{f}_h^{\mathfrak{n}}- f^\teq_h) \|_{L^2(\Omega)}.
\ee
Let $S^\teq (E^\teq)^\top  A_\EmAb E^\mathfrak{n}(B^\mathfrak{n})^{-1} =[l_1, \ldots, l_r] \in \mathbb{R}^{1 \times r}$ for scalars $l_i$ ($i=1,\ldots,r$).
By \eqref{Feq} and \eqref{Keqdecomp},
\be\label{halfmatodeKb4+}
\bal
\mathbf{K}^{\mathfrak{n}+1}
= &  U^\teq S^\teq (E^\teq)^\top  A_\EmAb E^\mathfrak{n}(B^\mathfrak{n})^{-1} + \frac{1}{\Delta t} B_{\mathbf{K}} =   \left[ l_1 U^{\teq} + \frac{1}{\Delta t} {b}_1, \ldots, l_r U^{\teq} + \frac{1}{\Delta t} {b}_r\right ],
\eal
\ee
where $b_i \in \mathbb{R}^{m} \ (i=1,\ldots, r)$ are the columns of $B_{\mathbf{K}}$. 
By \Cref{KstepQR} and \eqref{BK_bound},
\be\label{babbddU}
1-\|(U^{\mathfrak{n}+1})^\top U^\teq \|^2  \leq \frac{\|B_{\mathbf{K}}\|_{\rm{F}}^2}{\Delta t^2 \|l\|^2_\infty} \leq \frac{\chi_{\max}}{\Delta t^2 \|l\|^2_\infty \chi_{\min}^3} \|\varepsilon(\hat{f}_h^{\mathfrak{n}}- f^\teq_h) \|_{L^2(\Omega)}^2.  
\ee
By the assumption \eqref{UUEEasumptU} and the fact that $ \| (B^\mathfrak{n})^{-1} (E^\mathfrak{n})^\top A_\EmAb E^\teq\| =\| (E^\teq)^\top A_\EmAb P_{E^{\mathfrak{n}}}^\chi \|_{A_{\boldsymbol{1}}} $, 
\be
\|l\|_\infty 
= |S^\teq| \, \| (B^\mathfrak{n})^{-1} (E^\mathfrak{n})^\top A_\EmAb E^\teq\|_{\infty} 
\geq \frac{|S^\teq| \, \| (E^\teq)^\top A_\EmAb P_{E^{\mathfrak{n}}}^\chi \|_{A_{\boldsymbol{1}}} }{\sqrt{r}} \geq \frac{\alpha|S^\teq|}{\sqrt{r}}.
\ee
Thus, if $\Delta t \geq \frac{\sqrt{r}\chi_{\max}^{1/2}}{\alpha \delta \chi_{\min}^{3/2}}$, estimate \eqref{basebddU} holds.

The equality \eqref{UUEEasumpt1U} follows from \eqref{babbddU} when $\hat{f}_h^\mathfrak{n}=f_h^\teq$.
\end{proof}

\bibliography{references}

\begin{thebibliography}{10}

\bibitem{abdi2007}
Herv{\'e} Abdi.
\newblock Singular value decomposition ({SVD}) and generalized singular value
  decomposition ({GSVD}).
\newblock {\em Encyclopedia of measurement and statistics}, 907:912, 2007.

\bibitem{adams_2001}
Marvin~L Adams.
\newblock Discontinuous finite element transport solutions in thick diffusive
  problems.
\newblock {\em Nuclear science and engineering}, 137(3):298--333, 2001.

\bibitem{ascher1997implicit}
Uri~M Ascher, Steven~J Ruuth, and Raymond~J Spiteri.
\newblock Implicit-explicit {R}unge-{K}utta methods for time-dependent partial
  differential equations.
\newblock {\em Applied Numerical Mathematics}, 25(2-3):151--167, 1997.

\bibitem{ayuso_etal_2011}
B.~{Ayuso}, J.~A. {Carrillo}, and {Shu} C.-W.
\newblock Discontinuous {G}alerkin methods for the one-dimensional
  {V}lasov--{P}oisson system.
\newblock {\em Kinetic and Related Models}, 4(4):955--989, 2011.

\bibitem{baumann2023}
Lena Baumann, Lukas Einkemmer, Christian Klingenberg, and Jonas Kusch.
\newblock Energy stable and conservative dynamical low-rank approximation for
  the {S}u-{O}lson problem.
\newblock {\em arXiv preprint arXiv:2307.07538}, 2023.

\bibitem{beck2000}
Michael~H Beck, Andreas J{\"a}ckle, Graham~A Worth, and H-D Meyer.
\newblock The multiconfiguration time-dependent hartree (mctdh) method: a
  highly efficient algorithm for propagating wavepackets.
\newblock {\em Physics reports}, 324(1):1--105, 2000.

\bibitem{ceruti2022}
Gianluca Ceruti and Christian Lubich.
\newblock An unconventional robust integrator for dynamical low-rank
  approximation.
\newblock {\em BIT Numerical Mathematics}, 62(1):23--44, 2022.

\bibitem{cheng_etal_2013}
Y.~{Cheng}, I.~M. {Gamba}, and P.~J. {Morrison}.
\newblock {Study of conservation and recurrence of Runge--Kutta discontinuous
  Galerkin schemes for Vlasov--Poisson systems}.
\newblock {\em Journal of Scientific Computing}, 56:319--349, 2013.

\bibitem{ding2021dynamical}
Zhiyan Ding, Lukas Einkemmer, and Qin Li.
\newblock Dynamical low-rank integrator for the linear {B}oltzmann equation:
  error analysis in the diffusion limit.
\newblock {\em SIAM Journal on Numerical Analysis}, 59(4):2254--2285, 2021.

\bibitem{dirac1930}
Paul~AM Dirac.
\newblock Note on exchange phenomena in the thomas atom.
\newblock In {\em Mathematical proceedings of the Cambridge philosophical
  society}, volume~26, pages 376--385. Cambridge University Press, 1930.

\bibitem{einkemmer2022asymptotic}
Lukas Einkemmer, Jingwei Hu, and Jonas Kusch.
\newblock Asymptotic--preserving and energy stable dynamical low-rank
  approximation.
\newblock {\em arXiv preprint arXiv:2212.12012}, 2022.

\bibitem{einkemmer2021asymptotic}
Lukas Einkemmer, Jingwei Hu, and Yubo Wang.
\newblock An asymptotic-preserving dynamical low-rank method for the
  multi-scale multi-dimensional linear transport equation.
\newblock {\em Journal of Computational Physics}, 439:110353, 2021.

\bibitem{einkemmer2021efficient}
Lukas Einkemmer, Jingwei Hu, and Lexing Ying.
\newblock An efficient dynamical low-rank algorithm for the {B}oltzmann-{B}gk
  equation close to the compressible viscous flow regime.
\newblock {\em SIAM Journal on Scientific Computing}, 43(5):B1057--B1080, 2021.

\bibitem{einkemmer2018}
Lukas Einkemmer and Christian Lubich.
\newblock A low-rank projector-splitting integrator for the {V}lasov--{P}oisson
  equation.
\newblock {\em SIAM Journal on Scientific Computing}, 40(5):B1330--B1360, 2018.

\bibitem{einkemmer2022}
Lukas Einkemmer, Alexander Ostermann, and Carmela Scalone.
\newblock A robust and conservative dynamical low-rank algorithm.
\newblock {\em Journal of Computational Physics}, 484:112060, 2023.

\bibitem{frenkel1934}
J~Frenkel.
\newblock Wave {M}echanics, {C}larendon, 1934.

\bibitem{grasedyck2013}
Lars Grasedyck, Daniel Kressner, and Christine Tobler.
\newblock A literature survey of low-rank tensor approximation techniques.
\newblock {\em GAMM-Mitteilungen}, 36(1):53--78, 2013.

\bibitem{guermondKanschat_2010}
J.-L. {Guermond} and G.~{Kanschat}.
\newblock Asymptotic analysis of upwind discontinuous {G}alerkin approximation
  of the radiative transport equation in the diffusive limit.
\newblock {\em SIAM J. Numer. Anal.}, 48:53--78, 2010.

\bibitem{hesthaven2007}
Jan~S Hesthaven and Tim Warburton.
\newblock {\em Nodal discontinuous {G}alerkin methods: algorithms, analysis,
  and applications}.
\newblock Springer Science \& Business Media, 2007.

\bibitem{jahnke2008}
Tobias Jahnke and Wilhelm Huisinga.
\newblock A dynamical low-rank approach to the chemical master equation.
\newblock {\em Bulletin of mathematical biology}, 70(8):2283--2302, 2008.

\bibitem{kieri2019projection}
Emil Kieri and Bart Vandereycken.
\newblock Projection methods for dynamical low-rank approximation of
  high-dimensional problems.
\newblock {\em Computational Methods in Applied Mathematics}, 19(1):73--92,
  2019.

\bibitem{lubichkoch2007DLRA}
Othmar Koch and Christian Lubich.
\newblock Dynamical low-rank approximation.
\newblock {\em SIAM Journal on Matrix Analysis and Applications},
  29(2):434--454, 2007.

\bibitem{kusch2022}
Jonas Kusch, Gianluca Ceruti, Lukas Einkemmer, and Martin Frank.
\newblock Dynamical low-rank approximation for {B}urgers’ equation with
  uncertainty.
\newblock {\em International Journal for Uncertainty Quantification}, 12(5),
  2022.

\bibitem{larsenMorel_1989}
E.~W. {Larsen} and J.~E. {Morel}.
\newblock {Asymptotic Solutions of Numerical Transport Problems in Optically
  Thick, Diffusive Regimes II}.
\newblock {\em Journal of Computational Physics}, 83:212--236, 1989.

\bibitem{lubich2005}
Christian Lubich.
\newblock On variational approximations in quantum molecular dynamics.
\newblock {\em Mathematics of computation}, 74(250):765--779, 2005.

\bibitem{lubich2014}
Christian Lubich and Ivan~V Oseledets.
\newblock A projector-splitting integrator for dynamical low-rank
  approximation.
\newblock {\em BIT Numerical Mathematics}, 54(1):171--188, 2014.

\bibitem{lubich2014projector}
Christian Lubich and Ivan~V Oseledets.
\newblock A projector-splitting integrator for dynamical low-rank
  approximation.
\newblock {\em BIT Numerical Mathematics}, 54(1):171--188, 2014.

\bibitem{mihalasMihalas_1999}
D.~{Mihalas} and B.~W. {Mihalas}.
\newblock {\em {Foundations of radiation hydrodynamics}}.
\newblock Dover (New York), 1999.

\bibitem{pareschi2005implicit}
Lorenzo Pareschi and Giovanni Russo.
\newblock Implicit--explicit {R}unge--{K}utta schemes and applications to
  hyperbolic systems with relaxation.
\newblock {\em Journal of Scientific Computing}, 25:129--155, 2005.

\bibitem{peng2021high}
Zhuogang Peng and Ryan~G McClarren.
\newblock A high-order/low-order ({HOLO}) algorithm for preserving conservation
  in time-dependent low-rank transport calculations.
\newblock {\em Journal of Computational Physics}, 447:110672, 2021.

\bibitem{peng2022sweep}
Zhuogang Peng and Ryan~G McClarren.
\newblock A sweep-based low-rank method for the discrete ordinate transport
  equation.
\newblock {\em Journal of Computational Physics}, 473:111748, 2023.

\bibitem{peng2020}
Zhuogang Peng, Ryan~G McClarren, and Martin Frank.
\newblock A low-rank method for two-dimensional time-dependent radiation
  transport calculations.
\newblock {\em Journal of Computational Physics}, 421:109735, 2020.

\bibitem{riviere2008}
B{\'e}atrice Rivi{\`e}re.
\newblock {\em Discontinuous {G}alerkin methods for solving elliptic and
  parabolic equations: theory and implementation}.
\newblock SIAM, 2008.

\bibitem{schotthofer2022}
Steffen Schotth{\"o}fer, Emanuele Zangrando, Jonas Kusch, Gianluca Ceruti, and
  Francesco Tudisco.
\newblock Low-rank lottery tickets: finding efficient low-rank neural networks
  via matrix differential equations.
\newblock {\em Advances in Neural Information Processing Systems},
  35:20051--20063, 2022.

\bibitem{shu2009}
Chi-Wang Shu.
\newblock Discontinuous {G}alerkin methods: general approach and stability.
\newblock {\em Numerical solutions of partial differential equations}, 201,
  2009.

\end{thebibliography}
\bibliographystyle{plain}

\end{document}